\documentclass{amsart}

\usepackage{type1cm}              
\usepackage{graphicx}              
\usepackage{multicol}               
\usepackage[bottom]{footmisc}
\usepackage{newtxtext}        
\usepackage{newtxmath}       
\usepackage{amsmath,amsfonts}
\usepackage{algorithm}
\usepackage{algorithmic}
\usepackage{booktabs}
\usepackage{multirow}
\usepackage{xcolor}

\makeindex             

\newcommand{\NF}{N^{\textsc{f}}}
\newcommand{\NL}{N^{\textsc{l}}}
\newcommand{\R}{\mathbb{R}}
\newcommand{\target}{x^\tau}
\newcommand{\xf}{x}
\newcommand{\xl}{y}
\newcommand{\vf}{v}
\newcommand{\vl}{w}

\newcommand{\HFF}{H^{\textsc{f}}}
\newcommand{\HFL}{H^{\textsc{l}}}
\newcommand{\HLF}{K^{\textsc{f}}}
\newcommand{\HLL}{K^{\textsc{l}}}
\newcommand{\barS}{s}

\newcommand{\Ct}{C_{\tau}}
\newcommand{\Cchar}{C_{s}}
\newcommand{\Crf}{C^{\textsc{f}}_{r}}
\newcommand{\Crl}{C^{\textsc{l}}_{r}}
\newcommand{\Cat}{C_{at}}
\newcommand{\Caf}{C^{\textsc{f}}_{al}}
\newcommand{\Cal}{C^{\textsc{l}}_{al}}

\newcommand{\xx}{{\bf x}}
\newcommand{\yy}{{\bf y}}
\newcommand{\rhotop}{{\varrho_{\textrm{top}}}}

\newtheorem{alg}{Algorithm}
\theoremstyle{remark}
\newtheorem{remark}{Remark}

\begin{document}

\title[Optimized strategies for crowd evacuation with multiple exits]{Optimized leaders strategies for crowd evacuation in unknown environments with multiple exits }

\author{Giacomo Albi, Federica Ferrarese, and Chiara Segala}

\thanks{
\\
Giacomo Albi:
Dipartimento di Informatica, Universit\`a di Verona, Verona, Italy, e-mail: giacomo.albi@univr.it.
\\
Federica Ferrarese:
Dipartimento di Matematica, Universit\`a di Trento, e-mail: federica.ferrarese@unitn.it.
\\
Chiara Segala:
Dipartimento di Matematica, Universit\`a di Trento, e-mail: chiara.segala-1@unitn.it.
\\ \ \\
All authors acknowledge the support of the Italian Ministry of Instruction, University and Research (MIUR) with funds coming from PRIN Project 2017 (No. 2017KKJP4X entitled ``Innovative numerical methods for evolutionary partial differential equations and applications'') and from RIBA 2019 (prot. RBVR199YFL entitled ``Geometric Evolution of Multi Agent Systems'').}

\begin{abstract}
In this chapter, we discuss the mathematical modeling of egressing pedestrians in an unknown environment with multiple exits. We investigate different control problems to enhance the evacuation time of a crowd of agents, by few informed individuals, named leaders. Leaders are not recognizable as such and consist of two groups: a set of unaware leaders moving selfishly toward a fixed target, whereas the rest is coordinated to improve the evacuation time introducing different performance measures. Follower-leader dynamics is initially described microscopically by an agent-based model, subsequently a mean-field type model is introduced to  approximate the large crowd of followers. The mesoscopic scale is efficiently solved by a class of numerical schemes based on direct simulation Monte-Carlo methods. Optimization of leader strategies is performed by a modified compass search method in the spirit of metaheuristic approaches. Finally, several virtual experiments are studied for various control settings and environments.
\end{abstract}

\maketitle

\section{Introduction}\label{sec:intro}
Control methodologies for crowd motion are of paramount importance in real-life applications for the design of safety measures and risk mitigation. The creation of virtual models of a large ensemble of pedestrians is a first step for reliable predictions, otherwise not easily reproducible with real-life experiments.

Pedestrians have been properly modeled by means of different agent-based dynamics such as lattice models \cite{cirillo2013PhysA, guo2012TRB}, social force models \cite{helbing2000N, parisi2005PhysA}, or cellular automata models \cite{abdelghany2014EJOR, wang2015PhysA}. A different level of description is obtained using mesoscopic models \cite{agnelli2015M3AS, albi2016SIAP, FTW18} where the quantities of study are densities of agents; at a larger scale macroscopic models \cite{carrillo2016M3AS,di2011hughes,colombo2012nonlocal} describe the evolution of moments such as mass and momentum. {Multiscale models} have been also considered, to account for situations where different scales coexist, we refer in particular to \cite{cristiani2011MMS, cristiani2014book}. Such a hierarchy of models is able to capture coherent global behaviors emerging from local interactions among pedestrians. These phenomena are strongly influenced by the social rules, the \emph{ rationality} of the crowd, and the knowledge of the surrounding environment. In the case of egressing pedestrians in an unknown environment with limited visibility we expect people to follow basically an instinctive behavior \cite{carrillo2016M3AS, cirillo2013PhysA, guo2012TRB,bailo2018pedestrian}, whereas a perfectly rational pedestrian will compute an optimal trajectory towards a specific target (the exit), forecasting exactly the behavior of other pedestrians \cite{abdelghany2014EJOR,lachapelle2011mean}.

In this manuscript, we focus on the evacuation problem in an unknown environment with multiple exits. We aim at influencing their behavior towards the desired target with minimal intervention. Starting from the seminal work \cite{albi2016SIAP} we consider a bottom-up approach where few informed agents are acting minimizing verbal directives to individuals and preserving as much as possible their natural behavior. This approach is expected to be efficient in situations where direct communication is impossible, for example in the case of very large groups, emergencies, violent crowds reluctant to follow directions; or in panic situations where rational behavior is overtaken by instinctive decisions.
Furthermore, we consider few additional agents, who are informed about the position of some exit and acting as unaware leaders. Hence, their dynamics will influence the global behavior of the crowd, introducing inertia that may constitute an additional difficulty in the optimization problem, for example increasing congestions next to the exits or increasing the level of uncertainty.

The control problem associated to the evacuation of a crowd falls in the larger research field aimed at investigating the control of \textit{self-organizing agents}. From the mathematical view point, this type of problem is challenging due to the presence of non-local interaction terms and their high dimensionality. Control of alignment-type dynamics, such as the Cucker-Smale model \cite{cucker2007IEEE}, have risen a lot of interest in the mathematical community, where several strategies have been explored to enforce the emergence of consensus, see for example \cite{ BBCK18, BFK15, albi2017mean,HePaSt15}. At the same time, to cope with the high dimensionality of such optimal control problems, reduced approaches have been explored \cite{CFPT13,FS14,fornasier2014mean,bongini2017mean}, promoting sparsity of the control acting only on few agents. In biological models, it has been shown that a small percentage of individuals can influence the whole group towards a desired target, see \cite{couzin2005N}.  Similarly, leaders in crowd can act as control signals to enforce alignment towards a desired direction as recognizable leaders \cite{ABFHKPPS, albi2014Phil, borzi2015M3ASa, during2009PRSA, motsch2011JSP}, or moving undercover \cite{albi2016SIAP, couzin2005N, duan2014SR, han2006JSSC, han2013PLOSONE}, or even in a repulsive way \cite{burger2020instantaneous}. These strategies heavily rely on the power of the \emph{social influence} (or herding effect), namely the natural tendency of people to follow other mates in situations of emergency or doubt. 

Alternative control methodologies consist in optimal design of the surrounding enviroment such as obstacles \cite{cristiani2019AMM,cristiani2017AMM,albi2020Crowd}, or evacuation signage
\cite{zhang2017optimal,xie2012experimental}, or exit locations \cite{wang2015PhysA}.


The manuscript is organized as follows in Section \ref{sec:models} we introduce the mathematical framework for the microscopic dynamics of leader-follower type and we formulate different scenarios for the optimal control problem to be solved. In particular, we will distinguish between minimum time of evacuation, total mass evacuated, and optimal mass splitting among the multiple exits. In this work the word mass denotes the total amount of pedestrians.
Section \ref{sec:meanfield} is devoted to the description of the mesoscopic scale, first we introduce the mean-field type model, second we sketch an efficient Monte Carlo algorithm for its simulation. In Section \ref{sec:numerics} we focus on the numerical realization of the optimized strategies. We start introducing the algorithmic procedure used for the solution of the large-scale optimization problem, and we compare microscopic and mean-field dynamics in several scenarios and with different target functionals.
Finally in Section \ref{sec:conclusions} we outline possible extensions and further perspectives.

\section{Control of pedestrian dynamics through leaders}\label{sec:models}
In this section, we focus first on the mathematical description of pedestrian dynamics in complex environments. We consider an ensemble of agents,  {\em followers}, in an unknown environment trying to reach exit locations, at the same time the crowd population includes few informed agents, {\em leaders}, acting as controllers but not distinguishable from followers. In particular, we account for a mixed approach where leaders are either aware of their role, then responding to an optimal force as the result of an offline optimization procedure, {\em optimized leaders}, or unaware of their role and moving with a greedy strategy towards a target exit position, {\em selfish leaders}. The main mechanisms ruling the behaviors among the followers are isotropic interactions with other agents based on metrical short-range repulsion, induced by social distancing and collisional avoidance, and topological long-range alignment dynamics. Leaders instead consider only short-range repulsion. Additionally, for followers, we account self-driving forces describing the exploration phase, preferential direction, and desired speed.  The overall dynamics will be influenced by the surrounding environment when the exits are visible or close to obstacles.

In the following sections, we describe first the microscopic dynamics of the follower-leader system and later different control tasks for different applications.
\subsection{Microscopic model with leaders and multiple exits}\label{sec:micromodel}
Following the approach proposed in \cite{albi2016SIAP,albi2020Crowd} we model leaders by a first-order model and followers by a second-order one, where both positions and velocities are state variables. We denote by $d$ the dimension of the space in which the motion takes place (typically $d=2$), by $\NF$ the number of followers and by $\NL \ll \NF$ the number of leaders. We also denote by $\Omega\equiv\R^d$ the walking area, and we identify the different exits by $\target_e\in\Omega$ with $e=1,\ldots,N_e$.

To define each target's visibility area, we consider the set $\Sigma_e$, with $\target_e\in\Sigma_e\subset\Omega$, and we assume that the target is completely visible from any point belonging to $\Sigma_e$ and completely invisible from any point belonging to $\Omega\backslash\Sigma_e$, namely we also assume that visibility areas are disjoint sets, i.e $\Sigma_{e_i}\cap\Sigma_{e_j}=\varnothing$ for all $e_i,e_j\in\{1,..,N_e\}$ .

For every $i=1,\ldots,\NF$, let $(\xf_i(t),\vf_i(t))\in\R^{2d}$ denote position and velocity of the agents belonging to the population of followers at time $t\geq 0$ and, for every  $k=1,\ldots,\NL$, let $(\xl_k(t),\vl_k(t))\in\R^{2d}$  denote position and velocity of the agents among the population of leaders at time $t \geq 0$. Let us also define $\mathbf{\xf}:=(\xf_1,\ldots,\xf_{\NF})$ and $\mathbf{\xl}:=(\xl_1,\ldots,\xl_{\NL})$.

The microscopic dynamics described by the two populations is given by the following set of ODEs for $i = 1, \dots, \NF$ and $k = 1, \ldots, \NL$,
\begin{equation}\label{eq:micro}
\left\{
\begin{array}{l}
\dot{x}_i = \vf_i,\\ [1.5mm] 
\dot{v}_i = S(\xf_i,\vf_i) + \sum_{j=1}^{\NF} m^F_j \HFF(\xf_i,\vf_i,\xf_j,\vf_j;\mathbf{\xf},\mathbf{\xl}) +\sum_{\ell=1}^{\NL} m^L_\ell \HFL(\xf_i,\vf_i,\xl_\ell,\vl_\ell;\mathbf{\xf},\mathbf{\xl}),\\ [1.5mm]
\dot{y}_{k} = w_{k} = \sum_{j=1}^{\NF} m^F_j \HLF(\xl_k,\xf_j) +  \sum_{\ell=1}^{\NL}m^L_\ell \HLL(\xl_k,\xl_\ell) + \xi_k u^{\textrm{opt}}_{k}+(1-\xi_k)u^{\textrm{self}}_{k},\\[1.5mm]
\end{array}
\right.
\end{equation}
with initial data for followers $(x_i(0),v_i(0)) = (x^0_i,v^0_i)$  and leaders $(y_k(0),w_{k}(0))=(y_k^0,w_{k}^0)$. The quantities $m_i^F,m_k^L$ weight the interaction of followers and leaders, in what follows we will assume that $m_1^F=\ldots=m_{\NF}^F=m_1^L=\dots=m_{\NL}^L$ and the following mass constraint holds
\begin{align}\label{eq:scaled}
m_i^F =\frac{\rho^F}{\NF}, \qquad m_k^L =\frac{\rho^L}{\NL},\qquad \rho^F+\rho^L = 1,
\end{align}
for $\rho^F,\rho^L$ positive quantities.
 
\begin{enumerate}
\item $S$ is a self-propulsion term, given by the relaxation toward a random direction or the relaxation toward a unit vector pointing to the target (the choice depends on the position), plus a term which translates the tendency to reach a given characteristic speed $\barS \geq 0$ (modulus of the velocity), i.e.,
\begin{align} \label{eq:Ad}
S(x,v) := \Cchar(\barS^2-|v|^2)v+\sum_{e=1}^{N_e}\psi_e(x)\Ct\left(\frac{\target_e - x}{|\target_e - x|} - v\right),
\end{align}
where $\psi_e:\mathbb{R}^d \rightarrow [0,1]$ is the characteristic function of $\Sigma_e$, and $\Ct$, $\Cchar$ are positive constants.
\item
The interactions follower-follower and follower-leader account a repulsion and an alignment component, as follows
\begin{equation}
\begin{aligned}\label{eq:HFF}
&\HFF(x,v,x',v';\mathbf \xf,\mathbf \xl) :=
-\Crf R_{\gamma,r}(x,x')(x'-x)+(1-\psi(x)){\Caf}A(x,x';\mathbf\xf,\mathbf\xl)\left(v'-v\right) ,\\
&\HFL(x,v,y,w;\mathbf \xf,\mathbf \xl) := -\Crl R_{\gamma,r}(x,y)(y-x) +
(1-\psi(x)){\Cal}A(x,y;\mathbf\xf,\mathbf\xl)\left(w-v\right),
\end{aligned}
\end{equation}
for given positive constants $\Crf, \Caf, \Cal, \Cat, r,\gamma$, and where $1-\psi(x)$ is the characteristic function of the unknown environment $\Omega\backslash\cup_e\Sigma_e$,  such that  
\[
\psi(x):=\sum_{e=1}^{N_e}\psi_e(x).
\]

 The first term on the right hand side of \eqref{eq:HFF} represents the metrical repulsion force, where the intensity is modulated by the function $R_{\gamma,r}$ defined as
\begin{align}\label{eq:rep}
R_{\gamma,r}(x,y) & = \begin{cases}
\frac{e^{-|y-x|^\gamma}}{|y-x|} & \text{ if } y\in B_r(x)\backslash\{x\}, \\
0 & \text{ otherwise,}
\end{cases}
\end{align}
 where $B_r(x)$ is the ball of radius $r>0$ centered at $x\in\Omega$. 
The second term accounts for the (topological) alignment force, which vanishes inside the visibility regions, and where 
\begin{align}\label{eq:altop}
A(x,y;\mathbf\xf,\mathbf\xl) :=\chi_{\mathcal B_\mathcal{N}(x;\mathbf \xf,\mathbf \xl)}(y),
\end{align}
and by $\mathcal B_\mathcal N(x;\mathbf{\xf},\mathbf{\xl})$ the \emph{minimal} ball centered at $x$ encompassing at least $\mathcal N$ agents.

\item
The interactions leader-follower and leader-leader reduce to a mere (metrical) repulsion, i.e., $\HLF = \HLL = -\Crl R_{\zeta,r}$, where $\Crl>0$ and $\zeta>0$ are in general different from $\Crf$ and $\gamma$, respectively. 

\item
$u^\textrm{opt}_k,u^\textrm{self}_k:\R^+\to\R^{d\NL}$ characterize the strategies of the leaders and  are chosen in a set of admissible control functions. The parameter $\xi_k\in\left\{0,1\right\}$ identifies for $\xi_k=1$ leaders aware of their role, whose movements are the result of an optimization process, and alternatively for $\xi_k=0$ leaders moving ``selfishly'' towards a specific exit.
A specific description of leaders' strategy will be discussed in Section \ref{sec:numerics}.
Hence we account for situations where a small part of the mass is informed about exit positions, but policymakers have no control over them.

\end{enumerate}

\begin{remark}~
	\begin{itemize}
	\item Differently from the model proposed in \cite{albi2016SIAP,albi2020Crowd} the dynamics do not include random effects. However, we consider this uncertainty by assuming that the initial velocity directions of followers are distributed according to a prescribed density $v^0_i\sim p_v(\R^d)$, for example, a uniform distribution over the unitary sphere $\mathbb S_{d-1}$ .
\item The choice $\Caf=\Cal$ leads to $\HFF\equiv \HFL$ and, therefore, the leaders are not recognized by the followers as special. This feature opens a wide range of new applications, including the control of crowds not prone to follow authority's directives.
\item The pedestrian microscopic model \eqref{eq:micro}  allows agent movements in space without any constriction. However, in real applications, dynamics are constrained by walls or other kinds of obstacles. There are several ways of dealing with this feature in agent-based mode and we refer to \cite[Sect.\ 2]{cristiani2017AMM} for a review of obstacles handling techniques such as repulsive obstacle, rational turnaround, velocity cut-off. The choice for obstacle handling will be discussed in Section \ref{sec:numerics}.
%
			
\end{itemize}
\end{remark}

\subsection{Control framework for pedestrian dynamics}\label{sec:opt2}
In order to define the strategies of optimized leaders, we formulate an optimal control problem to exploit the tendency of people to follow group mates in situations of emergency or doubt.
The choice of a proper functional to be minimized constitutes a modeling difficulty, and it is typically a trade-off between a realistic task and a viable realization of its minimization. In general we 
will set up the following constrained optimal control  problem
\begin{equation}
\begin{aligned}\label{eq:opt}
\min_{{\bf u}^{\textrm {opt}}(\cdot)\in U_{adm} }  \mathcal J({\bf u}^{\textrm {opt}}),
\cr \textrm{s.t.} \qquad \eqref{eq:micro},
\end{aligned}
\end{equation}
where ${\bf u}^{\textrm {opt}} = (u_k^{opt}(\cdot))$ is the control vector associated to the optimized leaders, given a set of admissible controls $U_{adm}$. In what follows we will specify different functionals for different type of applications.
For later convenience  we introduced the empirical distributions defined as follows
\begin{align}\label{eq:empirical}
f^{\NF}(\cdot,x,v) = \sum_{i=1}^{\NF} m_i^{F}\delta(x-x_i(\cdot))\delta(v-v_i(\cdot)),\\
g^{\NL}(\cdot,x,v) = \sum_{j=1}^{\NL} m_j^{L}\delta(x-y_j(\cdot))\delta(v-w_j(\cdot)).
\end{align}
  
\begin{itemize}
	\item {\em Evacuation time.}  In a situation where egressing pedestrians are in an unknown enviroment the most natural functional is the { evacuation time}, that we may define as follows
\begin{equation}\label{eq:minevactime}
\mathcal J(\textbf{x},{\bf y},\textbf{u}^{\textrm {opt}}) =\left\lbrace t > 0 \ | \ (x_i(t),y_j(t)) \notin \Omega\ \forall i = 1, \dots N^F,\forall j=1,\ldots,\NL \right\rbrace,
\end{equation}
where we explicit the dependency on the states vector of follower positions ${\bf x}\in \R^{d\NF}$.
This cost functional is extremely irregular, therefore the search of minima is particularly difficult, additionally the evacuation of the total mass in some situations
can not be completely reached.

\item {\em Total mass with multiple exits.} Instead of minimizing the total evacuation, we fix a final time $T>0$ and we aim to minimize the total mass inside the computational domain $\Omega\setminus\cup_e\Sigma_e$, which coincides with maximizing the mass inside the visibility areas. The functional reads
\begin{equation} \label{eq:test2intro} 
\mathcal J ({\bf x},{\bf y},{\bf u}^{\textrm {opt}}) = \int_{\mathbb{R}^d} \int_{\Omega\setminus\cup_e\Sigma_e} (f^{\NF}(T,x,v) + g^{\NL}(T,x,v))dx dv.
\end{equation}
\item {\em Optimal mass splitting over multiple exits.} In complex environments, it may happen that total mass does not distribute in an optimal way between the target exits.  This may lead to problems of heavy congestions and overcrowding around the exits that, in real-life situations, can cause injuries due to overcompression and suffocation. Hence we ask to distribute the total evacuated mass at final time $T$ among the exits according to a given desired distribution. To this end we set
\begin{equation}\label{eq:test3intro}
\mathcal J ({\bf x},{\bf y},{\bf u}^{\textrm{opt}}) = \sum_{e=1}^{N_e} \left|\mathcal M^F_e(T)-\mathcal M^{\textrm{ des}}_e \right| ^2,  
\end{equation}
where $\mathcal M_e^{\textrm{des}}$ is the desired mass to be reached in the visibility area $\Sigma_e$ and $\mathcal M^F_e(T)$ is the total mass of followers and leaders who reached exit $\target_e$ up to final time $T$.
\end{itemize}

\section{Mean-field approximation of follower-leader system}\label{sec:meanfield}
Mean-field scale limit for large number of interacting individuals has been investigated in several directions for single and multiple population dynamics, see for example \cite{carrillo2010particle,di2013measure}, and it is a fundamental step to tame the curse of dimensionality arising for coupled systems of ODEs.

In the current setting, we want to give a statistical description of the followers-leaders dynamics considering a continuous density for followers and maintaining leaders microscopic. Hence, we introduce  the non-negative distribution function of followers  $f=f(t,x,v)$ with $x\in\R^d, v\in\R^d$ at time $t\geq0$,  the meso-micro system corresponding to \eqref{eq:micro} reads as follows 
\begin{equation}\label{eq:MFmodelStrong}
\begin{aligned}
&\partial_t f + v \cdot \nabla_x f = - \nabla_v \cdot \left(f\left(S(x,v)+\mathcal{H}^F[f,g^{\NL}]+\mathcal{H}^L[f,g^{\NL}]\right)\right),\\
\displaystyle
&\dot{\xl}_k = \vl_k = \int_{\R^{2d}} \HLF(\xl_k,x) f(t,x,v) \ dx\ dv + \sum_{\ell = 1}^{\NL}m^L_\ell \HLL(\xl_k,\xl_{\ell})\cr
 &\qquad\qquad\qquad\qquad+ \xi_k u^{\textrm{opt}}_{k}+(1-\xi_k)u^{\textrm{self}}_{k},
\end{aligned}
\end{equation}
where the followers dynamics is described by a kinetic equation of Vlasov-type, and where we use the corresponding empirical distribution for leaders $g^{\NL}$.
Furthermore we assume that the follower and leader densities are such that their number densities are
\[
\varrho^F =  \int_{\R^{2d}} f(t,x,v)\ dx\ dv,\qquad \varrho^L =  \int_{\R^{2d}} g^{\NL}(t,x,v)\ dx\ dv.
\]
We observe that the terms $S(\cdot),\HLF(\cdot)$ and $\HLL(\cdot)$  are defined respectively as in the microscopic setting, whereas
the non-local operators $\mathcal{H}^F,\mathcal{H}^L$ correspond to the following integrals
\begin{align}
&\mathcal{H}^F[f,g^{\NL}](t,x,v) = -\Crf\int_{\R^{d}}\int_{B_r(x)}R_{\gamma,r}(x,x')(x'-x)f(t,x',v')\ dx' dv'\cr
&\qquad +\Caf(1-\psi(x))\int_{\R^{d}}\int_{\mathcal B_{r_*}(t,x)}(v'-v)f(t,x',v')\ dx' dv',
\\
&\mathcal{H}^L[f,g^{\NL}](t,x,v) = -\Crl\int_{\R^{d}}\int_{B_r(x)}R_{\gamma,r}(x,x')(x'-x)g^{\NL}(t,x',v')\ dx' dv'\cr
&\qquad +\Cal(1-\psi(x))\int_{\R^{d}}\int_{\mathcal B_{r_*}(t,x)}(v'-v)g^{\NL}(t,x',v')\ dx' dv',
\end{align}
where the first term corresponds to the metrical repulsion as in \eqref{eq:rep}, and the second part accounts the topological ball $\mathcal B_{r_*}(t,x)\equiv\mathcal B_{r^*}(t,x;f,g^{\NL})$ whose radius is defined for a fixed $t\geq 0$ by the following variational problem
\begin{align}\label{eq:topradius}
r^*(t,x) = \arg\min_{\alpha>0}\left\{ \int_{\R^d}\int_{B_\alpha(x)} \left(f(t,x,v) +g^{\NL}(t,x,v) \right)\ dx \ dv\geq \mathcal \rhotop\right\},
\end{align}
where $\rhotop>0$ is the target topological mass.

\begin{remark}~
\begin{itemize}
	\item Rigorous derivation of the mean-field limit \eqref{eq:MFmodelStrong} from \eqref{eq:micro} is a challenging task due to the strong irregularities induced by the behavior of topological-type interactions. We refer to \cite{haskovec13topological} for possible regularization in the case of Cucker-Smale type dynamics, and to \cite{degond17topological,degond19topological} for alignment driven by jump-type processes.  
	\item Alternative derivation of mesoscopic models in presence of diffusion has been obtained in \cite{albi2016SIAP}, where the authors derived a Fokker-Planck equation of the original microscopic system via quasi-invariant scaling of binary Boltzmann interactions. This technique, analogous to the so-called grazing collision limit in plasma physics, has been thoroughly studied in \cite{Vil02} and allows to pass from a Boltzmann description to the mean-field limit, see for example \cite{PT:13}.
	\item For optimal control of large interacting agent systems, the derivation of a mean-field approximation involves the convergence of minimizers from microscopic to mesoscopic scale. This problem has been addressed from different directions, and we refer to \cite{FS14,bongini2017mean}. 
	
 \end{itemize}	
\end{remark}
\begin{remark}~
In order to obtain a closed hydrodynamic system for \eqref{eq:MFmodelStrong} a standard assumption is to assume the velocity distribution to be mono-kinetic, i.e.\ $f(t,x,v)=\rho(t,x)\delta(v-V(t,x))$, and the fluctuations to be negligible. Hence, computing the moments of \eqref{eq:MFmodelStrong} leads to the following macroscopic system for the density $\rho$ and the bulk velocity $V$,
\begin{equation}
\left\{ 
\begin{array}{l}
\partial_t \rho + \nabla_x\cdot (\rho V) = 0,\\ [2mm]
\partial_t (\rho V)+ \nabla_x\cdot (\rho V\otimes V) = \mathcal{G}_m\left[\rho,\rho^{\textsc{l}},V,V^{\textsc{l}}\right]\rho,\\ [1mm]
\displaystyle
\dot{\xl}_k = \vl_k = \int_{\R^{d}} \HLF(\xl_k,x) \rho(t,x) \ dx+ \sum_{\ell = 1}^{\NL} \HLL(\xl_k,\xl_{\ell}) \\\qquad\qquad\qquad\qquad+ \xi_k u^{\textrm{opt}}_{k}+(1-\xi_k)u^{\textrm{self}}_{k},
\end{array}
\right.
\end{equation}
where $\rho^{\textsc{l}}(x,t), V^{\textsc{l}}(x,t)$ represent the leaders macroscopic density and bulk velocity, respectively, and $\mathcal{G}_m$ the macroscopic interaction operator associated to the followers, we refer to \cite{albi2013AML,carrillo2010particle} for further details. 
\end{remark}

\subsection{MFMC algorithms }
For the numerical solution of the mean-field followers dynamics in \eqref{eq:MFmodelStrong} we employ
mean-field Monte-Carlo methods (MFMCs) generalizing the approaches proposed in \cite{albi2013MMS,PT:13}. These methods fall in
the class of fast algorithms developed for interacting particle systems such as direct
simulation Monte-Carlo methods (DSMCs), and they are strictly related to more recent class of algorithms named Random
Batch Methods (RBMs) \cite{jin2020random}. 

In order to approximate the evolution of the followers density, first we sample $N^F_s$ particles from the initial distribution $f^0(x,v)$ in the phase space, i.e. $\{(x_i^0,v_i^0)\}_{i=1}^{N^F_s} $.
Furthermore we consider a subsample of $M$ particles, $j_1,\ldots,j_M$  uniformly without repetition such that $1\leq M\leq N_s^F$.
In order to approximate the non-local terms $\mathcal H^F,\mathcal H^L$ we evaluate the interactions with a subsample of size $M$ at every time step. Hence we define the discretization step as
\begin{align}\label{eq:vnp1}
v_i^{n+1} &= v_i^{n}+\Delta t S(x_i^n,v_i^n)-\Delta t\left[\hat{R}^{F,n}_i (\hat{X}^{n}_i-x^n_i)+\hat{R}^{L,n}_i (\hat{Y}^{n}_i-x^n_i)\right]
\cr
&\qquad+\Delta t (1-\psi(x_i^n))\left[\hat{A}^{F,n}_i (\hat{V}^{n}_i-v^n_i)+\hat{A}^{L,n}_i (\hat{W}^{n}_i-v^n_i)\right],
\end{align}
where we defined the following auxiliary variables for the repulsion term \eqref{eq:rep},
\begin{equation}\label{eq:RX}
\begin{split}
& \hat{R}^{F,n}_i = \frac{\Crf\varrho^F}{M} \sum_{k=1}^M R_{\gamma,r}(x^n_i,x^n_{j_k}), 
\quad \hat{X}^{n}_i = \frac{\Crf\varrho^F}{M} \sum_{k=1}^M\frac{R_{\gamma,r}(x^n_i,x^n_{j_k})}{\hat{R}_i^{F,n}}x^n_{j_k},\\ 
&\hat{R}^{L,n}_i = \frac{\Crl\varrho^L}{\NL} \sum_{\ell=1}^{\NL} R_{\gamma,r}(x^n_i,y^n_{\ell}), 
\quad \hat{Y}^{n}_i = \frac{\Crl\varrho^L}{\NL} \sum_{k=1}^{\NL}\frac{R_{\gamma,r}(x^n_i,y^n_{\ell})}{\hat{R}_i^{L,n}}y^n_{\ell}.
\end{split}
\end{equation}
For the topological alignment we have
\begin{equation}\label{eq:PV}
\begin{split}
& \hat{A}^{F,n}_i = \frac{\Caf\varrho^F}{M} \sum_{k=1}^M \chi_{\mathcal{B}_{r^*_M}(x_i;\xx,\yy)} (x_{j_k}), 
\quad \hat{V}^{n}_i = \frac{\Caf\varrho^F}{M} \sum_{k=1}^M\frac{\chi_{\mathcal{B}_{r^*_M}(x_i;\xx,\yy)} (x_{j_k})}{\hat{A}_i^{F,n}}v^n_{j_k},\\ 
&\hat{A}^{L,n}_i = \frac{\Cal\varrho^L}{\NL} \sum_{\ell=1}^{\NL}\chi_{\mathcal{B}_{r^*_M}(x_i;\xx,\yy)}(x_{j_k}), 
\quad \hat{W}^{n}_i = \frac{\Cal\varrho^L}{\NL} \sum_{k=1}^{\NL}\frac{\chi_{\mathcal{B}_{r^*_M}(x_i;\xx,\yy)} (x_{j_k})}{\hat{A}_i^{L,n}}w^n_{\ell},
\end{split}
\end{equation}
where, the topological ball $\mathcal B_{r^*_M}(x)$ is the topological ball defined over the subsample of $M$ agents, with radius such that
\begin{align}\label{eq:topradius2}
r_M^*(t,x_i) = \arg\min_{\alpha>0}\left\{ \frac{\varrho^F}{M}\sum_{k=1}^M\chi_{{B}_{\alpha}(x_i)}(x_{j_k})+\frac{\varrho^L}{\NL}\sum_{\ell=1}^{\NL}\chi_{{B}_{\alpha}(x_i)}(y_{\ell})\geq \mathcal \rhotop\right\}.
\end{align}
From the above considerations we obtain the following Algorithm in the time interval $[0,T]$. 

\begin{alg}[MFMC follower-leader]~ \label{alg_MC}
	\begin{enumerate}
		\item[\texttt 1.] Given $N^F_s$ samples $v_i^0$, with $i=1,\ldots,N^F_s$ computed from the initial distribution $f(x,v)$ and $M\leq N^F_s$;
		\item[\texttt 2.] \texttt{for} $n=0$ \texttt{to} $n_{\textrm{tot}}$
		\begin{enumerate}
			\item \texttt{for} $i=1$ \texttt{to} $N^F_s$
			\begin{enumerate}
				\item sample $M$ particles $j_1,\ldots,j_M$  uniformly without repetition among all particles;
				\item compute the quantities $\hat{R}_i^{L,n},\hat{R}_i^{F,n},\hat{X}_i^n$ and $\hat Y^{n}_i$ from \eqref{eq:RX};
				\item compute the quantities $\hat{A}_i^{L,n},\hat{A}_i^{F,n},\hat{V}_i^n$ and $\hat W^{n}_i$ from \eqref{eq:PV};
				\item compute the velocity change $v_i^{n+1}$ according to \eqref{eq:vnp1};
				\item compute the position change
				\[
				x^{n+1}_i =x_i^n+\Delta t v_i^{n+1}.
				\]
			\end{enumerate}
		
		\end{enumerate}
		\texttt{end for}
	\end{enumerate}
	\texttt{end for}
\end{alg}
\begin{remark}~
	\begin{itemize}
\item By using this Monte Carlo algorithm we can reduce the computational cost due to the computation of the interaction term from the original $\mathcal{O}\left({N_s^F}^2\right)$ to $\mathcal{O}\left(M N_s^F\right)$. For $M=N_s^F$ we obtain the explicit Euler scheme for the original $N_s^F$ particle system. 
  
\end{itemize}
\end{remark}

\section{Numerical optimization of leaders strategies} \label{sec:numerics}
In this section we focus on the numerical realization of the general optimal control problem of type
\begin{equation}
\begin{aligned}\label{eq:opt2}
\min_{{\bf u}^{\textrm {opt}}(\cdot)\in U_{adm} }  \mathcal J({\bf u}^{\textrm {opt}}),
\end{aligned}
\end{equation}
constrained to the evolution of microscopic \eqref{eq:micro} or mean-field system \eqref{eq:MFmodelStrong}. We observe that the minimization task for evacuation time or total mass can be extremly difficult, due to the strong irregularity and the presence of many local minima.


In order to optimize $\eqref{eq:opt2}$ we propose instead an alternative suboptimal, but computationally efficient strategy, named modified Compass Search (CS). This method falls in the class of metaheuristic algorithms, it ensures the convergences towards local minima, without requiring any regularity of the cost functional \cite{audet2014MPC}. 

We use the CS method in order to optimize the trajectory of the {\em aware } leaders. 
The idea is to start from an initial guess $u_k^{\textrm{opt},(0)}$ which produces an admissible trajectory toward a target exit, for example as follows
\begin{equation}\label{eq:gototarget_beta}
u^{\rm opt, (0)}_k(t) = \beta \frac{\Xi_k(t)-y_k(t)}{\Vert  \Xi_k(t)-y_k(t) \Vert}+(1-\beta)(m_F(t) - y_k(t)),
\end{equation}
where $\Xi_k(t)$ is the target position at time $t$, depending on the environment and such that  $\Xi_k(t)=x^\tau_e$ for $t>t_*$ . The parameter $\beta\in[0,1]$ measures the tendency of leaders to move toward the target $\Xi_k(t)$ or staying close to followers center of mass $m_F(t)$.

We will refer  to \eqref{eq:gototarget_beta} as ``go-to-target'' strategy. Then CS method iteratively modifies the current best control strategy found so far computing small random piecewise constant variation of points on the trajectories. 
Then, if the cost functional decreases, the variation is kept, otherwise it is discarded. We consider piecewise constant trajectories, introducing suitable switching times for the leaders controls.

We summurize this procedure in the following algorithm.
\begin{alg}[Modified Compass Search]~ \label{alg_CS}
	\begin{enumerate}
		\item[\texttt 1.] Select a discrete set of sample times $S_M=\{t_1,t_2,\ldots,t_M\}$, the parameters $j =0$, $j_{\texttt {max}}$ and $J_E$.
		\item[\texttt 2.] Select an initial strategy $u^*$  piecewise constant over the set $S_M,$
		e.g. constant direction and velocity speed towards a fixed target $\Xi_k(t)$, $k=1,\dots,N^L$, see Equation \eqref{eq:gototarget_beta}.
		Compute the functional $\mathcal J(\textbf{x},\textbf{u}^*)$.
		\item[\texttt 3.] Perform a perturbation of the trajectories over  a fixed set of points $P^*(t)$ on current optimized leader trajectories with small random variations over the time-set $S_M$,
		\begin{equation}\label{eq:randvar}
			P^{(j)}(t_m) = P^*(t_m)+B_m,\quad m=1,\ldots,M, \tag{$\mathcal P$}
		\end{equation}
		where $B_m\sim \textrm{Unif}([-1,1]^d)$ is a random perturbation and set for $m=1,\ldots,M,$
		\[
		u^{\textrm{opt},(j)}(t) = \frac{P^*(t_{m+1}) - P^*(t_m)}{\Vert P^*(t_{m+1}) - P^*(t_m) \Vert},\quad  t\in[t_m,t_{m+1}].
		\]
		Finally compute $ \mathcal J(\textbf{x},\textbf{u}^{(j)})$.
		\item[\texttt 4.] \texttt{while} $j<j_{\texttt{max}}$ AND $\mathcal J(\textbf{x},\textbf{u}^*)<\mathcal J_E$
		\begin{enumerate}
			\item Update $ j\leftarrow j+1$.
			\item Perform the perturbation \eqref{eq:randvar} and compute $\mathcal J(\textbf{x},\textbf{u}^{(j)})$.
			\item \texttt{If} $\mathcal J(\textbf{x},\textbf{u}^{(j)})\leq \mathcal J(\textbf{x},\textbf{u}^{*})$
			\\ \qquad set $u^{*}\leftarrow u^{(j)}$ and $\mathcal J(\textbf{x},\textbf{u}^{*})\leftarrow \mathcal J(\textbf{x},\textbf{u}^{(j)})$.
		\end{enumerate}
		\texttt{repeat}
	\end{enumerate}
\end{alg}
\begin{remark}~
\begin{itemize}
	\item Compass search does not guarantee the convergence to a global minimizer, on the other hand it offers a good compromise in terms of computational efficiency. 
	\item Alternative metaheuristic schemes can be emploied to enanche leader trajectories and improbing the convergence towards the global minimizer, among several possibilities we refer to genetic algorithms, and particle swarm based optimizations.
	\item The synthesis of control strategies via compass search for the microscopic and the mean-field dynamics can produce different results, due to the strong non-linearities of the interactions, and the non-convexity of the functional considered, such as the evacuation time. However, in any case, the solutions retrived by this approach statisfy a local optimality criteria by construction.

\end{itemize}	
\end{remark}

\subsection{Numerical experiments}
We present three different numerical experiments at microscopic and mesoscopic levels, corresponding to the minimization of cost functionals presented in Section \ref{sec:opt2}.

{\em Numerical discretization.}
The dynamics at microscopic level is discretized by a forward Euler scheme with a time step $\Delta t = 0.1$, whereas the evolution of the mean-field dynamics is approximated by MFMCs algorithms. We choose a sample of $\mathcal{O}(10^3)$ particles for the approximation of the density and we reconstruct their evolution in the phase space by kernel density estimator with a multivariate standard normal density function with bandwidth $h=0.4$. Table \ref{tab:all_parameters} reports the parameters of the model for the various scenarios unchanged for every test. The number of leaders instead changes and it will be specified later.
\begin{table}[!h]
	\label{tab:all_parameters}
	\begin{center}
		\caption{Model parameters for the different scenarios.}
		\begin{tabular}{ccccccccccccccccc}
			 &$N^F$ &$\mathcal{N}$ & $C_r^F$ & $C_r^L$ & $C_{a}^L$ &$C_{a}^F$ & $C_{\tau}$ & $C_s$ & $s^2$ & $r\equiv\zeta$ & $\gamma$\\
			\cmidrule(r){2-2}\cmidrule(r){3-3}\cmidrule(r){4-4}\cmidrule(r){5-5}\cmidrule(r){6-6}\cmidrule(r){7-7}\cmidrule(r){8-8}\cmidrule(r){9-9}\cmidrule(r){10-10}
			\cmidrule(r){11-11}\cmidrule(r){12-12}\cmidrule(r){13-13}\\[-1em]
			\textbf{}
	&150   & 20 & 2 & 1.5 & 3 & 3 & 1 & 0.5 & 0.4 & 1 &1 \\
			\hline
		\end{tabular}
	\end{center}
\end{table}

{\em Obstacles handling.} In order to deal with obstacles we use a {\em cut-off velocity} approach, namely we  compute the velocity field first neglecting the presence of the obstacles, then nullifying the component of the velocity vector which points inside the obstacle. This method is used in, e.g., \cite{albi2016SIAP, cristiani2011MMS, cristiani2015SIAP} and requires additional conditions to avoid situations where pedestrians stop walking completely because both components of the velocity vector vanish, e.g. in presence of corners, or when obstacles are very close to each other. We refer to \cite{cristiani2017AMM} for more sophisticated approaches of obstacles handling.

\subsubsection{Test 1: Minimum time evacuation with multiple exits}\label{sec:test1}
In this first test, leaders aim to minimize the time of evacuation \eqref{eq:minevactime}, hence trying to enforce crowd towards the exit avoiding congestion and ease the outflux of the pedestrian. We assume that leaders informed about exits position follow `go-to-target' strategy defined as in \eqref{eq:gototarget_beta}, where the target is defined by the different exits and will be specified for each leader. In what follows we account for two different settings comparing microscopic and mesoscopic dynamics.

\begin{figure}[h!]
	\centering
	\includegraphics[width=0.45\linewidth]{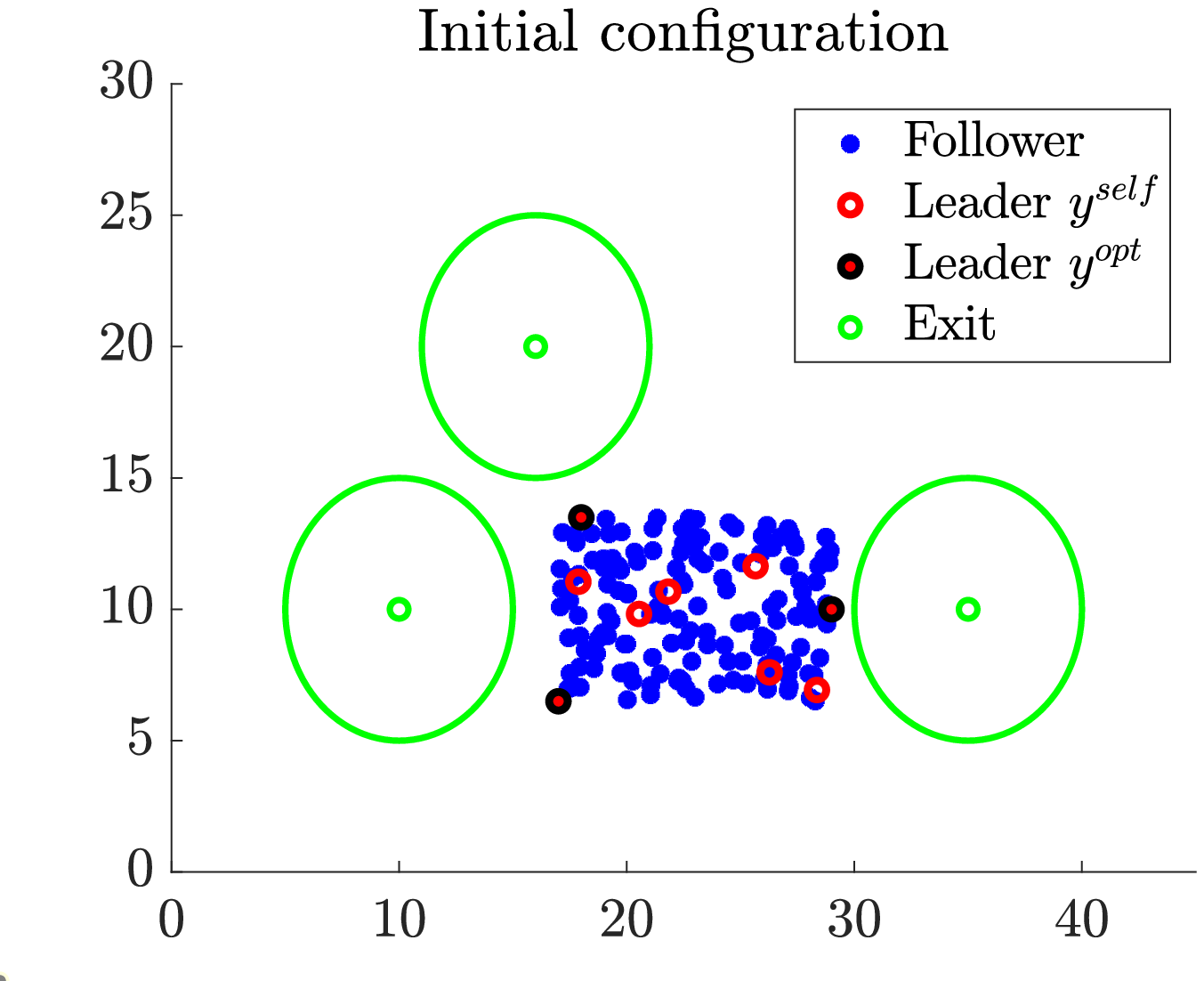}
		\includegraphics[width=0.45\linewidth]{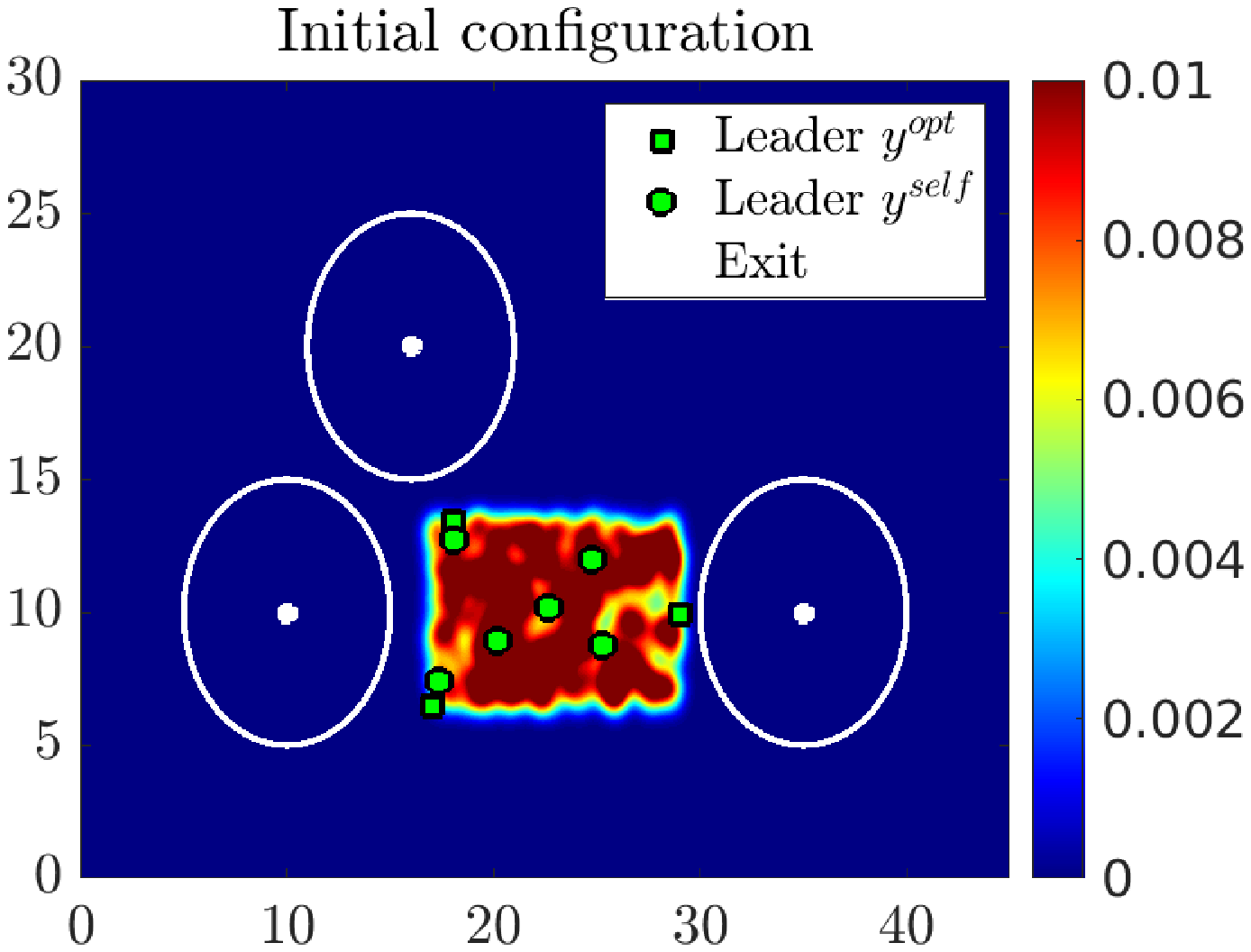}
	\caption{{\em Test 1a}. Minimum time evacuation with multiple exits, initial configuration.   }
	\label{fig:test1_initconf}
\end{figure}

\paragraph{\bf Setting a) Three exits} We consider the case of a room with no obstacles and three exits located at $x_1^\tau = (35, 10)$, $x_2^\tau = (16, 20)$, $x_3^\tau = (10, 10)$ with visibility areas $\Sigma_e = \left\lbrace x \in \mathbb{R}^2 : |x - x_e^\tau| < 5 \right\rbrace$. We consider two different types of leaders, we call selfish leaders $y^{\rm self}$ the agents who do not care about followers and follow the direction that connects their positions to the exits.
While the optimized leaders $y^{\rm opt}$ are aware of their role and they move with the aim to reach the exits and to maintain contact with the crowd, only the trajectories of this type of leaders will be optimized.
The admissible leaders trajectories are defined as in Equation \eqref{eq:gototarget_beta}, we choose $\beta = 1$ for selfish leaders, $\beta = 0.6$ for optimized leaders and the target position as $\Xi_k(t)= x_e^\tau \ \forall t$ and for every leader $k$.
At initial time leaders and followers are uniformly distributed in the domain $[17,29]\times[6.5,13.5]$ where followers velocities are sampled from a normal distribution with average $-0.5$ and variance $0.1$, hence biased towards the wrong direction. We report in Figure \ref{fig:test1_initconf}  the initial configuration for both microscopic and mesoscopic dynamics.

{\em Microscopic case.}
We consider $N^L=9$ leaders, three optimized and six unaware leaders. Each leader is associated with an exit: unaware leaders move towards the nearest exit, whereas each optimized leader is assigned to a different exit.

\begin{figure}[h!]
	\centering
	\includegraphics[width=0.328\linewidth]{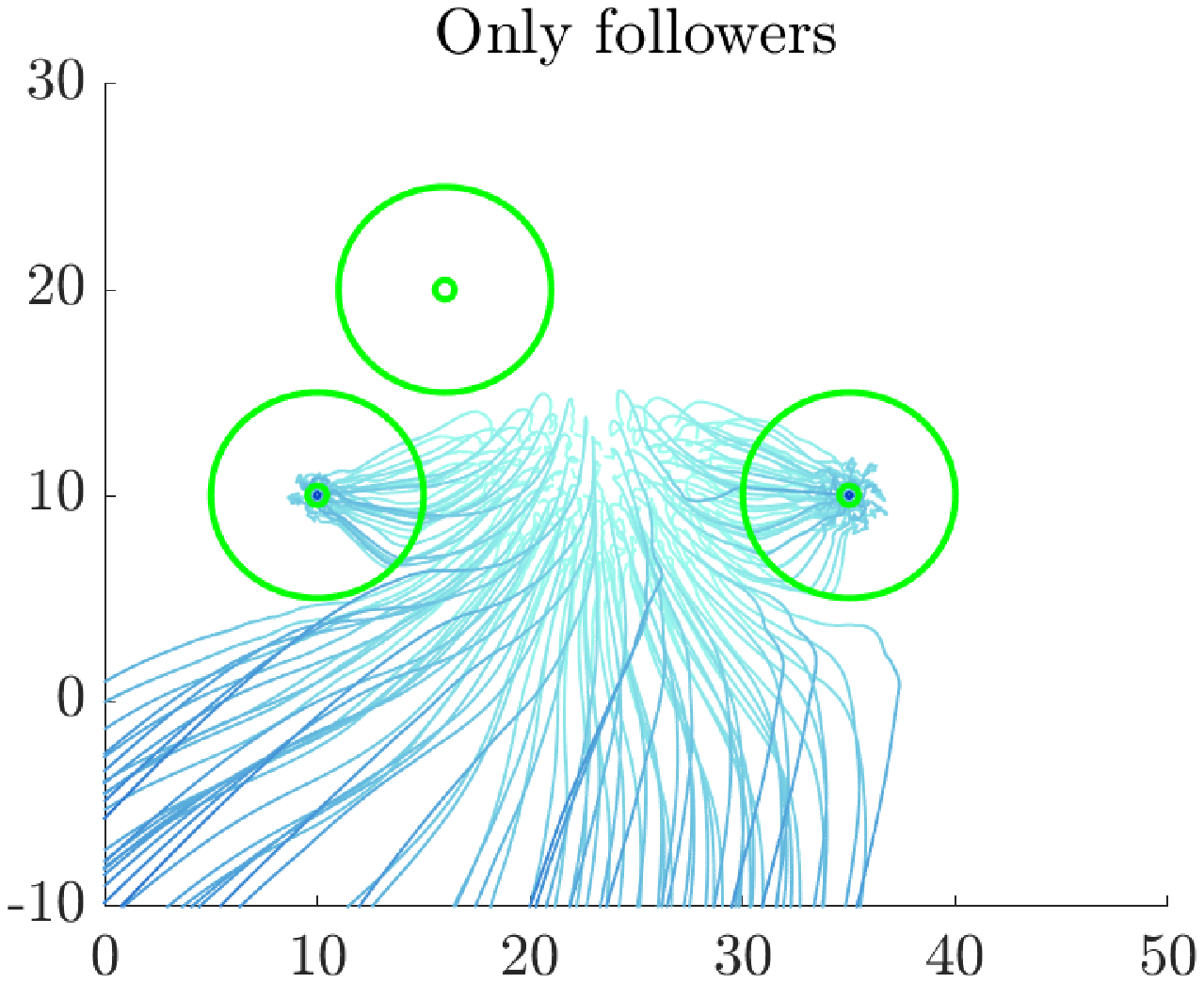}
	\includegraphics[width=0.328\linewidth]{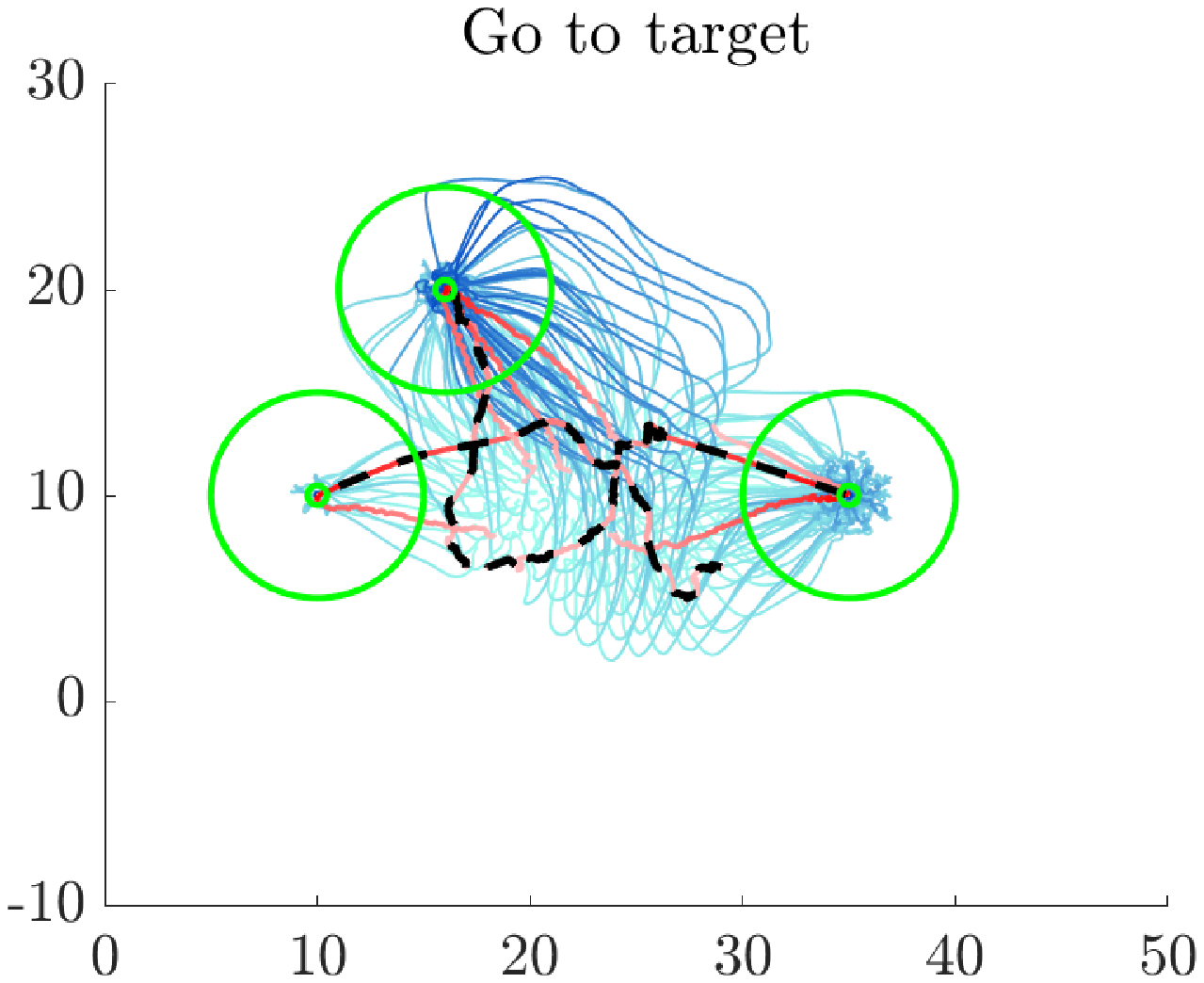}
	\includegraphics[width=0.328\linewidth]{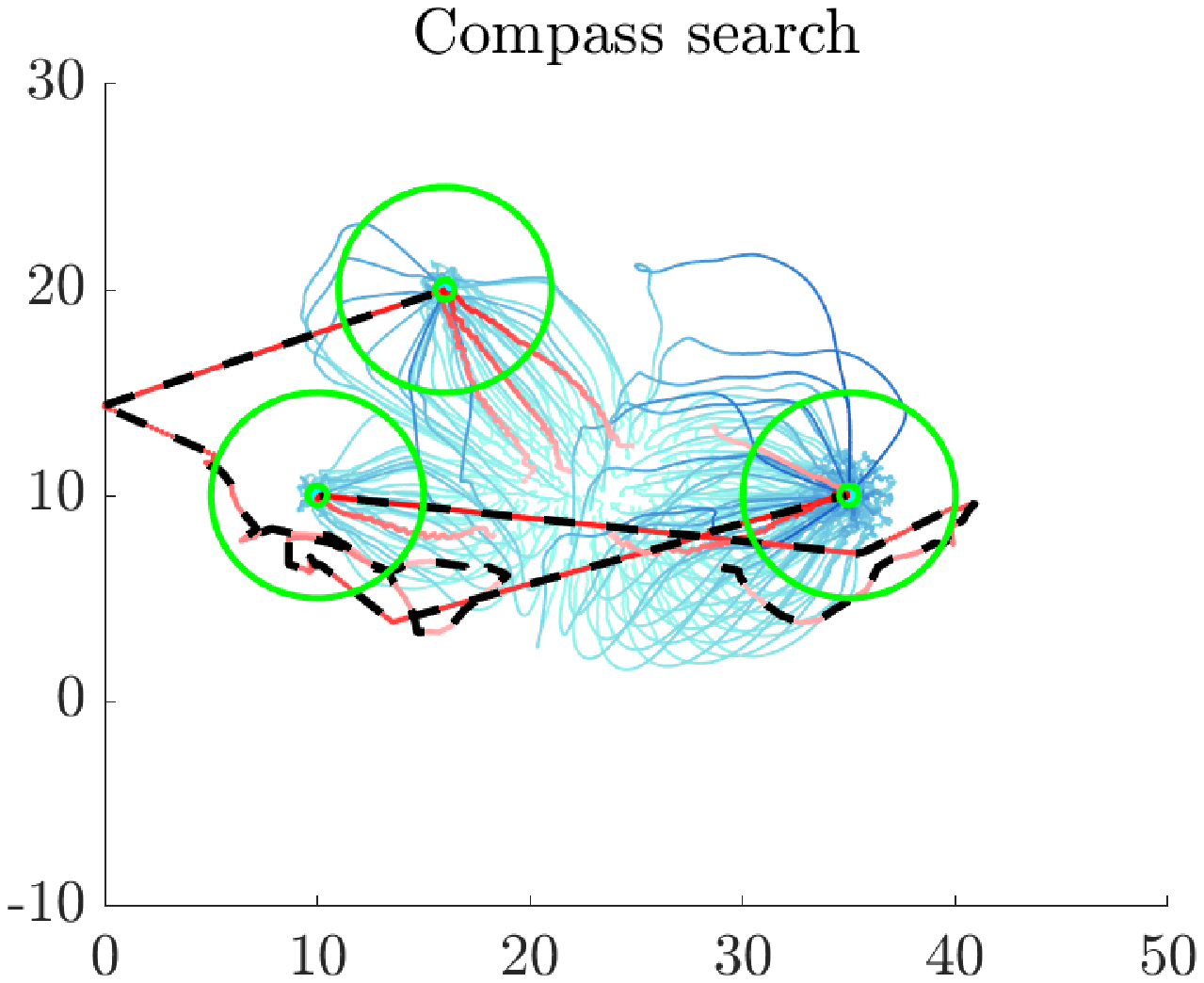}
	\caption{{\em Test 1a}. Microscopic case: minimum time evacuation with multiple exits. On the left the uncontrolled case, in the centre the go-to-target and on the right the optimal compass search strategy.}
	\label{fig:test1_micro_1}
\end{figure}

\begin{figure}[h]
	\centering
	\includegraphics[width=0.45\linewidth]{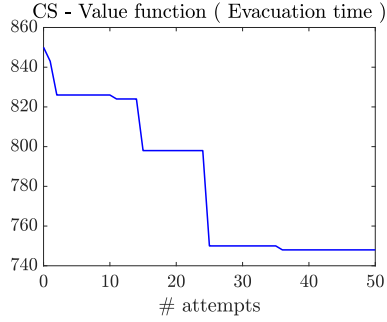}
	\caption{{\em Test 1a}. Microscopic case:  decrease of the value function \eqref{eq:minevactime} as a function of compass search iteration. }
	\label{fig:test1_micro_3}
\end{figure}

\begin{figure}[h!]
	\centering
	\includegraphics[width=0.328\linewidth]{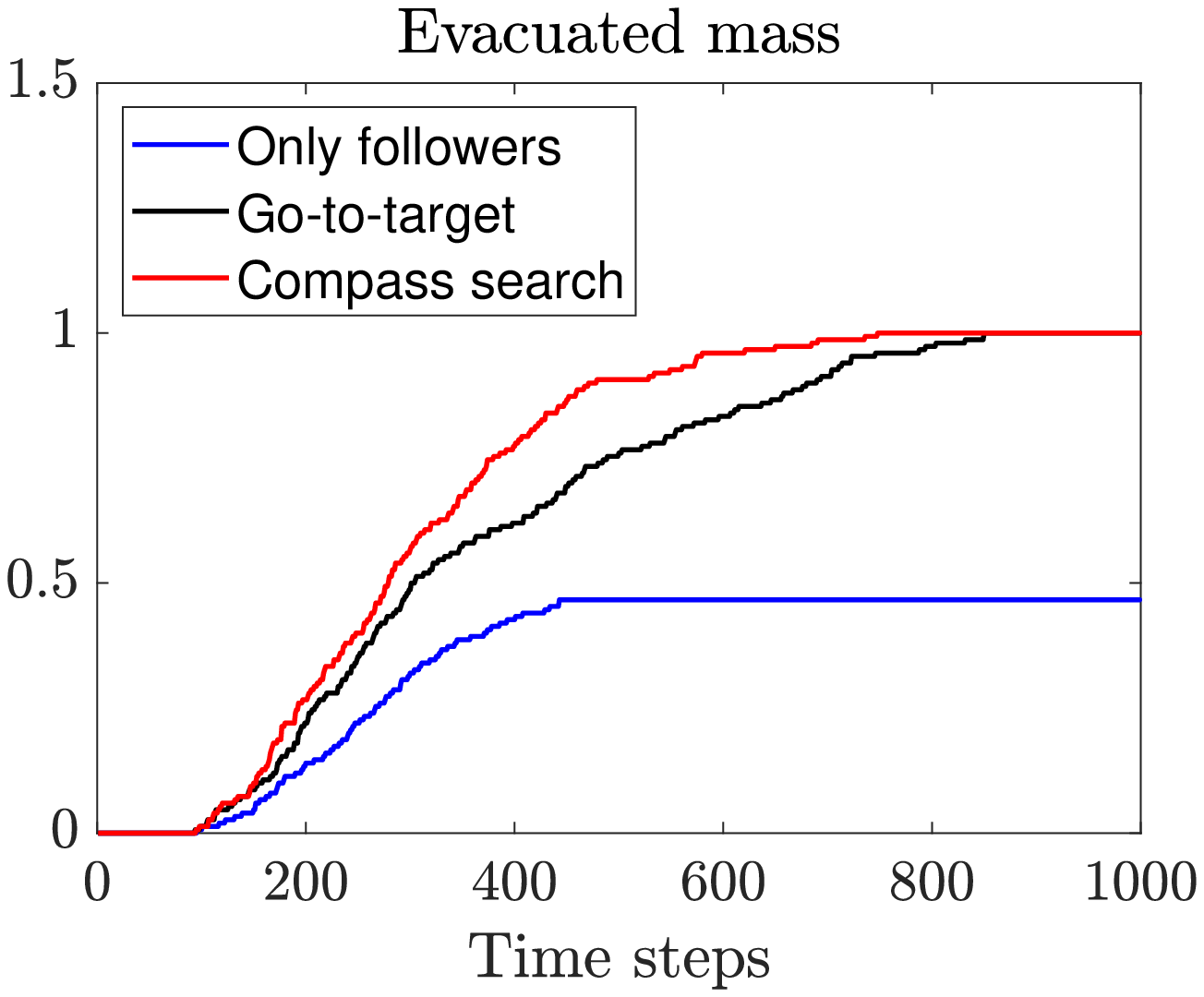}
	\\
	\includegraphics[width=0.328\linewidth]{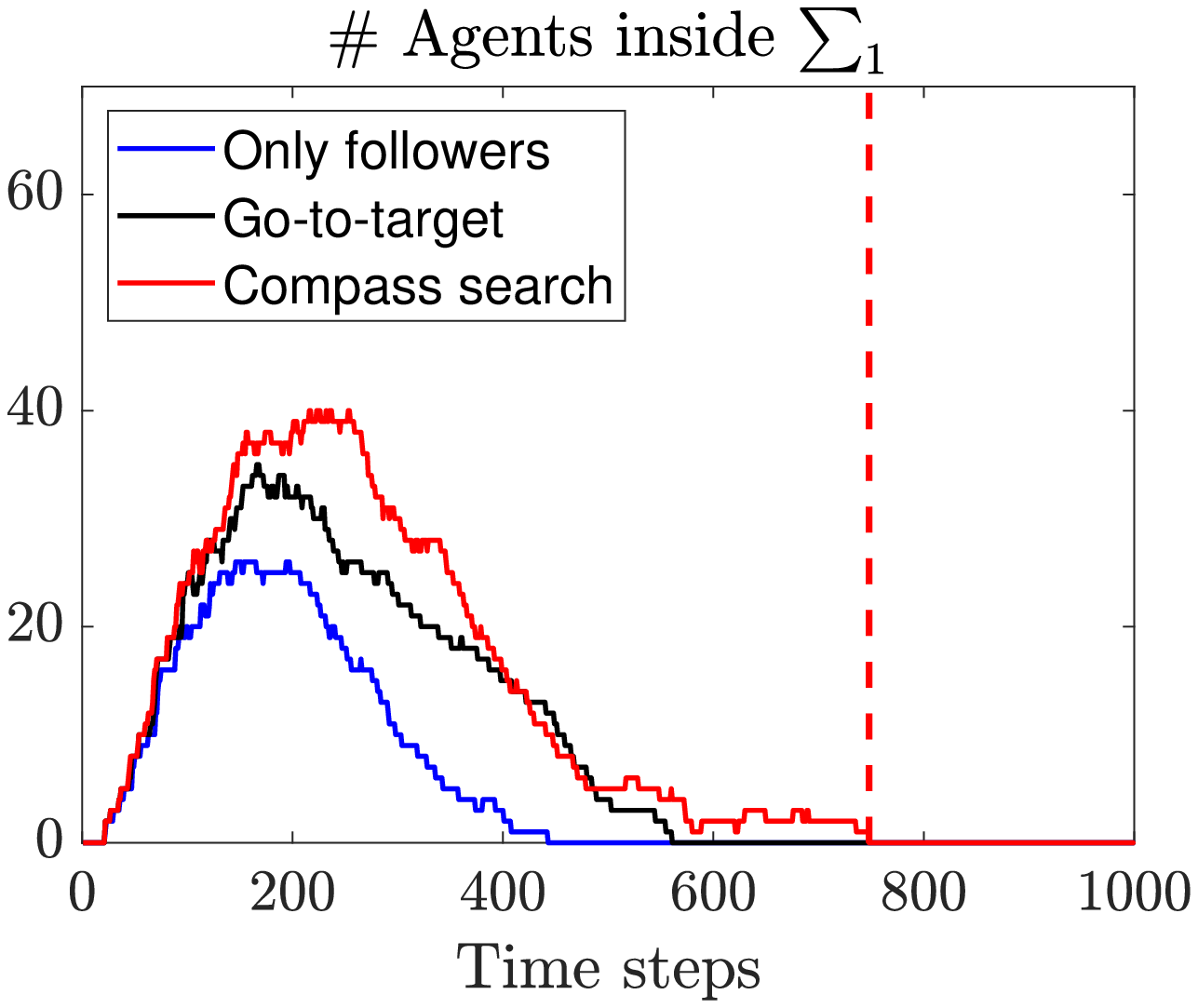}
	\includegraphics[width=0.328\linewidth]{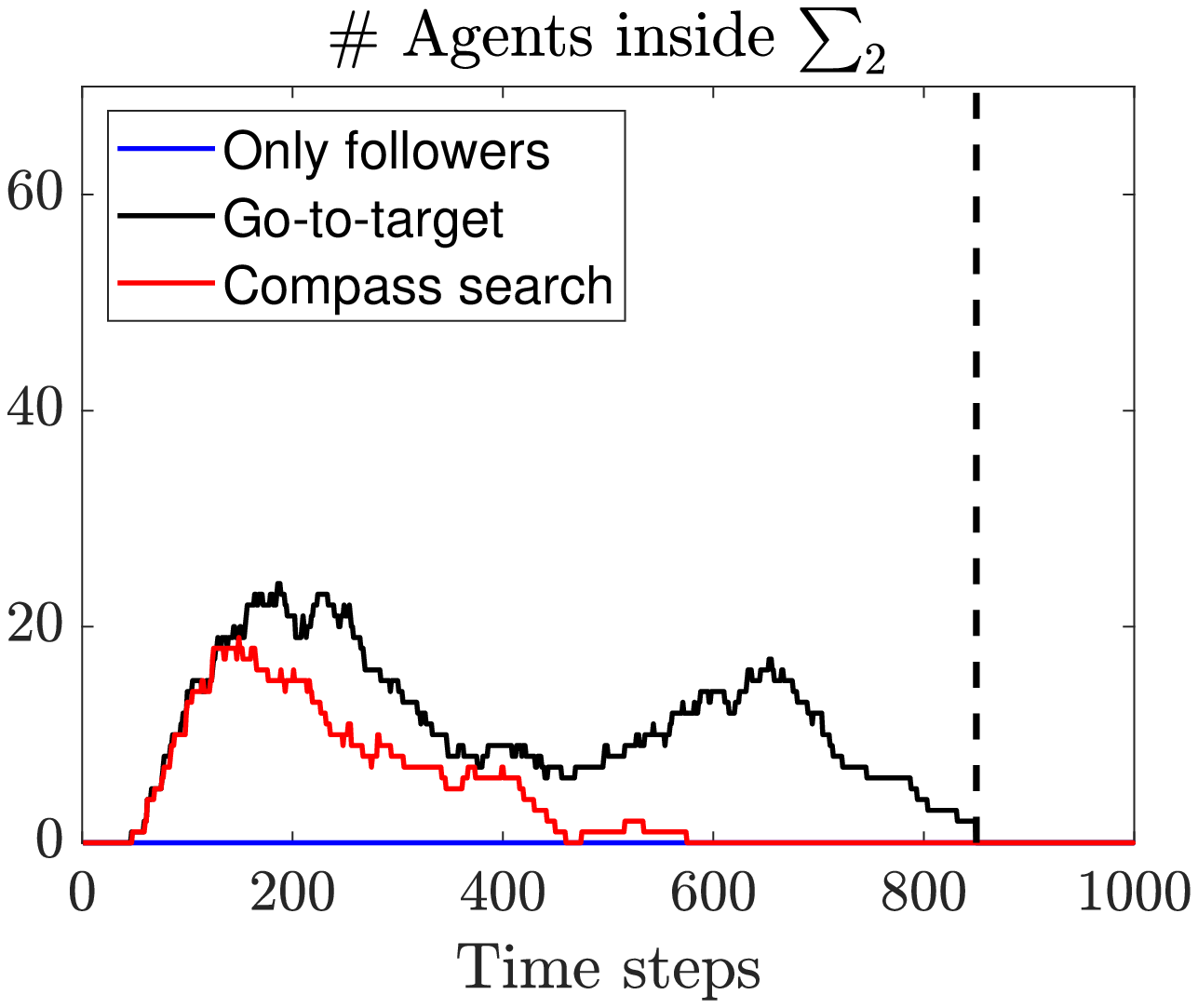}
	\includegraphics[width=0.328\linewidth]{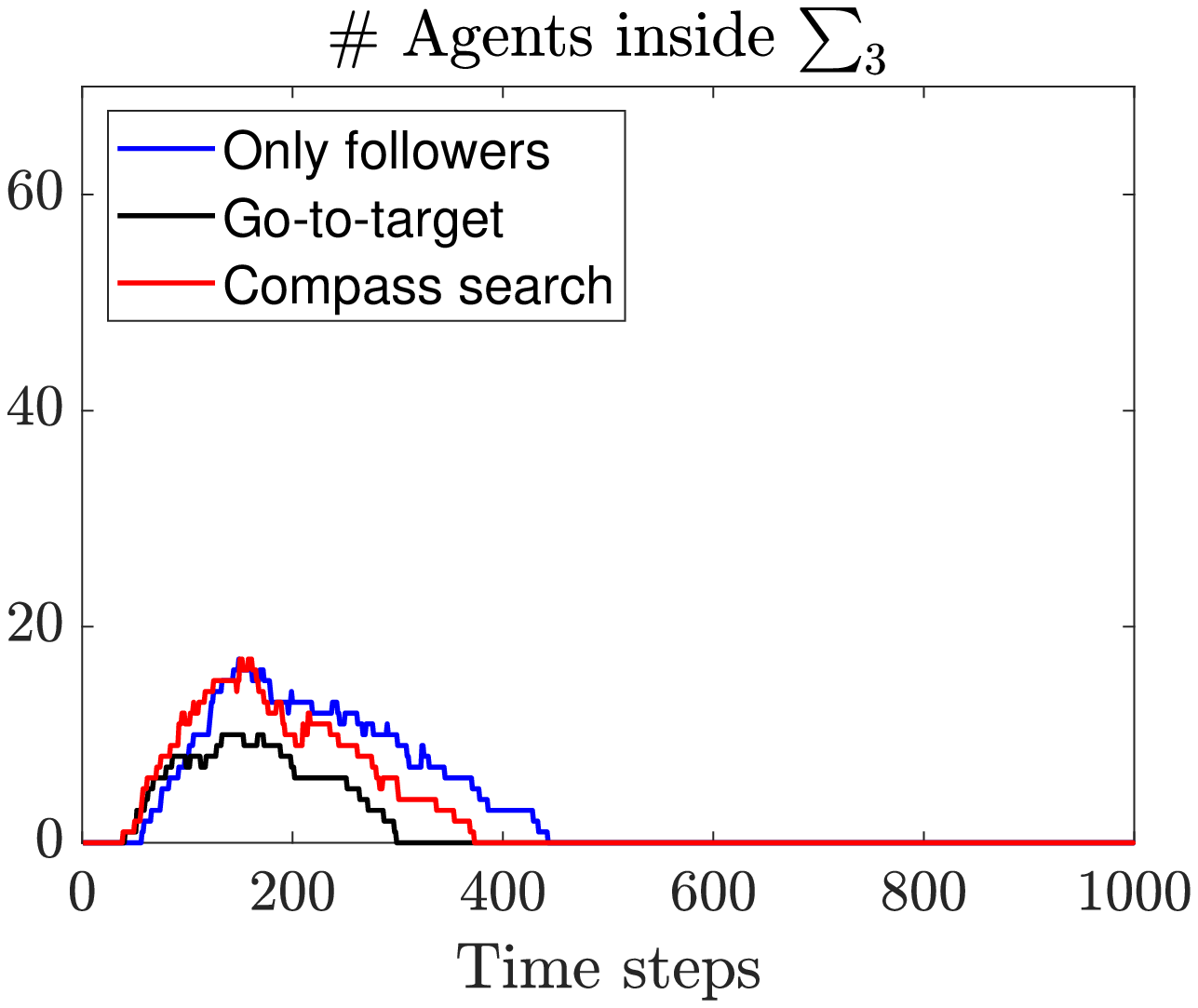}
	\caption{{\em Test 1a.} Microscopic case: minimum time evacuation with multiple exits. Evacuated mass (first row), occupancy of the visibility area $\Sigma_1$ (second row, left), $\Sigma_2$ (second row, centre) and $\Sigma_3$ (second row, right) as a function of time for uncontrolled, go-to-target and optimal compass search strategies. The dot line denotes the time step in which the whole mass is evacuated, the line is black for the go-to-target and red for the optimal compass search strategy.}
	\label{fig:test1_micro_2}
\end{figure}

Figure \ref{fig:test1_micro_1} shows the evolution of the agents with the go-to-target strategy on the left and with the optimal strategy obtained by the compass search algorithm on the right. As it can be seen in Figure \ref{fig:test1_micro_1} with the go-to-target strategy the whole crowd reaches the exit, after $850$ time steps. We distinguish optimized leaders $y^{opt}$ with a dashed black line. 
Optimized movements for leaders are retrived by means of Algorithm \ref{alg_CS}, with initial guess go-to-target strategy, we report in Figure \ref{fig:test1_micro_3} the decrease of the performance function \eqref{eq:minevactime} as a function of the iterations of compass search.
Eventually optimized leaders influence the crowds for a larger amount of time and the total mass is evacuated after 748 time steps, as shown in Table \ref{eq:test1_table1}. 
\begin{table}
	\centering
	\caption{{\em Test 1a}. Performance of leader strategies over microscopic dynamics.}\label{eq:test1_table1}
\begin{tabular}{ccccc}
& uncontrolled & go-to-target  & CS ($50$ it) \\
\cmidrule(r){2-2}\cmidrule(r){3-3}\cmidrule(r){4-4}
Evacuation time (time steps)& > 1000 & 850  & 748\\
Evacuated mass (percentage) & 46\% & 100\% & 100\%\\
\hline
\end{tabular}
\end{table}
Figure \ref{fig:test1_micro_2} compares the evacuated mass and the occupancy of the exits visibility zone as a function of time for the uncontrolled case, the go-to-target strategy and the optimal compass search strategy. Dashed lines indicate times of total mass evacuation.

\begin{figure}[h!]
	\centering
	\includegraphics[width=0.45\linewidth]{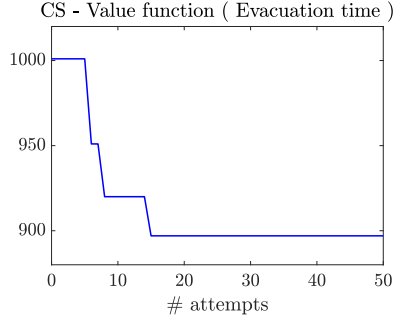}
	\caption{{\em Test 1a.} Mesoscopic case: minimum time evacuation with multiple exits. Decrease of the value function \eqref{eq:minevactime} as a function of attempts.}
	\label{fig:test1nvl_time_meso_4}
\end{figure}

\begin{figure}[h!]
	\centering
	\includegraphics[width=0.328\linewidth]{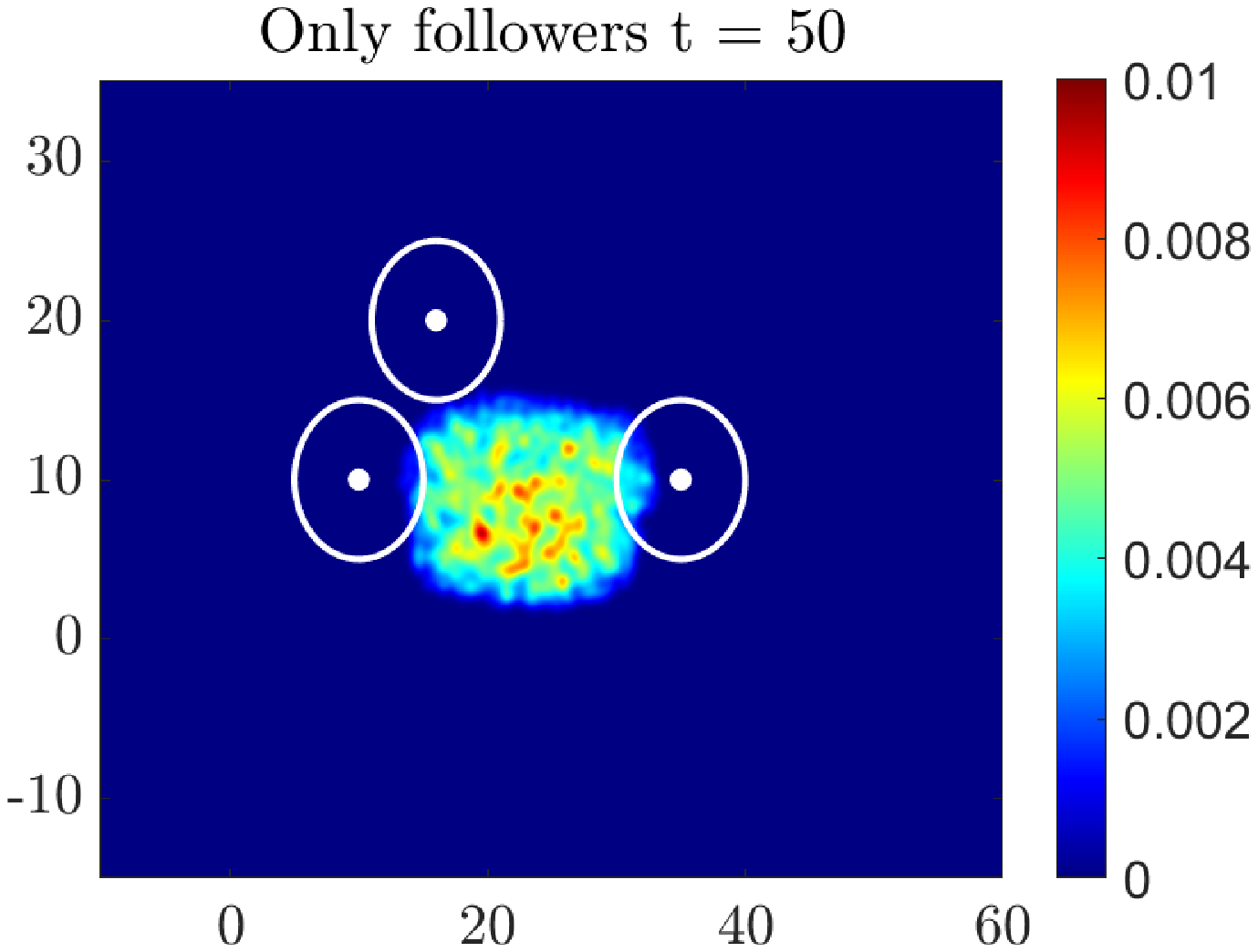}
	\includegraphics[width=0.328\linewidth]{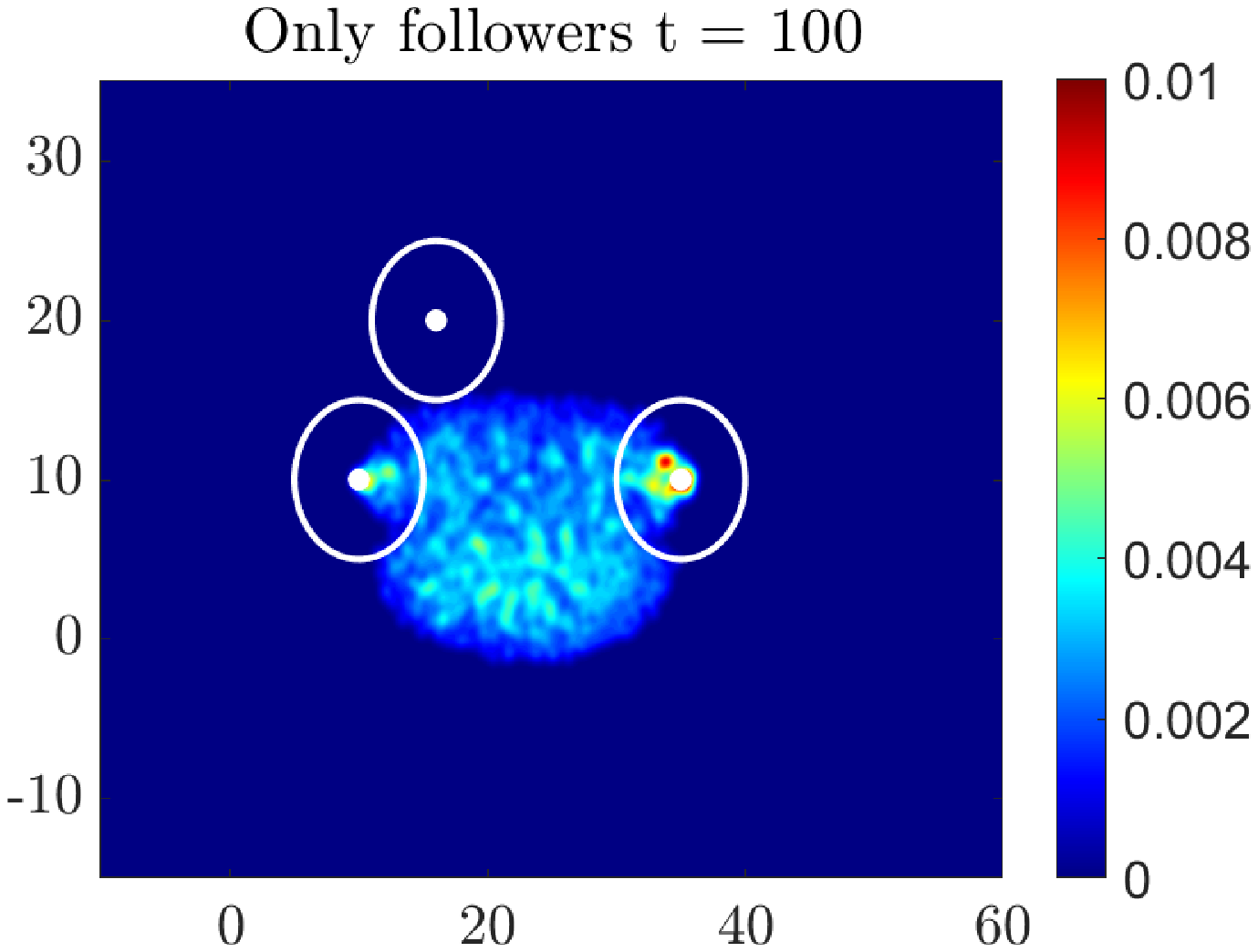}
	\includegraphics[width=0.328\linewidth]{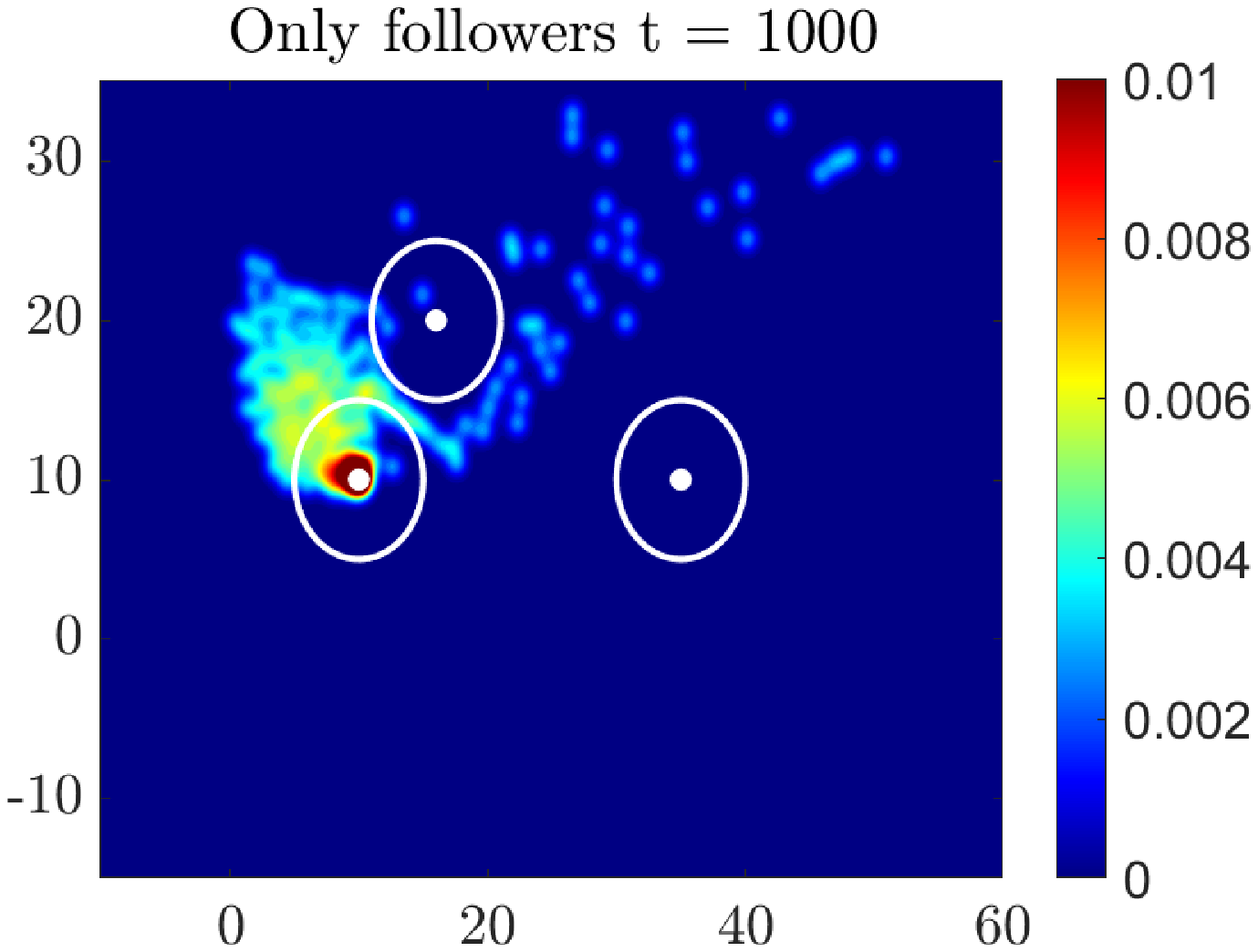}
	\\
	\includegraphics[width=0.328\linewidth]{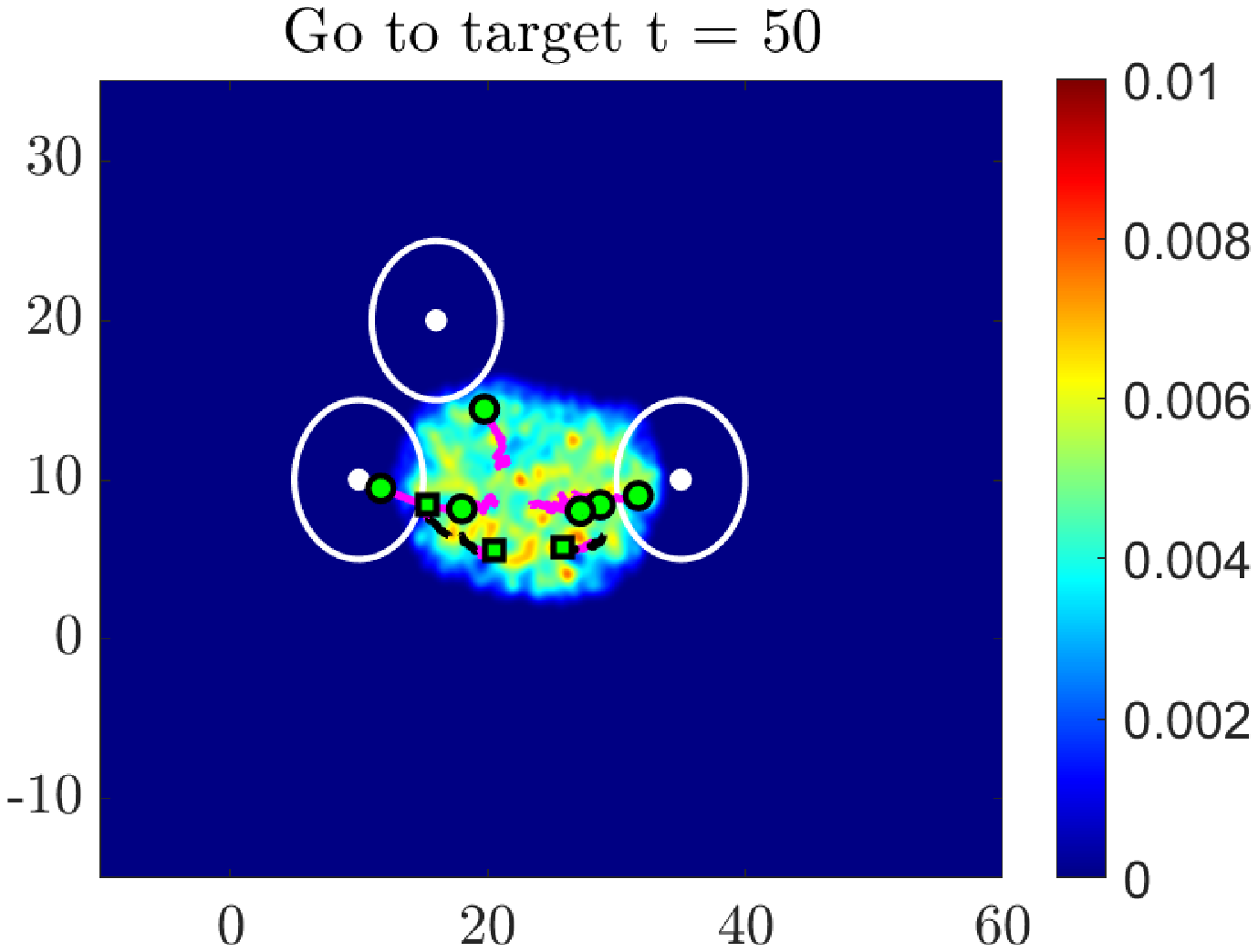}
	\includegraphics[width=0.328\linewidth]{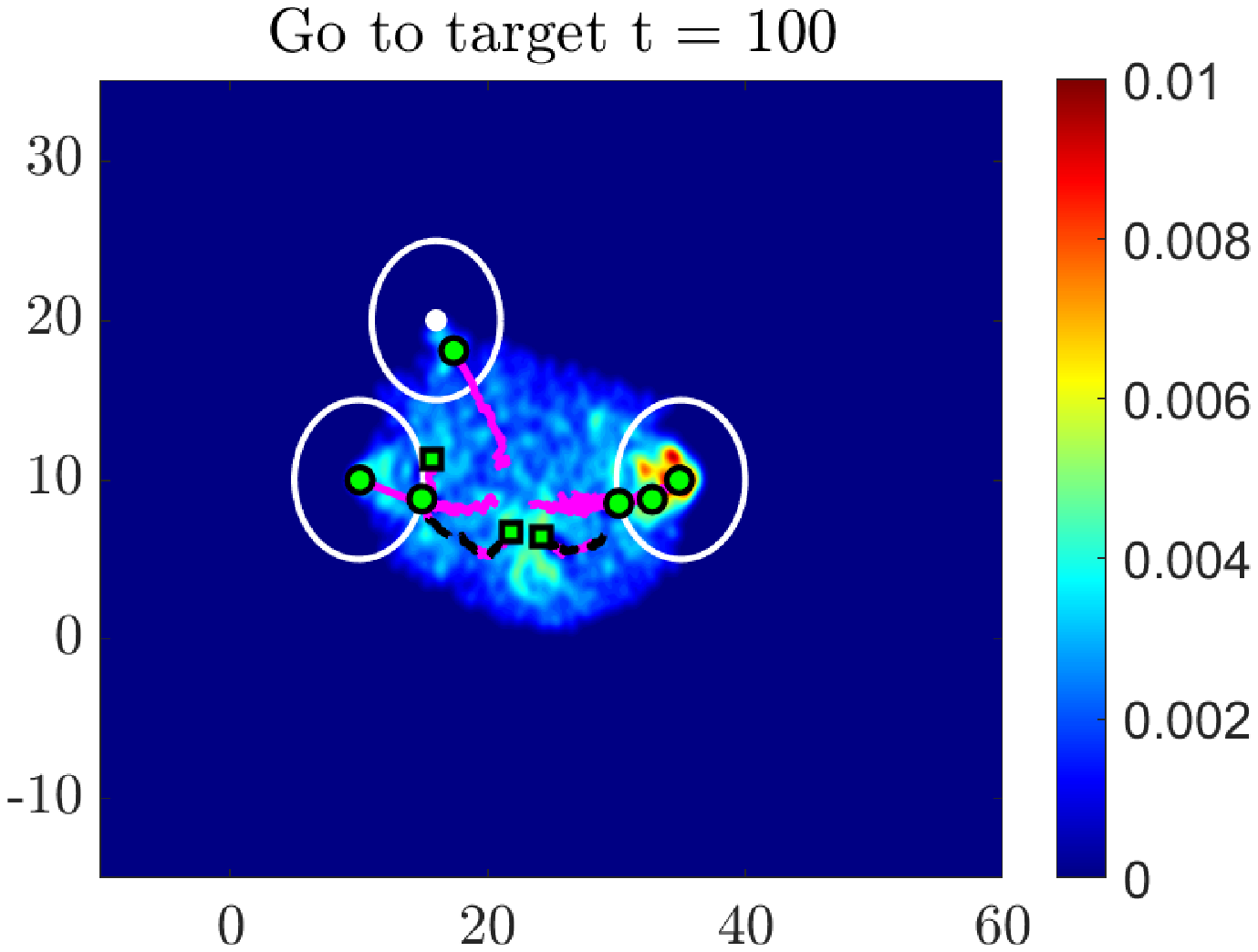}
	\includegraphics[width=0.328\linewidth]{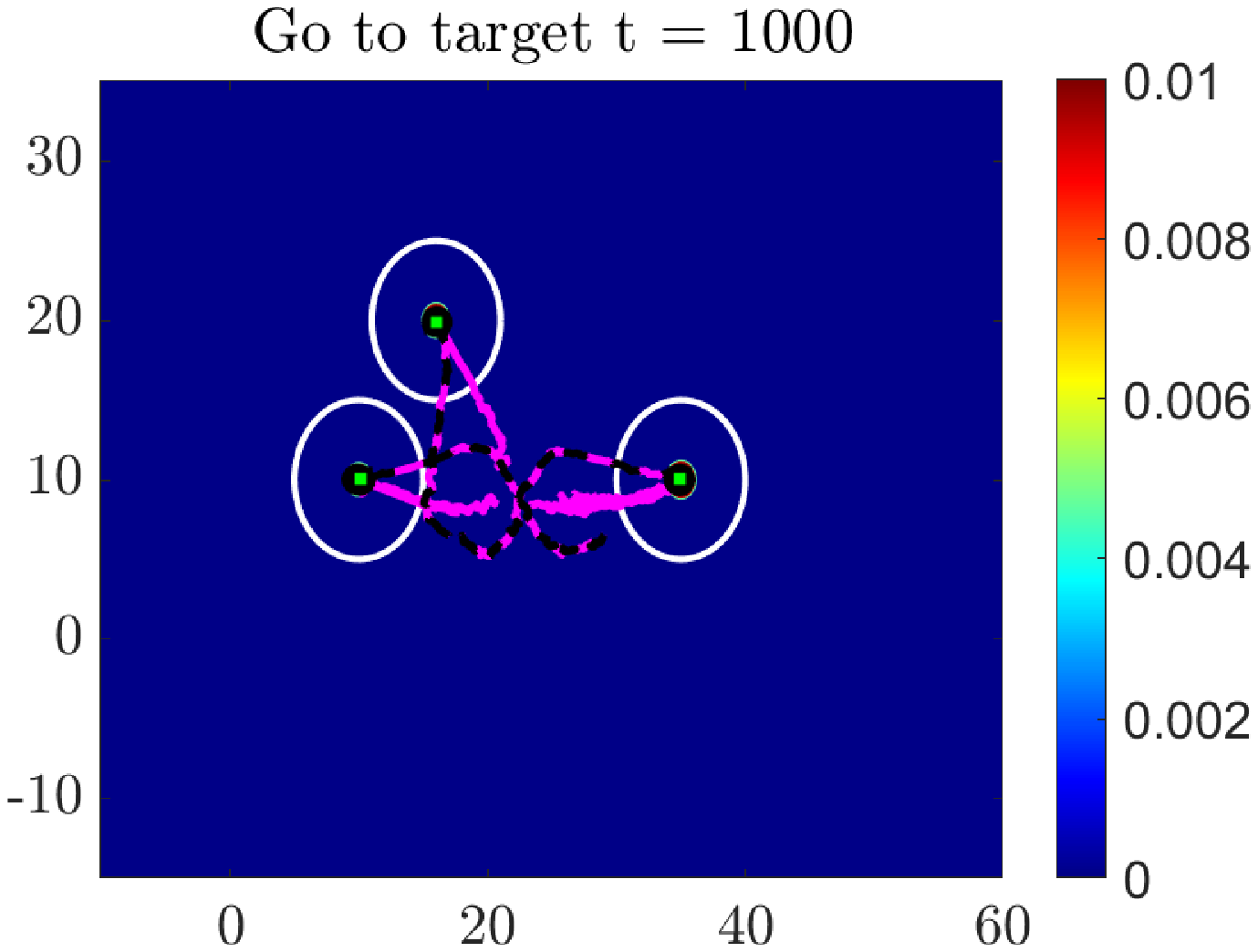}
	\\
		\includegraphics[width=0.328\linewidth]{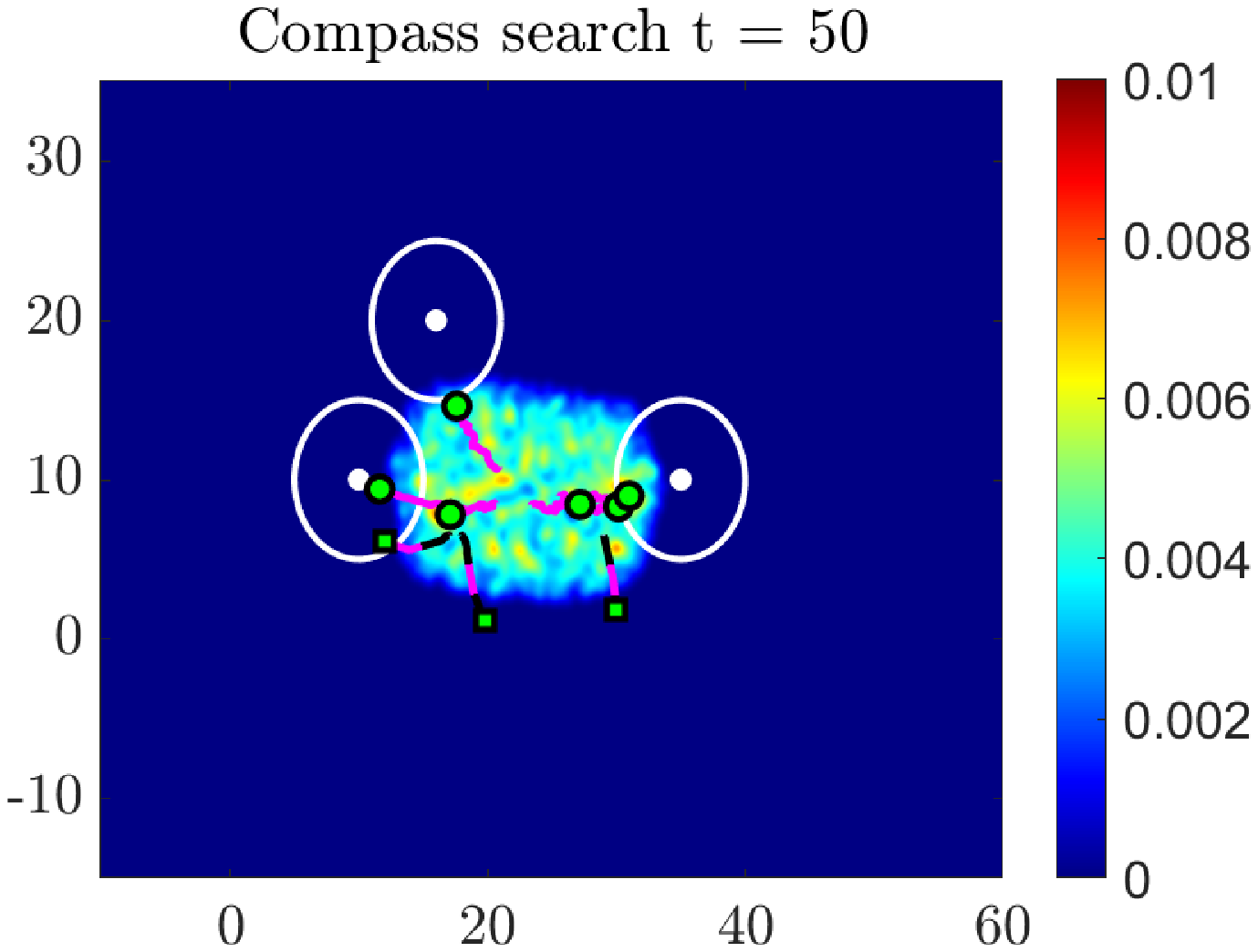}
	\includegraphics[width=0.328\linewidth]{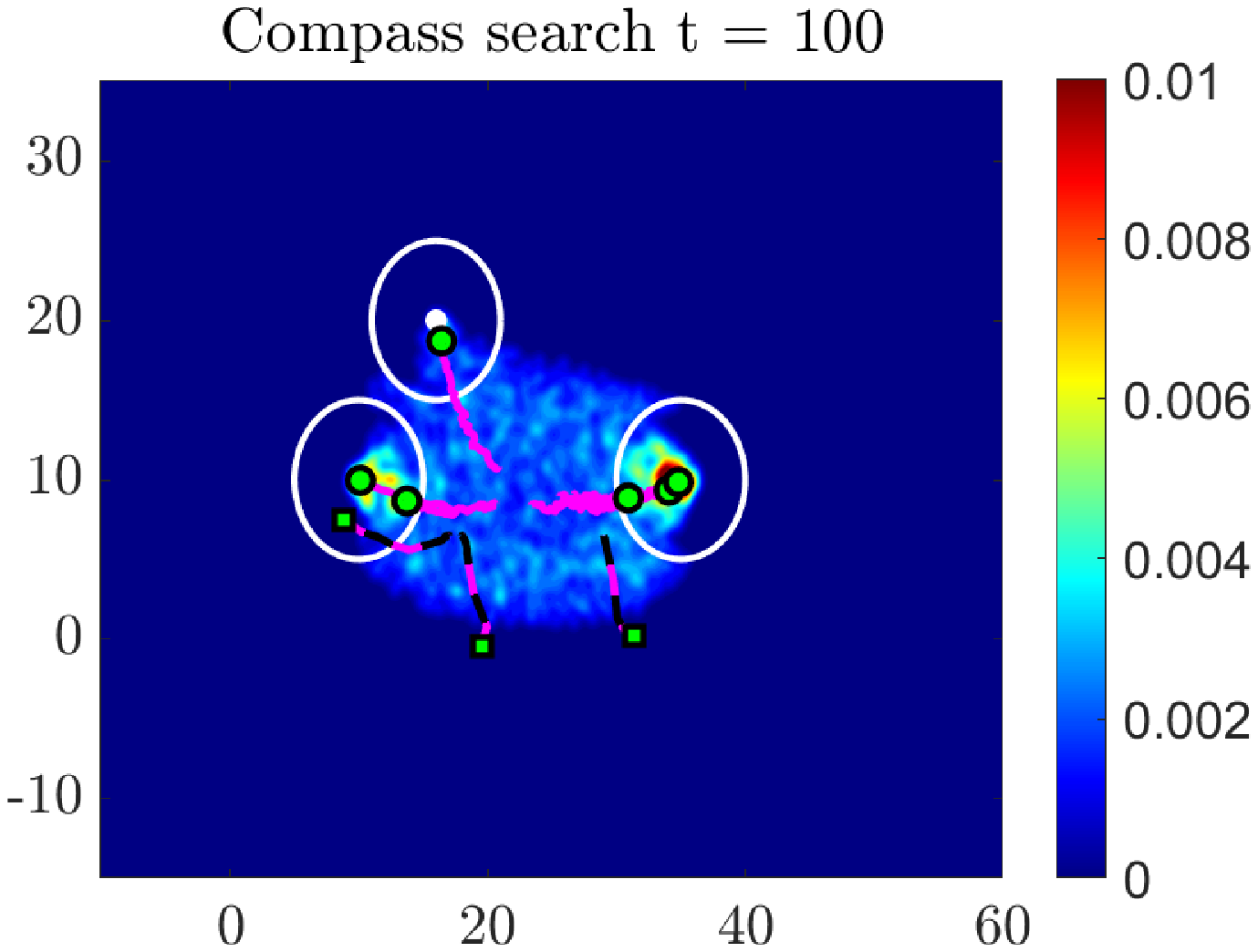}
	\includegraphics[width=0.328\linewidth]{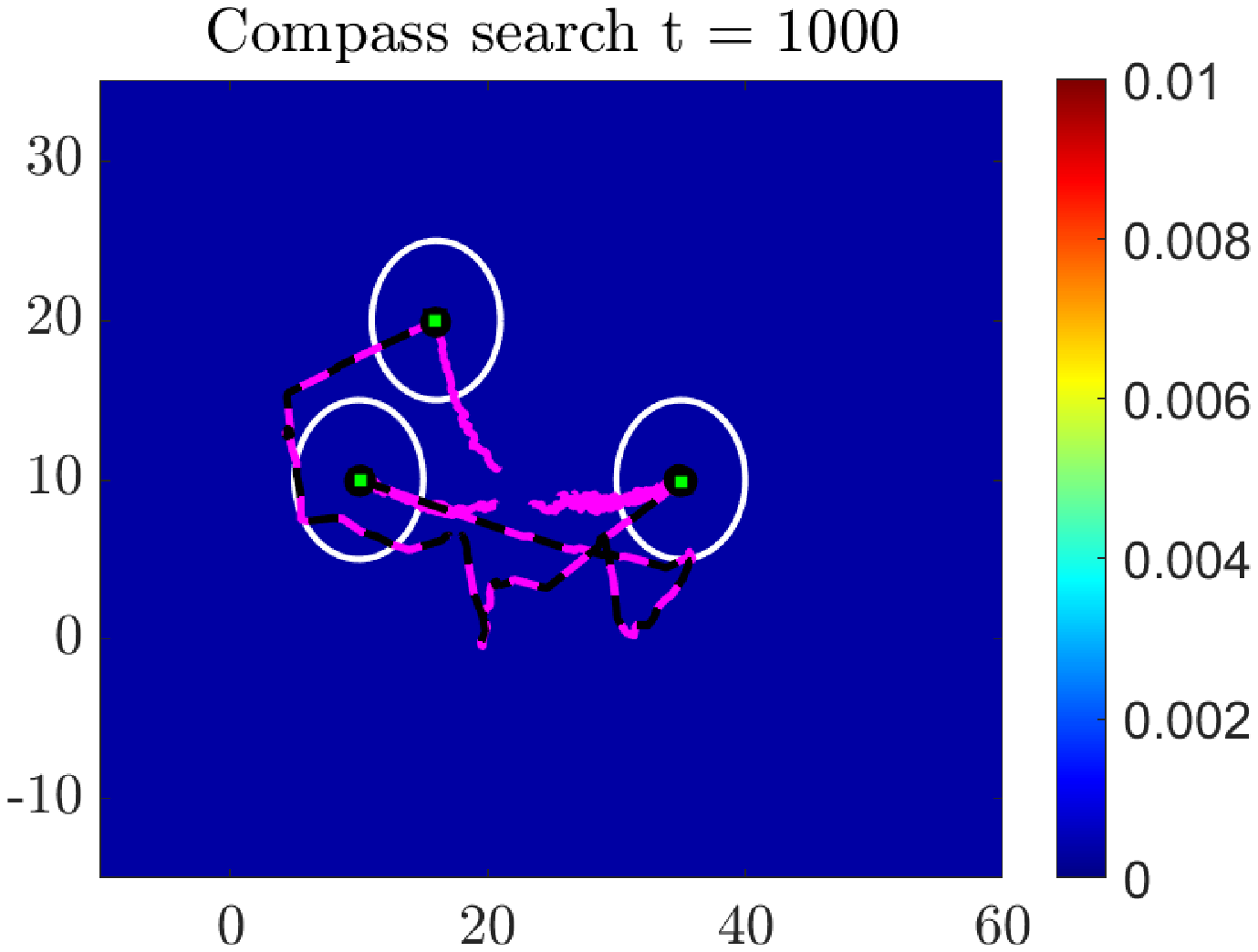}	
\caption{{\em Test 1a.} Three snapshots taken at time $t=50$, $t=100$, $t=1000$ of the mesoscopic densities for the  minimum time evacuation with multiple exits. In the upper row the uncontrolled case, in the central row the three aware leaders follows a go-to-target strategy, whereas in the bottom row their trajectories are optimized according to CS algorithm.}
	\label{fig:test1nvl_time_meso}
\end{figure}

{\em Mesoscopic case.} We consider now a continuous density of followers, in the same setting of the previous microscopic case: we account for $N^L=9$ microscopic leaders moving in a room with no obstacles and three exits. Hence we compare uncontrolled dynamics, go-to-target strategies, and optimized strategies with compass search.  In Table \ref{tab_nvl_meso} we show that without any control followers are unable to reach the total evacuation reaching $84\%$ of total mass evacuated. Go-to-target strategy improves total mass evacuated, however, a small part of the mass spreads around the domain and is not able to reach the target exit. Eventually, with optimized strategies, we reach the evacuation of the total mass in $897$ simulation steps. 
\begin{table}
	\centering
	\caption{{\em Test 1a}. Performance of leader strategies over mesoscopic dynamics.}\label{tab_nvl_meso}
\begin{tabular}{cccc}
& uncontrolled & go-to-target  & CS ($50$ it)\\
\cmidrule(r){2-2}\cmidrule(r){3-3}\cmidrule(r){4-4}
Evacuation time (time steps)  & > 1000 & > 1000 & 897\\
Evacuated mass (percentage)  & 84\% & 99\% & 100\%\\
\hline
\end{tabular}
\end{table}
The better performance of the optimized strategy can be observed directly from Figure \ref{fig:test1nvl_time_meso_4}, where functional \eqref{eq:minevactime} is evaluated at subsequent iterations of Algorithm \ref{alg_CS}.
In Figure \ref{fig:test1nvl_time_meso} we show three snapshots of the followers density comparing leaders with different strategies and the uncontrolled case.
In the upper row, we report the evolution without any control. The middle row shows leaders driven by a go-to-target strategy promoting evacuation of followers density. At time $t=50$ leaders are moving to influence the followers towards the three exits. At time $t=100$, the followers mass splits and starts to reach the exits.  At time $t=1000$, complete evacuation is almost reached. The bottom row depicts improved strategies of leaders, where total mass is evacuated at time step $912$.

\begin{figure}[h!]
	\centering
	\includegraphics[width=0.328\linewidth]{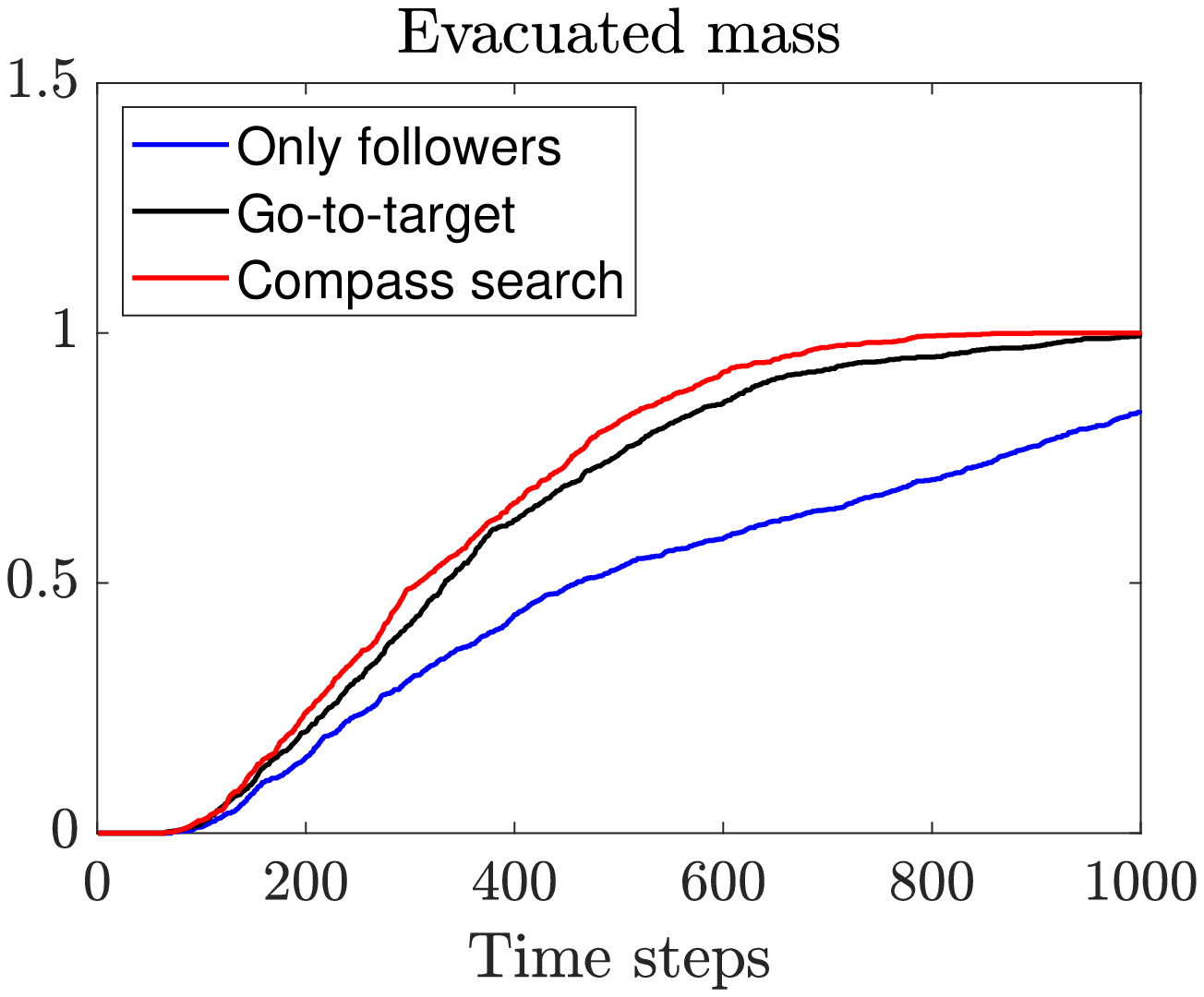}
	\\
	\includegraphics[width=0.328\linewidth]{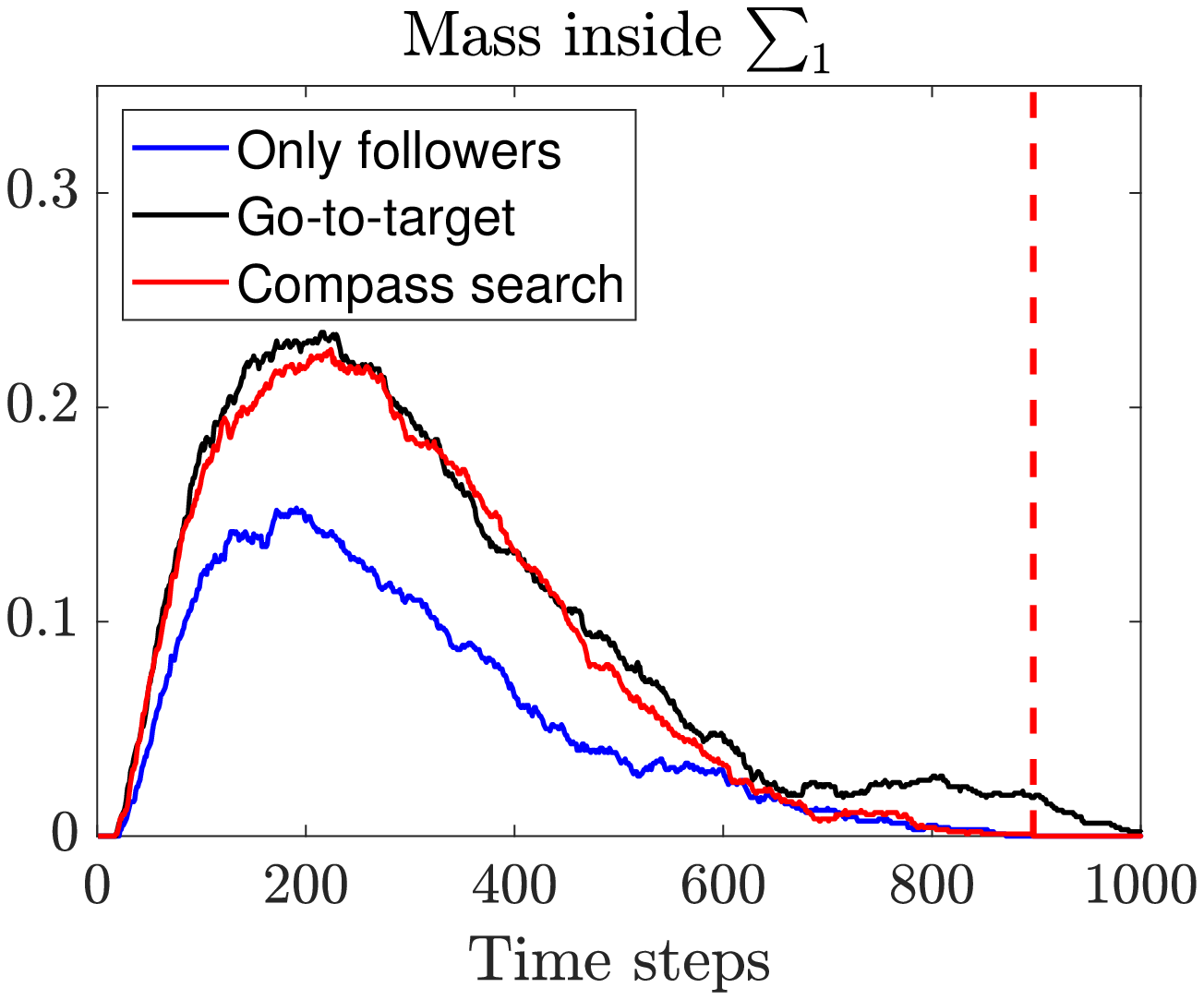}
	\includegraphics[width=0.328\linewidth]{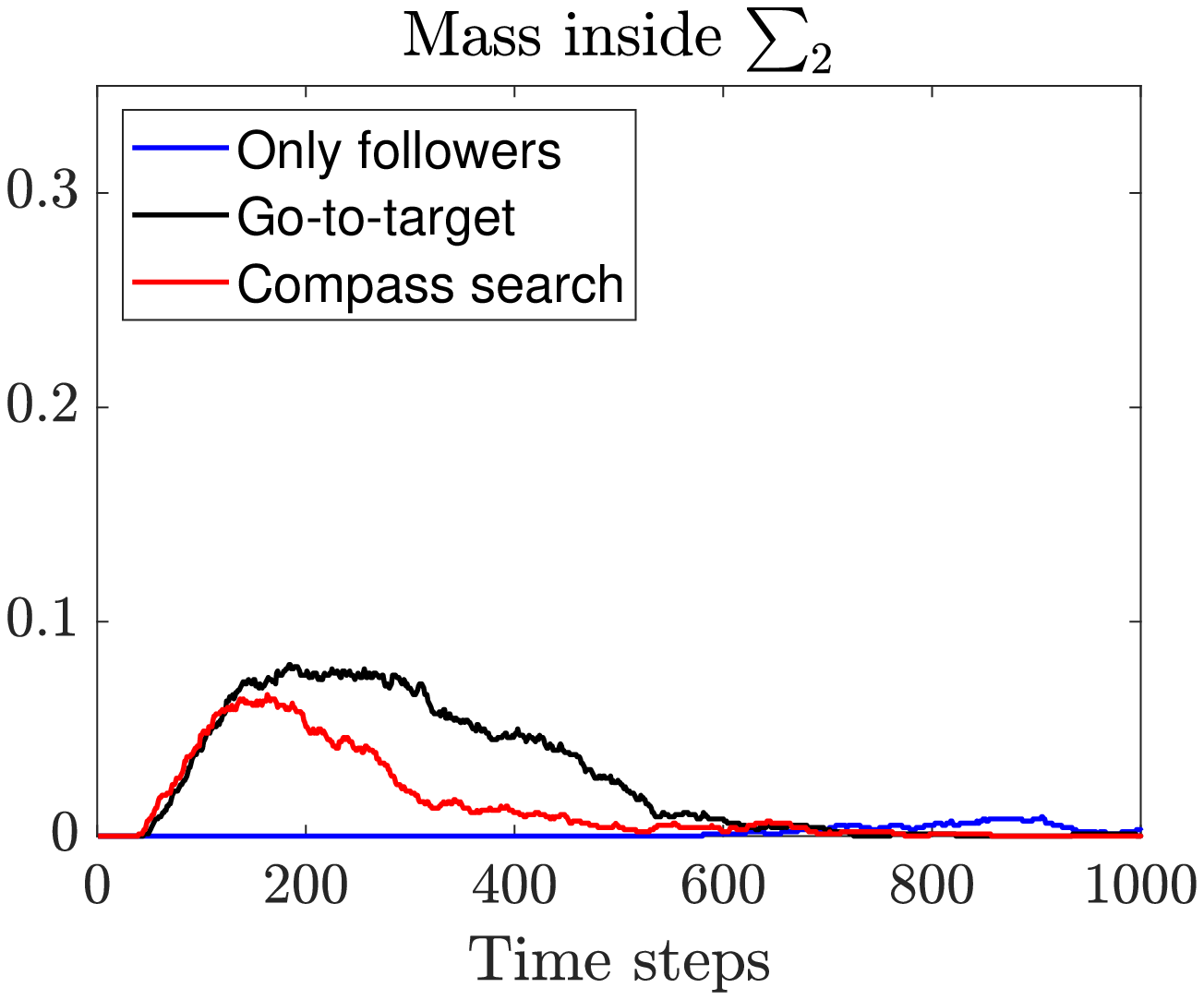}
	\includegraphics[width=0.328\linewidth]{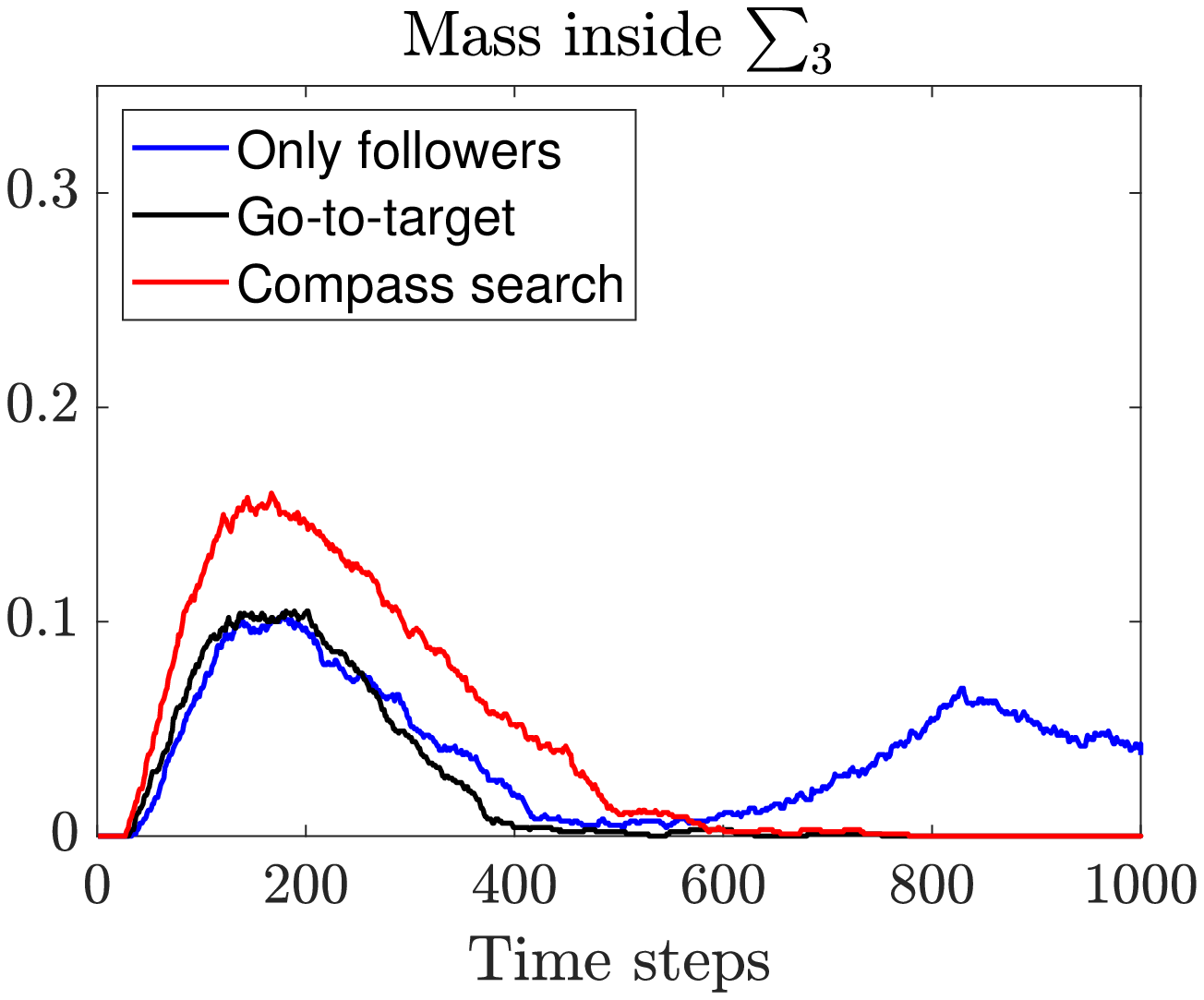}
	\caption{{\em Test 1a}. Mesoscopic case: minimum time evacuation with multiple exits. Evacuated mass (first row), occupancy of the visibility area $\Sigma_1$ (second row, left), $\Sigma_2$ (second row, centre) and $\Sigma_3$ (second row, right) as a function of time for go-to-target and optimal compass search strategies. The dot line denotes the time step in which the whole mass is evacuated with the optimal compass search strategy.}
	\label{fig:test1nvl_time_meso_3}
\end{figure}

Finally in Figure \ref{fig:test1nvl_time_meso_3} we summarize the  results showing the evacuated mass as the cumulative distribution of agents who reached the exit, and the occupancy of the visibility areas in terms of total mass percentage for the various exits.  Dashed red line indicates time of complete evacuation.

\begin{figure}[h!]
	\centering
	\includegraphics[width=0.485\linewidth]{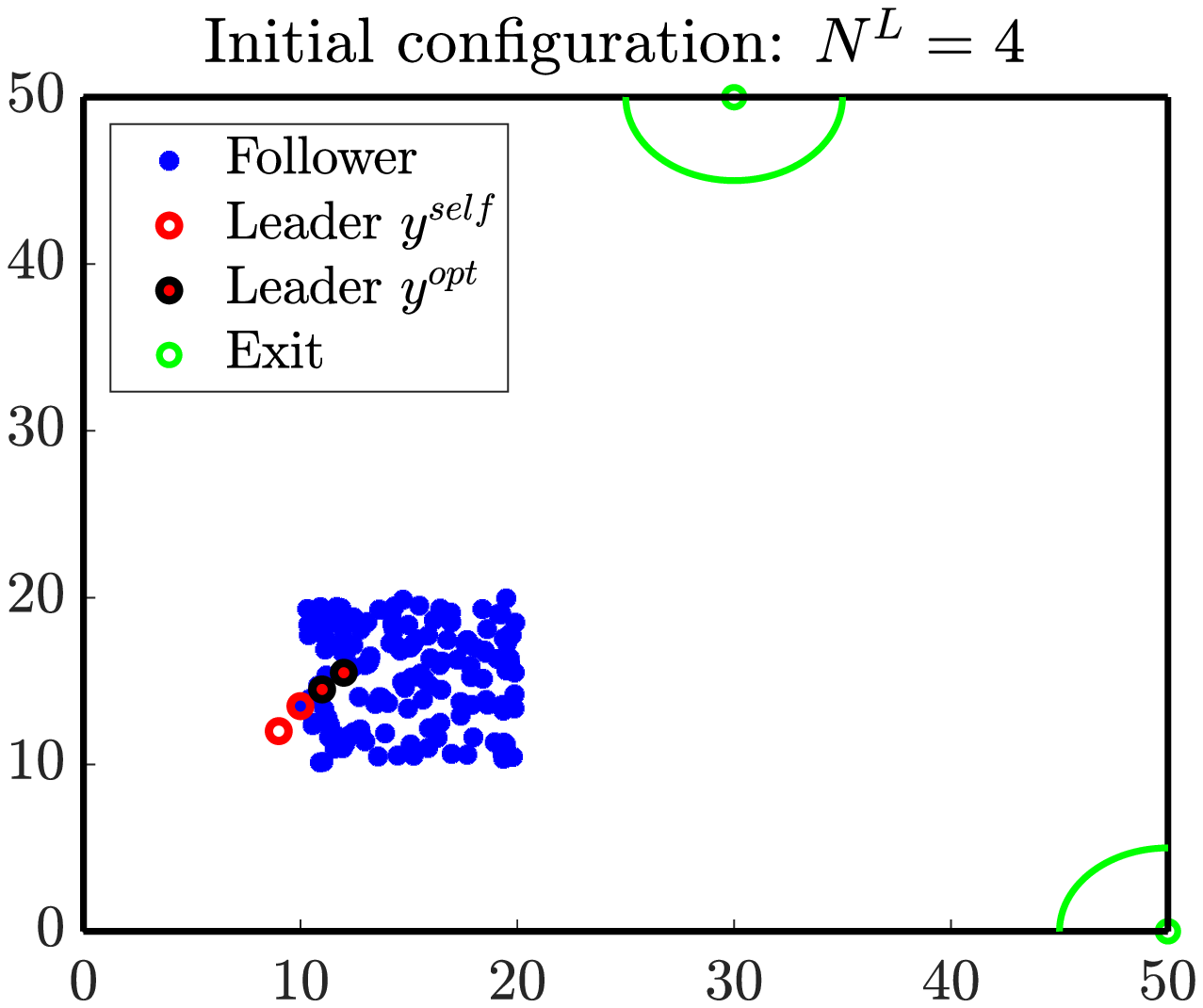}
	\includegraphics[width=0.485\linewidth]{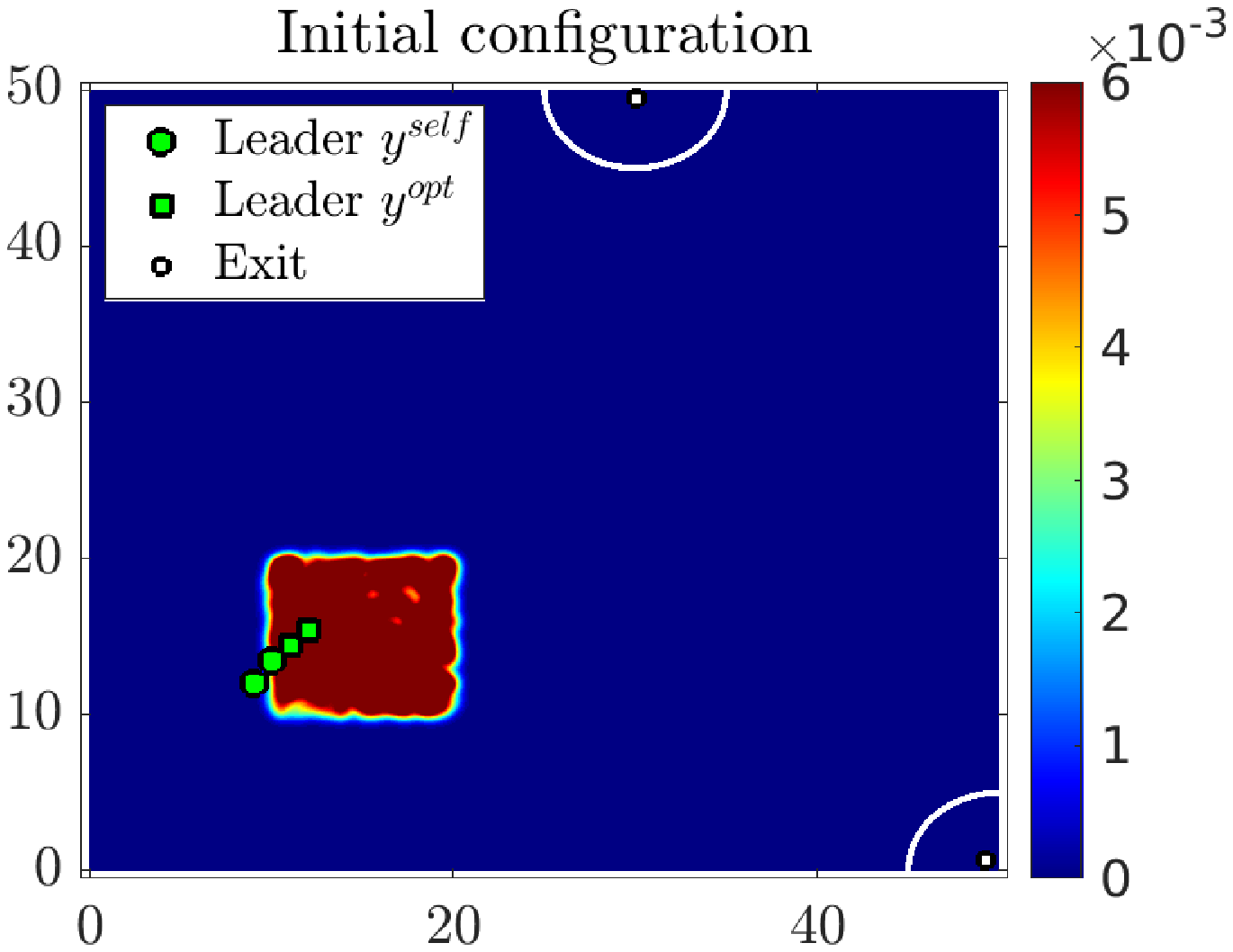}
	\caption{{\em Test 1b.} Minimum time evacuation with multiple exits and obstacles, initial configuration for microscopic and mesoscopic case.}
	\label{fig:test1_2exits_initconf}
\end{figure}

\paragraph{\bf Setting b) Two exits in a closed environment.}
Assume now to have a room with walls that contains two exits, $x_1^\tau=(50,0)$ and $x_2^\tau=(30,50)$. Followers are uniformely distributed in $[0,10]\times [0,10]$. Assume that initially two unaware leaders $y^{self}$ move towards exit $x_1^\tau$ with selfish strategy, i.e. $\beta =0$ in \eqref{eq:gototarget_beta}.
Hence the goal is to minimize the total evacuation time as reported in \eqref{eq:minevactime} introducing two additional leaders $y^{opt}$ moving towards exit $x_2^\tau$, for this two leaders we choose the parameter $\beta = 0.6$ in \eqref{eq:gototarget_beta}.
The target position is $\Xi_k(t)= x_e^\tau \ \forall t$ and for every leader $k$.
Figure \ref{fig:test1_2exits_initconf} shows the initial configuration in the microscopic and mesoscopic case, and with an initial position of $N^L=4$  unaware and aware leaders.

\begin{figure}[h!]
 	\centering
 	\includegraphics[width=0.328\linewidth]{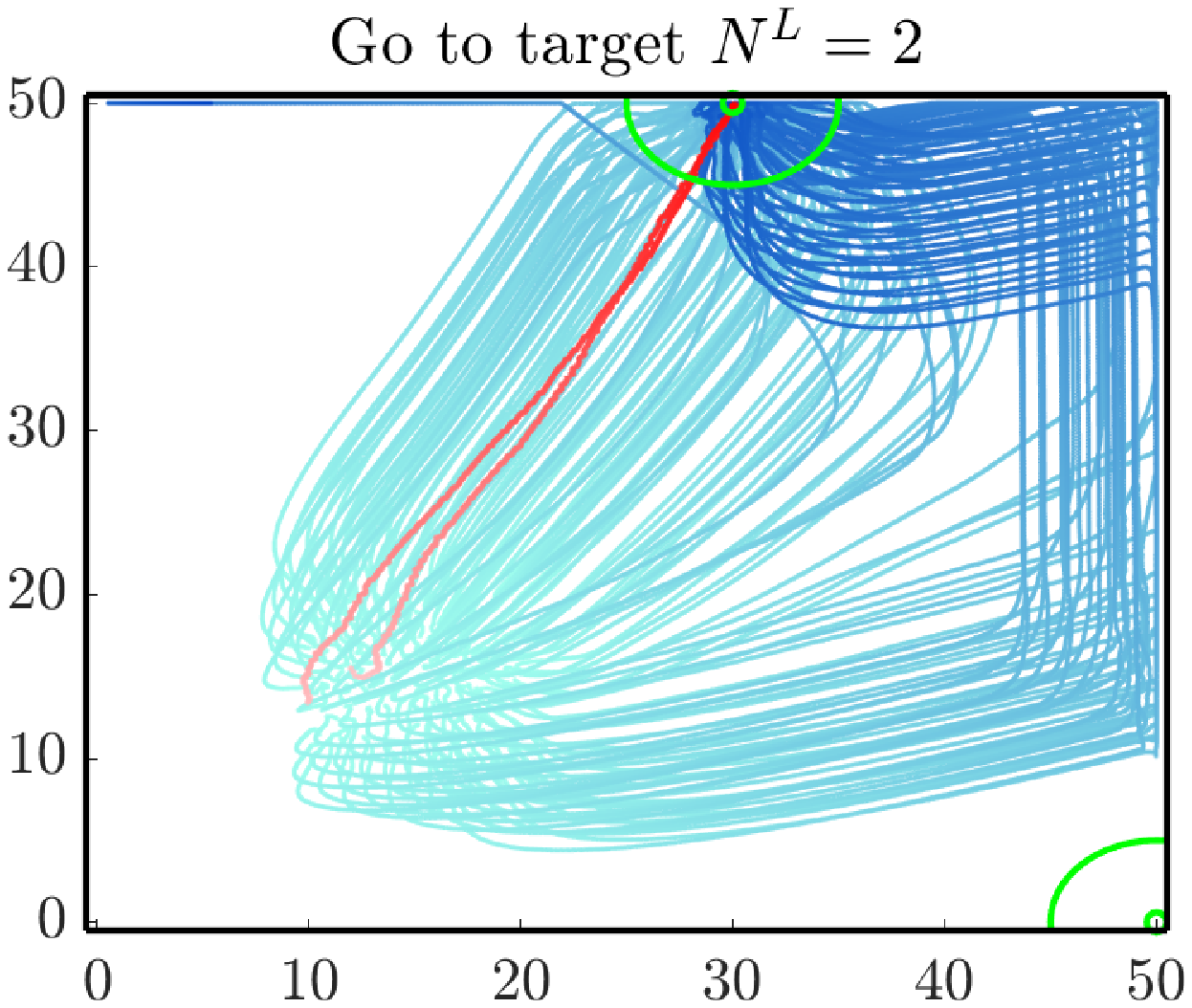}
 	\includegraphics[width=0.328\linewidth]{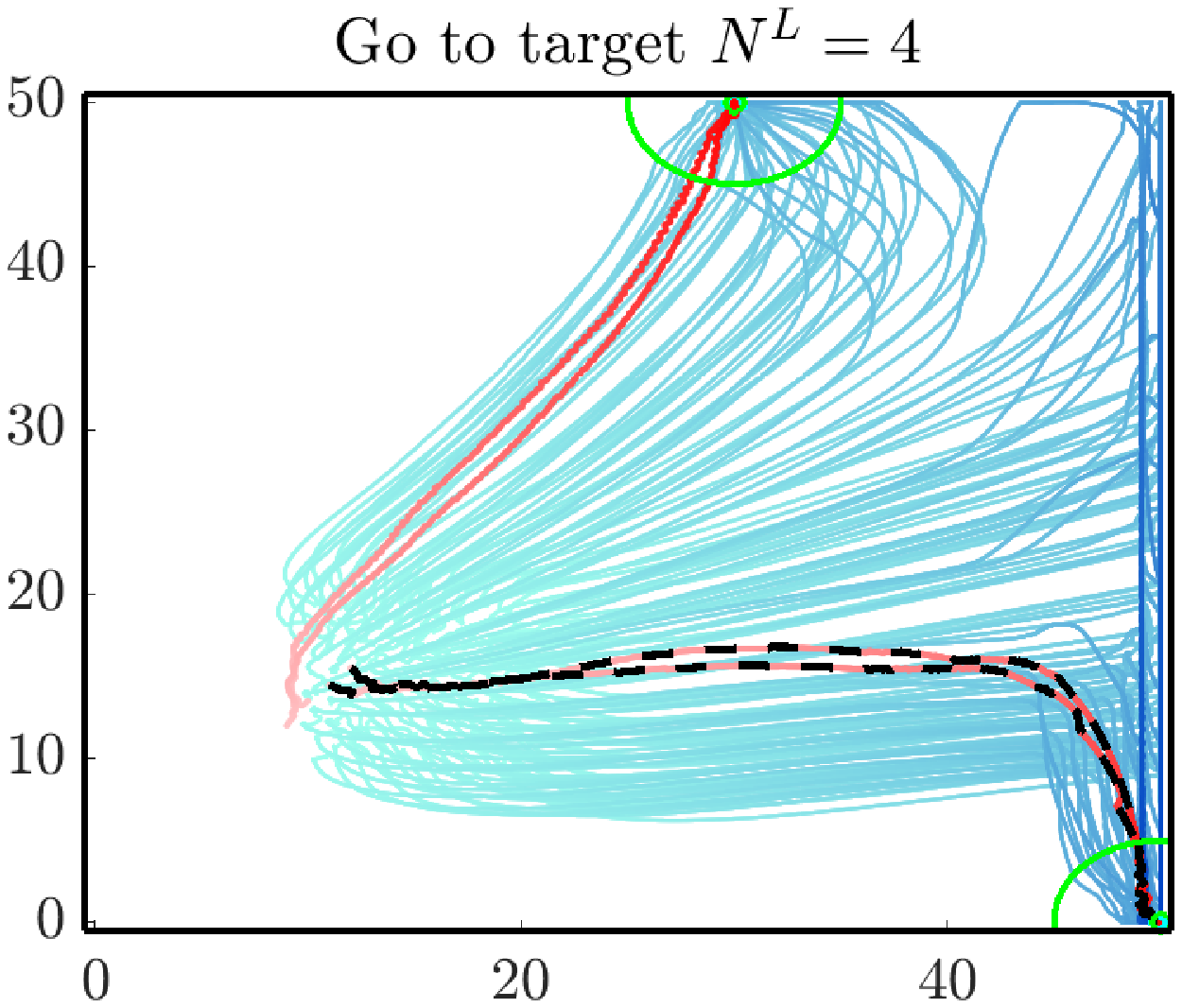}
 	 	\includegraphics[width=0.328\linewidth]{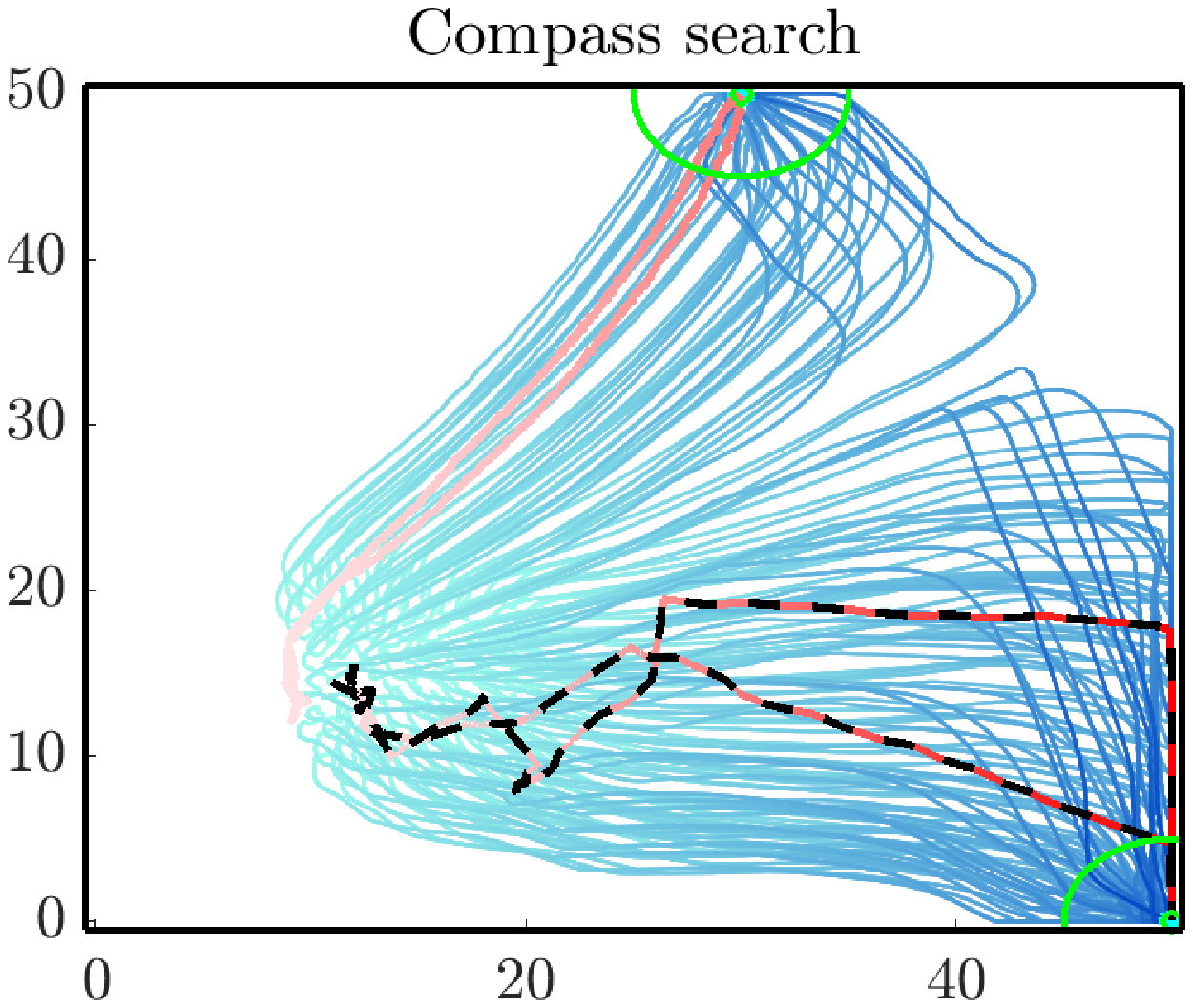}
 	\caption{{\em Test 1b.} Microscopic case: minimum time evacuation with multiple exits and obstacles. Go-to-target $N^L=2$ (left), go-to-target $N^L=4$ (centre), optimal compass search (right).  }
	\label{fig:test1_2exits_micro_1}
 \end{figure}

\begin{figure}[h!]
	\centering
 	\includegraphics[width=0.328\linewidth]{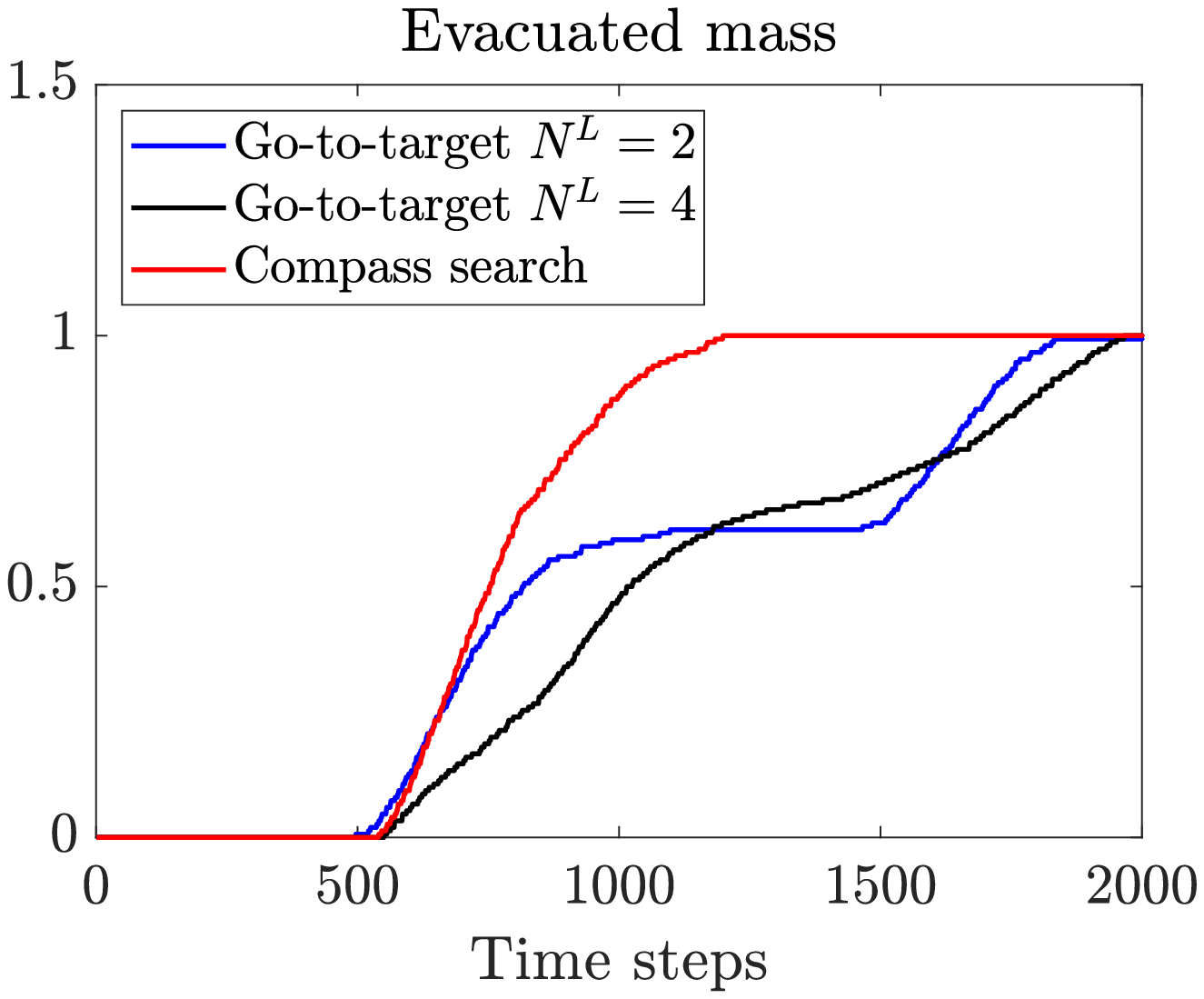}
 	\includegraphics[width=0.328\linewidth]{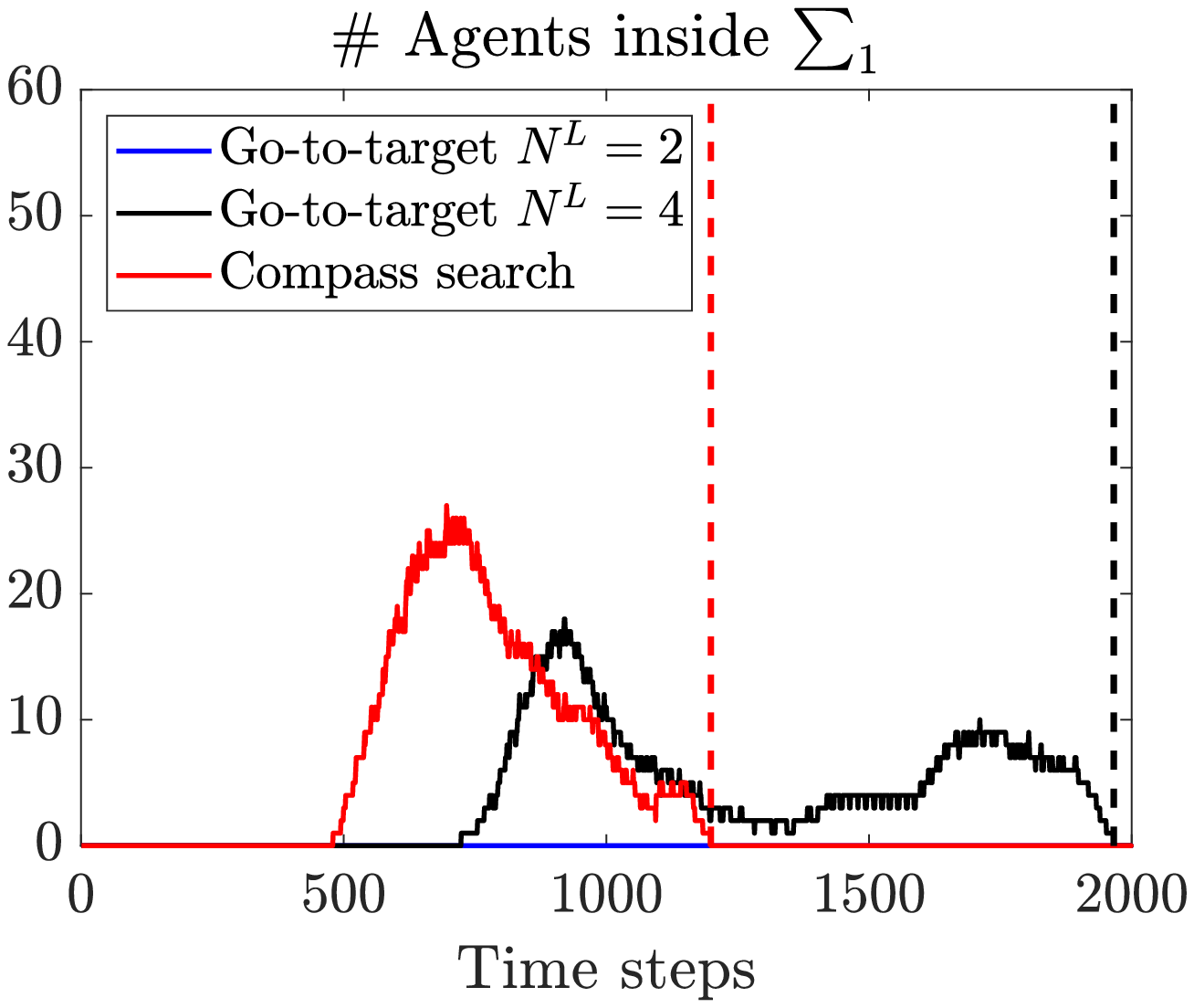}
 	\includegraphics[width=0.328\linewidth]{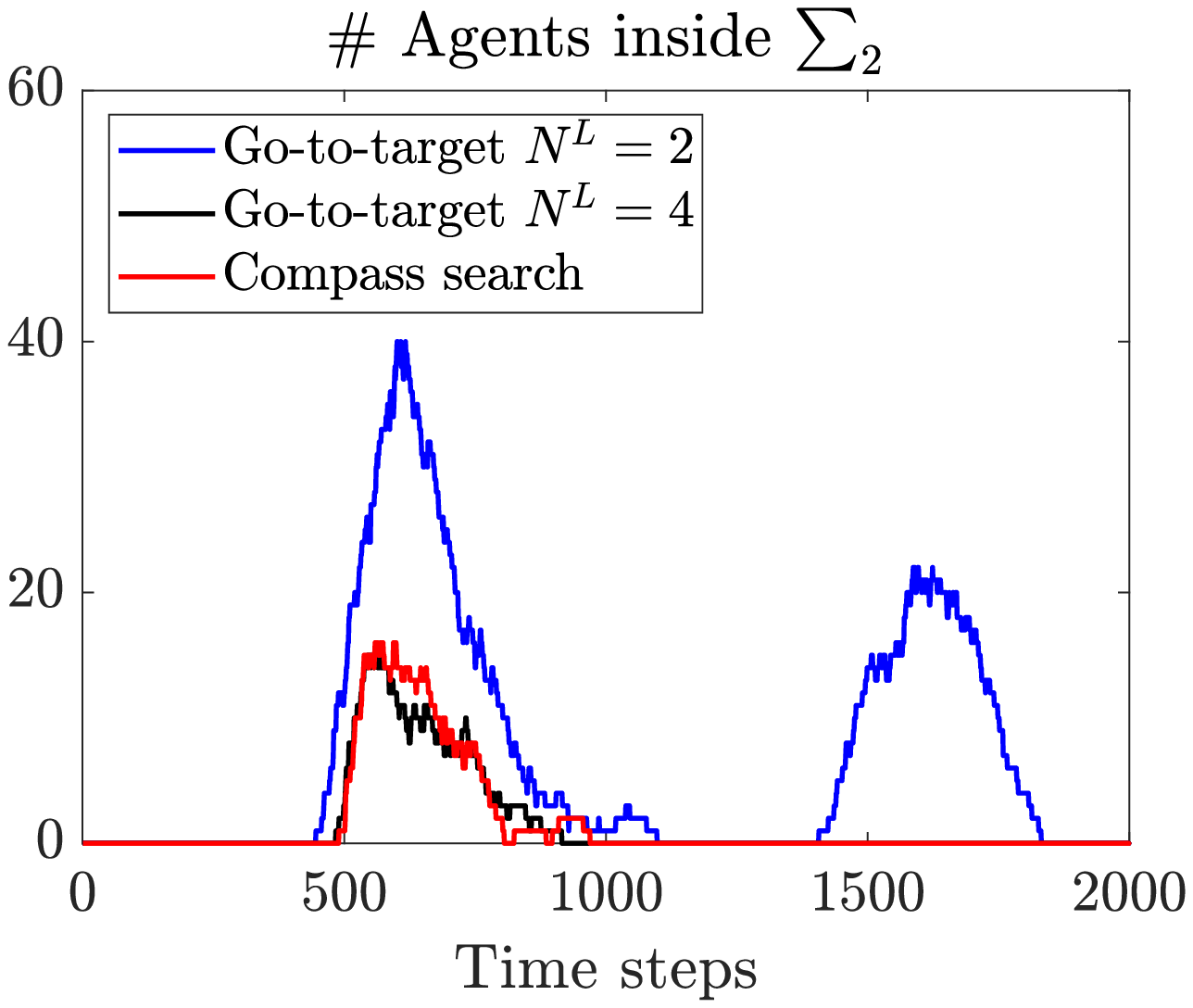}
	\caption{{\em Test 1b.} Microscopic case: minimum time evacuation with multiple exits and obstacles. Evacuated mass (left), occupancy of the visibility area $\Sigma_1$ (centre) and $\Sigma_2$ (right) as a function of time for go-to-target and optimal compass search strategies. The black and red dot lines denote the time step in which the whole mass is evacuated with the go-to-target $(N^L=4)$ and optimal compass search strategy, respectively.}
	\label{fig:test1_2exits_micro_2}
\end{figure}

{\em Microscopic case.}
In Figure \ref{fig:test1_2exits_micro_1} we report the crowd's evolution in various scenarios: left plot shows the trajectories where only unaware leaders are present, in this case, the whole crowd reaches the exit $x_1^\tau$; central and right plots show the influence of two aware leaders moving to $x_2^\tau$ respectively with fixed and optimized strategies. 
Unaware leaders influence the whole crowd to move towards the exit, however generating overcrowding at $x_1^\tau$ and leaving some agents getting lost. 
Introducing two aware leaders with fixed strategies the whole mass is evacuated in $1966$ time steps, with optimized strategies evacuation time is further reduced to $1199$ time steps. In these last cases, the mass is split between the two exits and hence overcrowding phenomena are reduced. In Table \ref{tab_2exits_micro} the total evacuation time and the corresponding evacuated mass for the three scenarios are reported, where we indicate that optimized strategy is obtained after 50 iterations of compass search. Finally, in Figure \ref{fig:test1_2exits_micro_2} we report the occupancy of the visibility areas and the cumulative distribution of the mass evacuate as a function of time for the various scenarios.

 \begin{table}
 	\caption{{\em Test 1b.} Performance of leader strategies over microscopic dynamics.}
 	\label{tab_2exits_micro}
 	\centering
 	\begin{tabular}{ccccc}
 		& go-to-target $N^L=2$  & go-to-target $N^L=4$ & CS (50 it)\\
 		\cmidrule(r){2-2}\cmidrule(r){3-3}\cmidrule(r){4-4}
 		Evacuation time (time steps)  & >2000 & 1966 & 1199\\
 		Evacuated mass (percentage)  & 99\% & 100\% & 100\%\\
 		\hline
 	\end{tabular}
 \end{table}

\begin{figure}[h]
	\centering
	\includegraphics[width=0.328\linewidth]{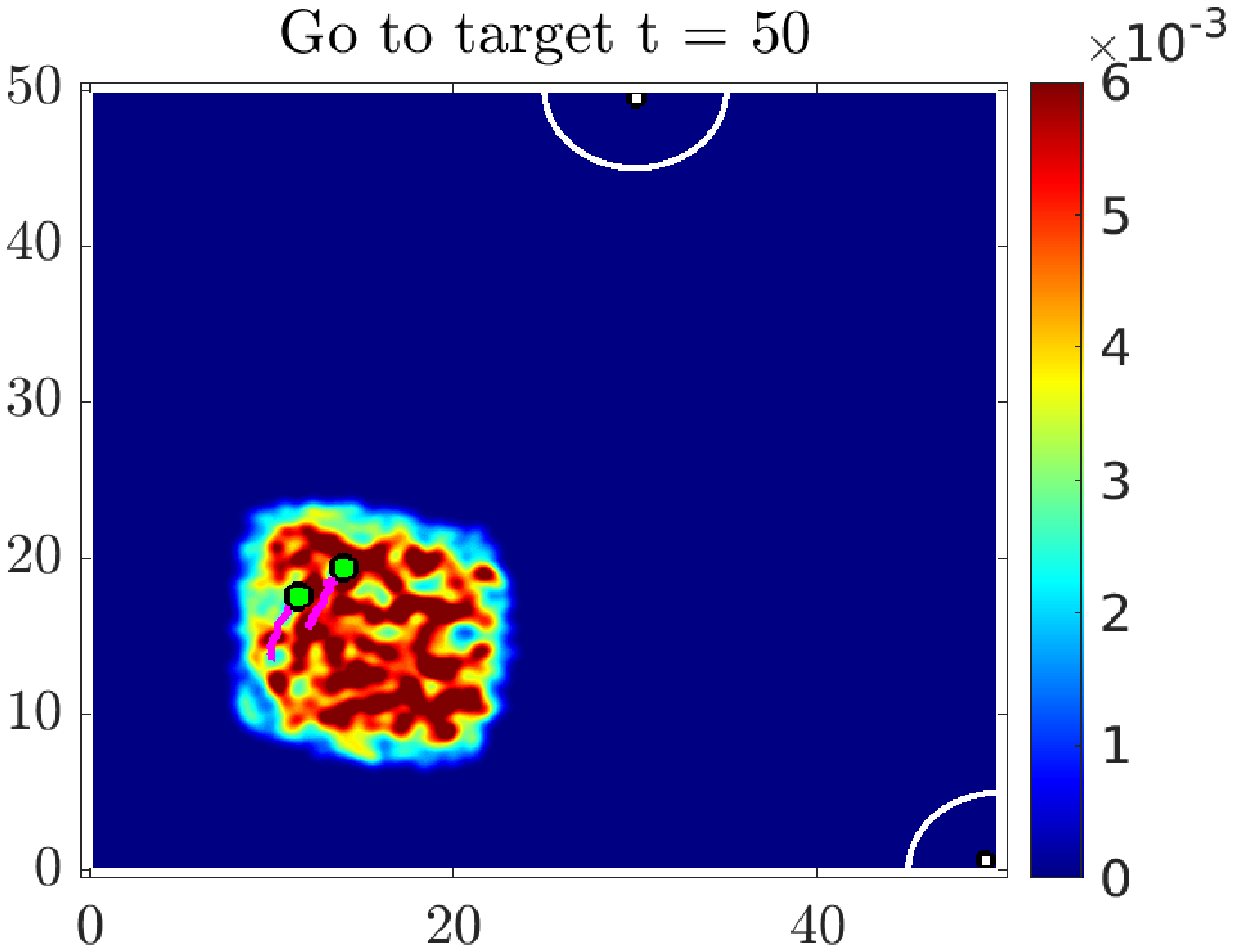}
	\includegraphics[width=0.328\linewidth]{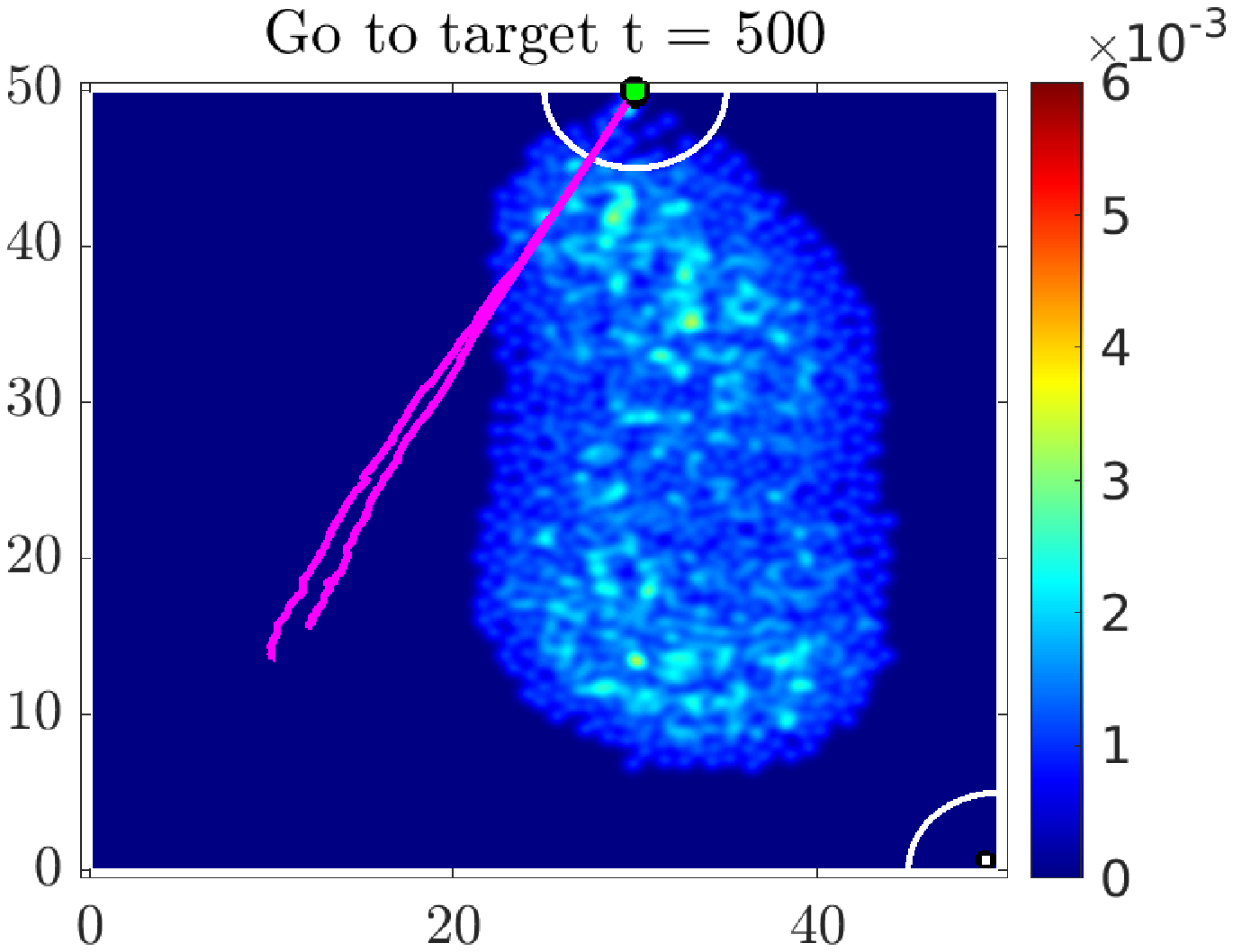}
	\includegraphics[width=0.328\linewidth]{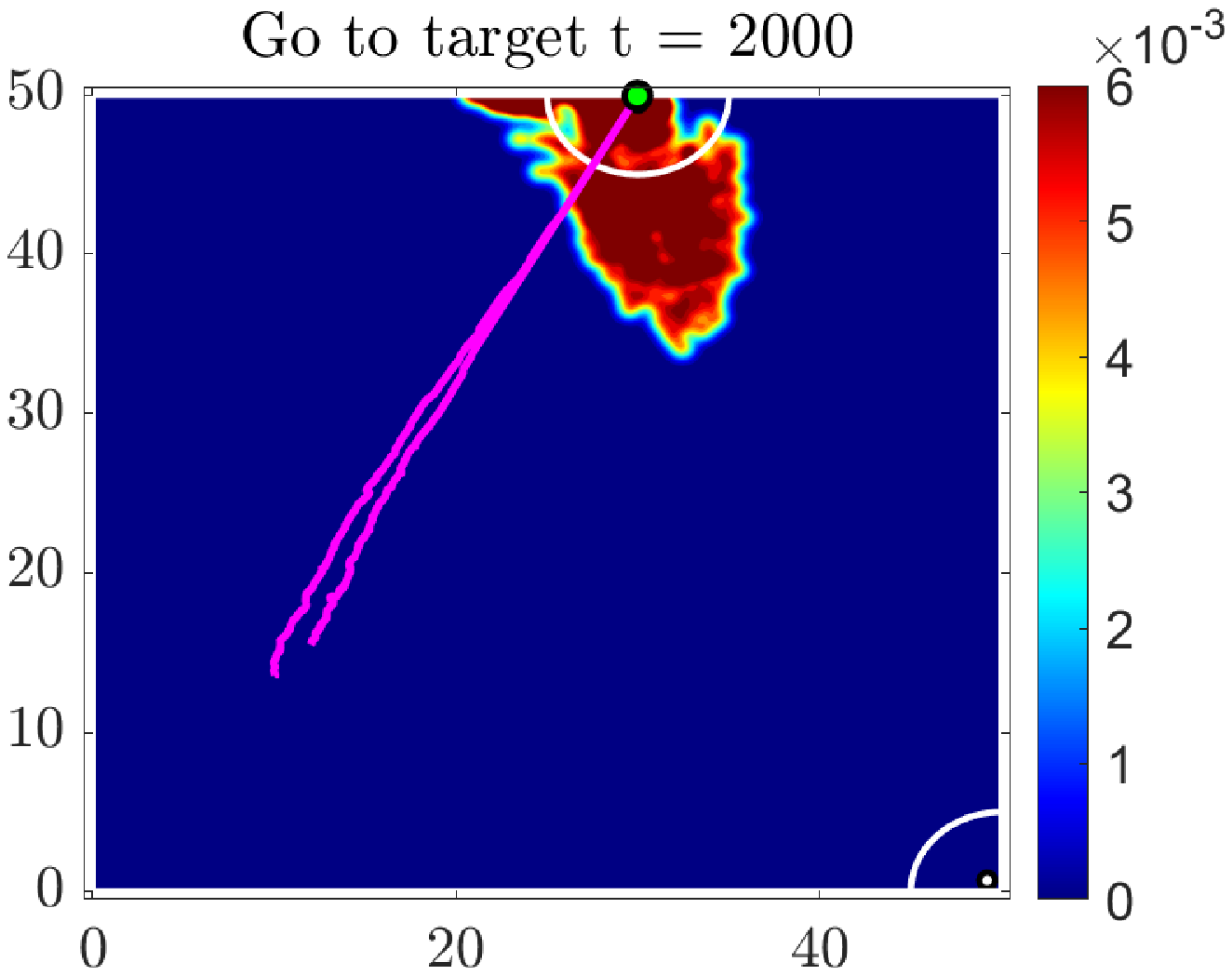}
	\\
    \includegraphics[width=0.328\linewidth]{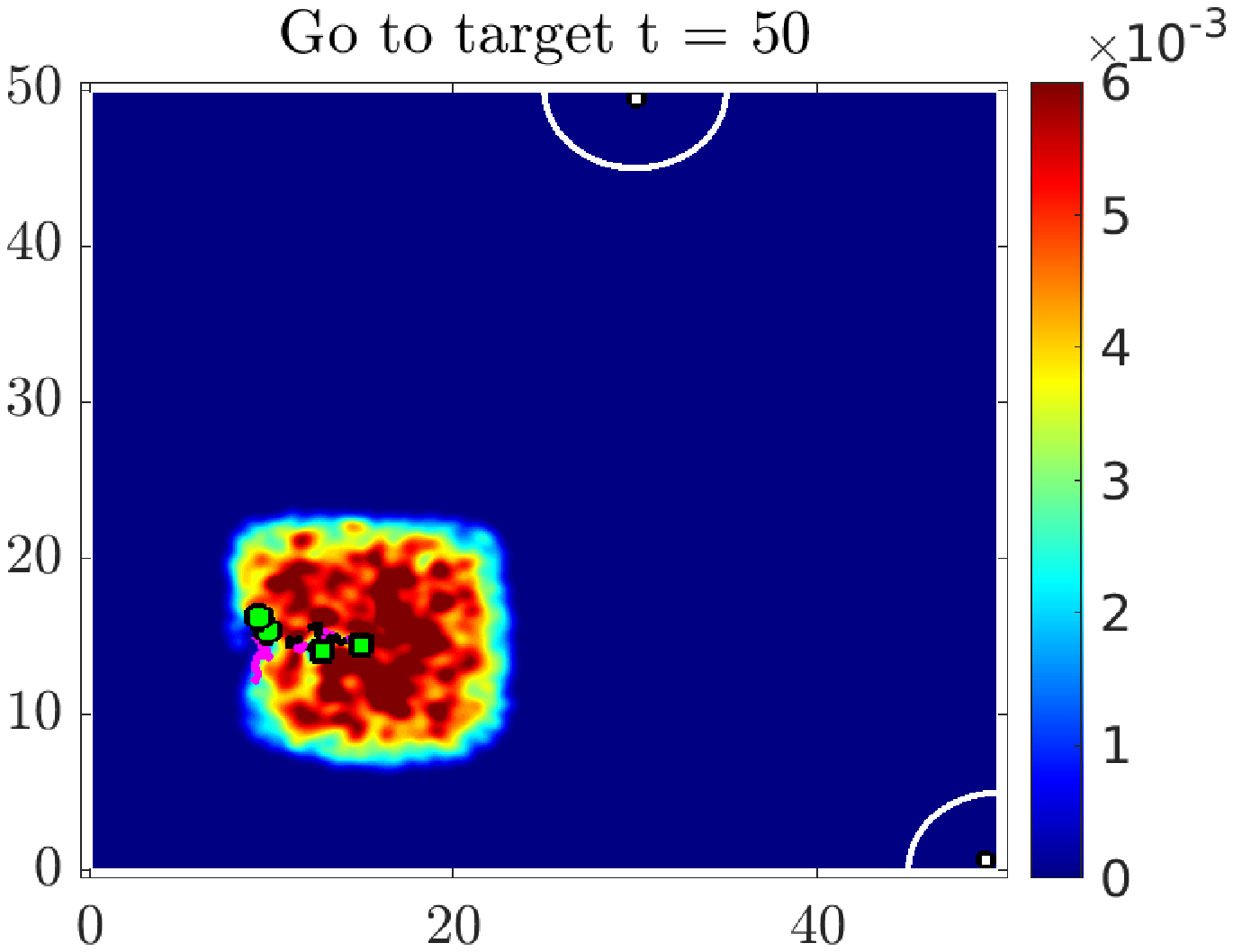}
	\includegraphics[width=0.328\linewidth]{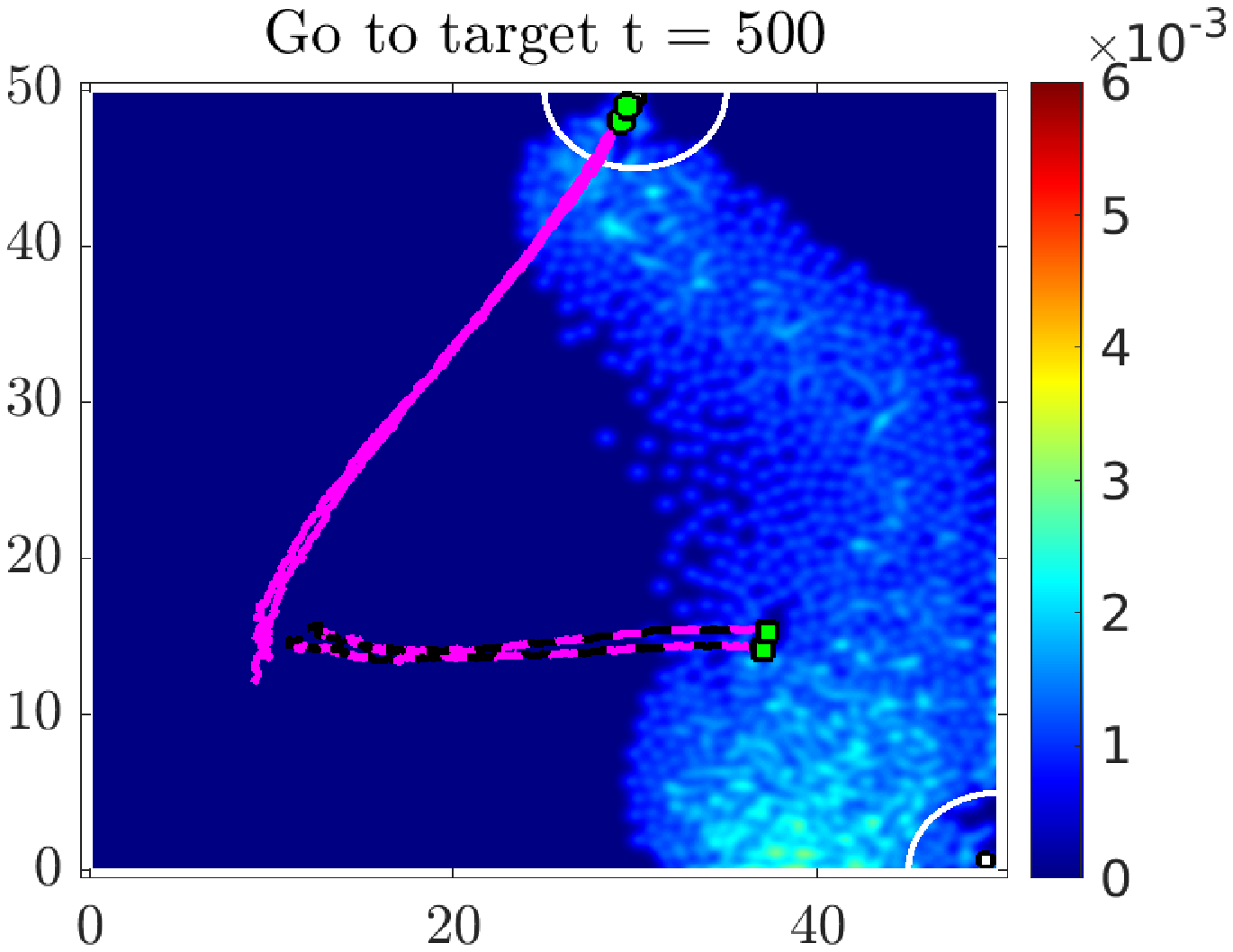}
	\includegraphics[width=0.328\linewidth]{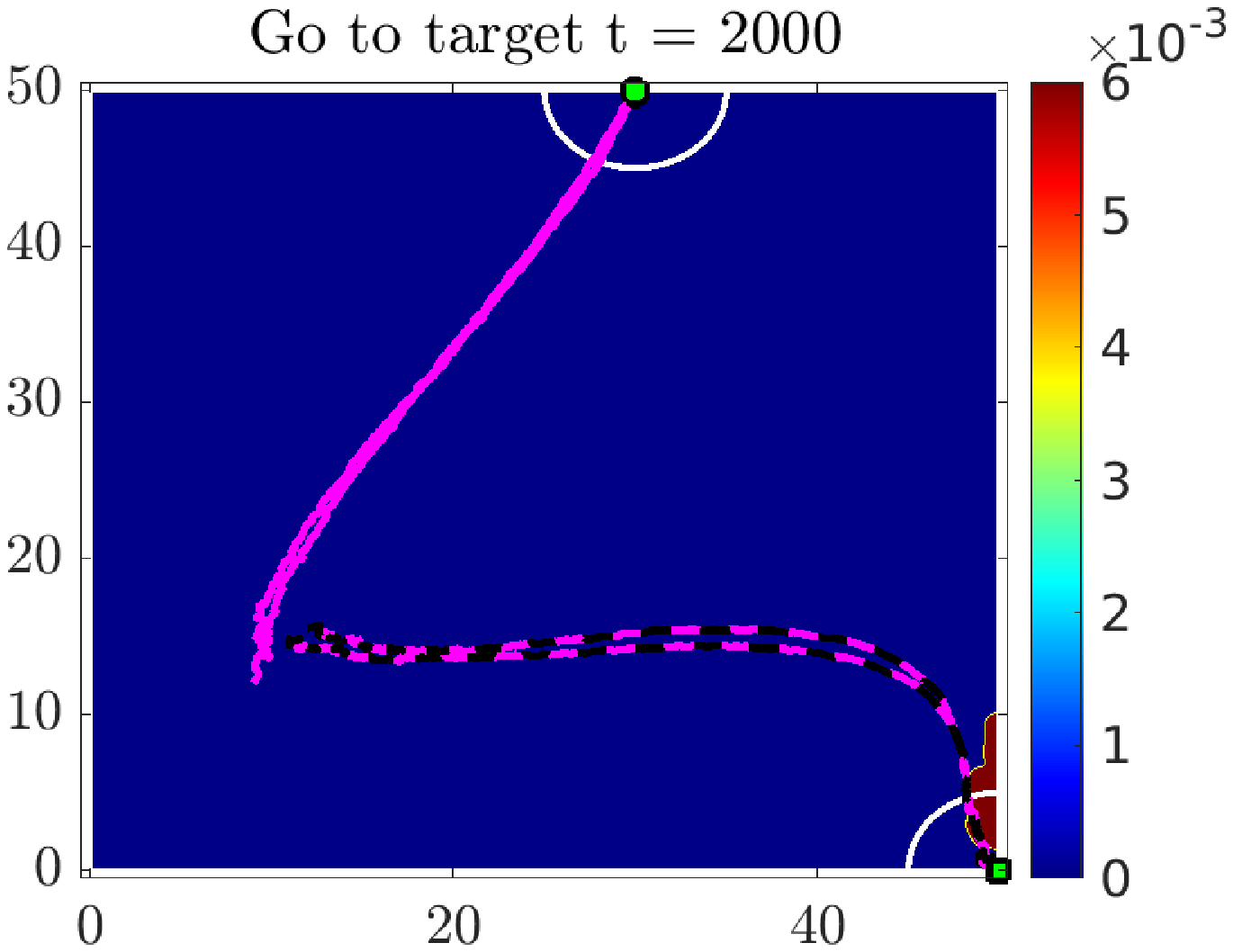}
	\\
	\includegraphics[width=0.328\linewidth]{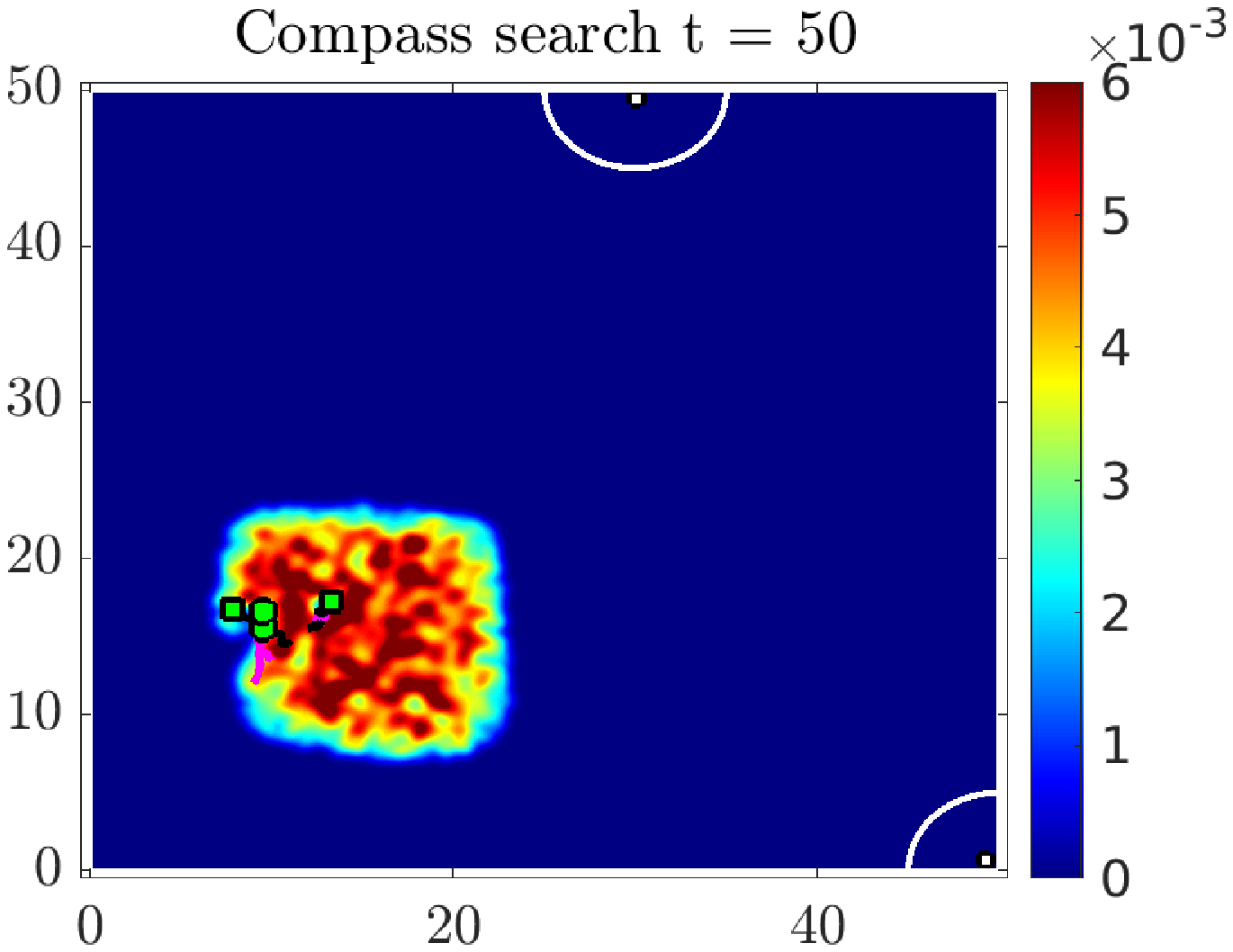}
	\includegraphics[width=0.328\linewidth]{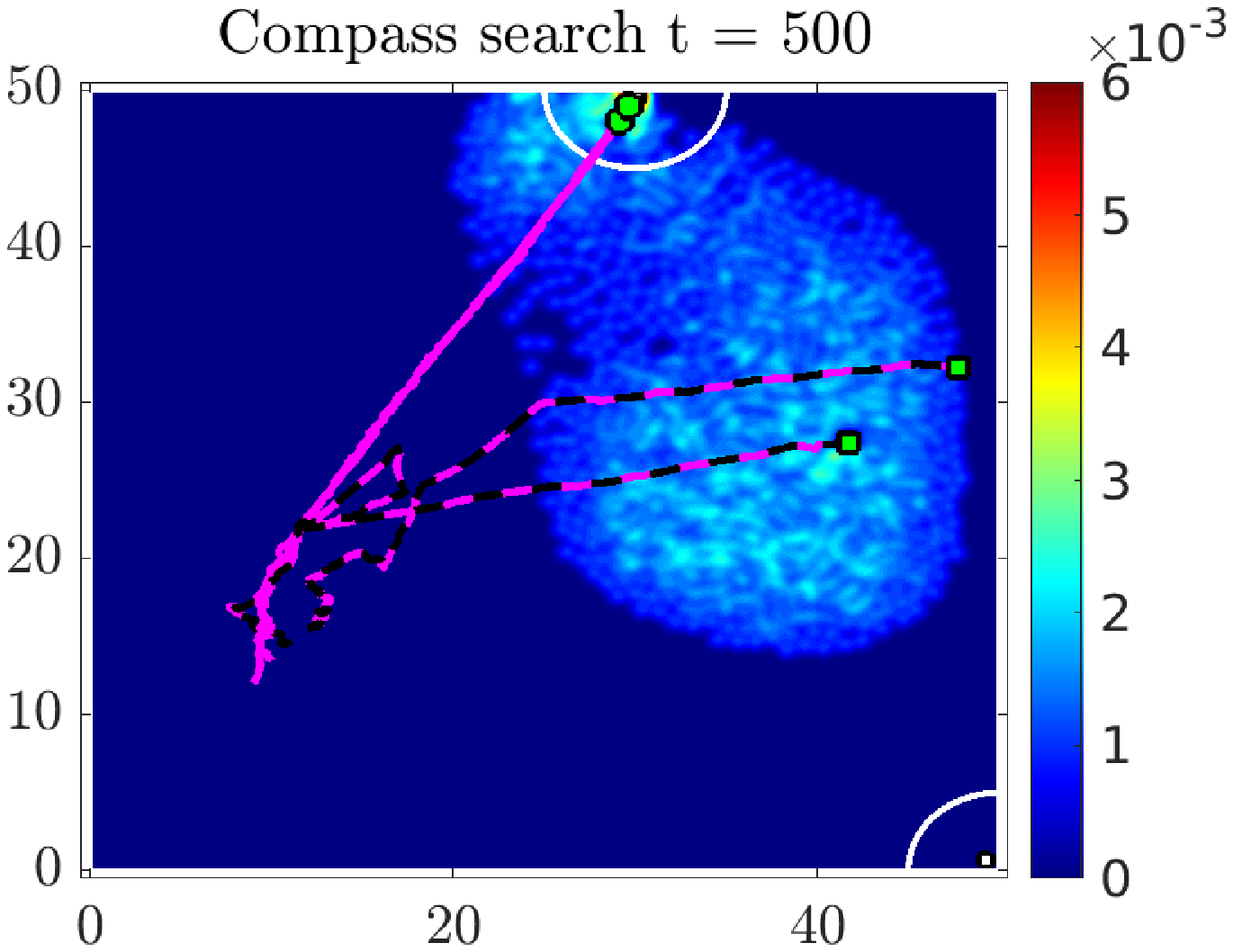}
	\includegraphics[width=0.328\linewidth]{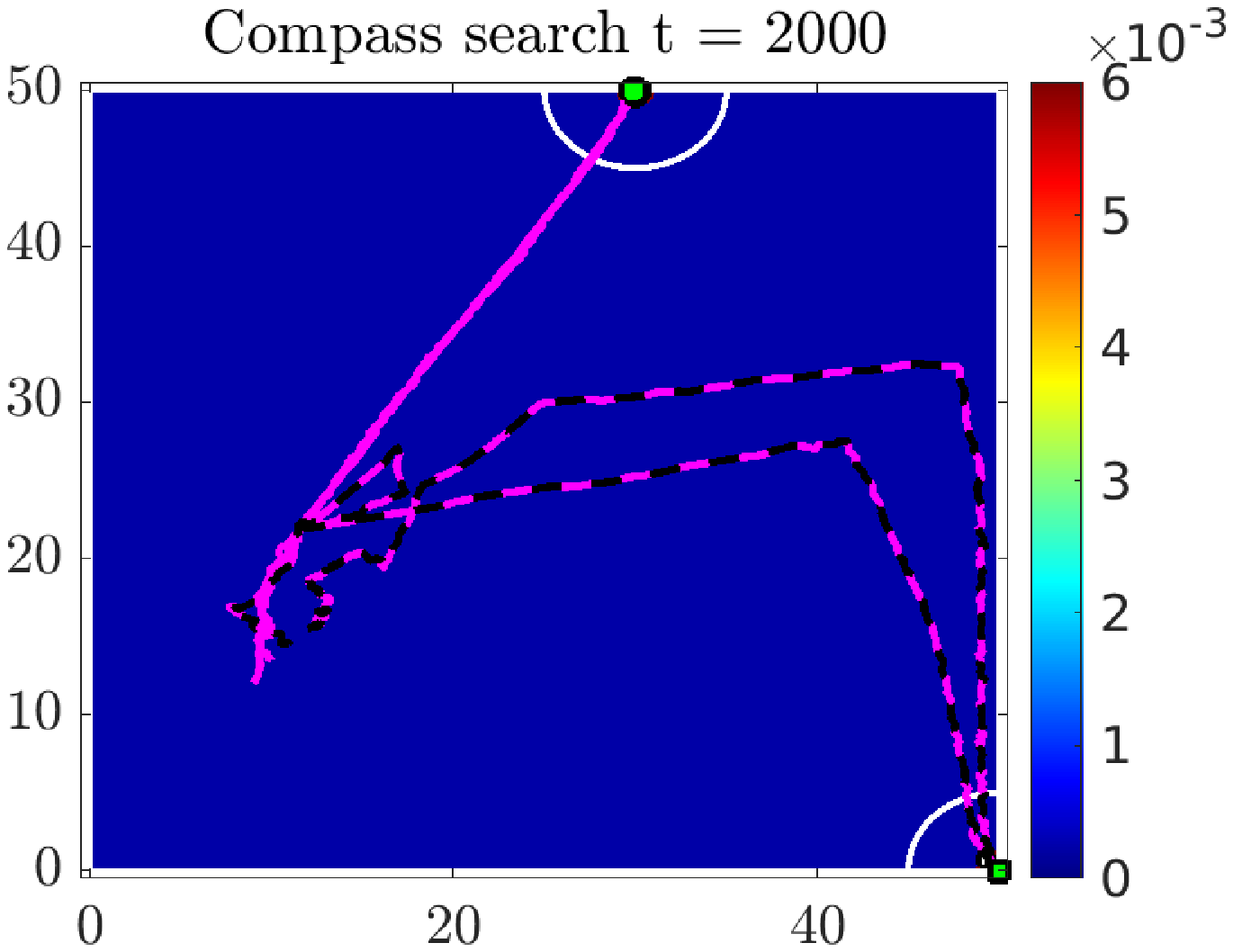}
	\caption{{\em Test 1b.} Mesoscopic case: minimum time evacuation with multiple exits and obstacles. Three snapshots taken at time $t=50$, $t=500$, $t=2000$ with the go-to-target strategy in the case $N^L = 2$ (upper row), $N^L=4$ (central row) and with the optimized compass-search strategy (lower row).}
	\label{fig:test12exits_time_meso_1}
\end{figure}

\begin{figure}[h!]
	\centering
	\includegraphics[width=0.328\linewidth]{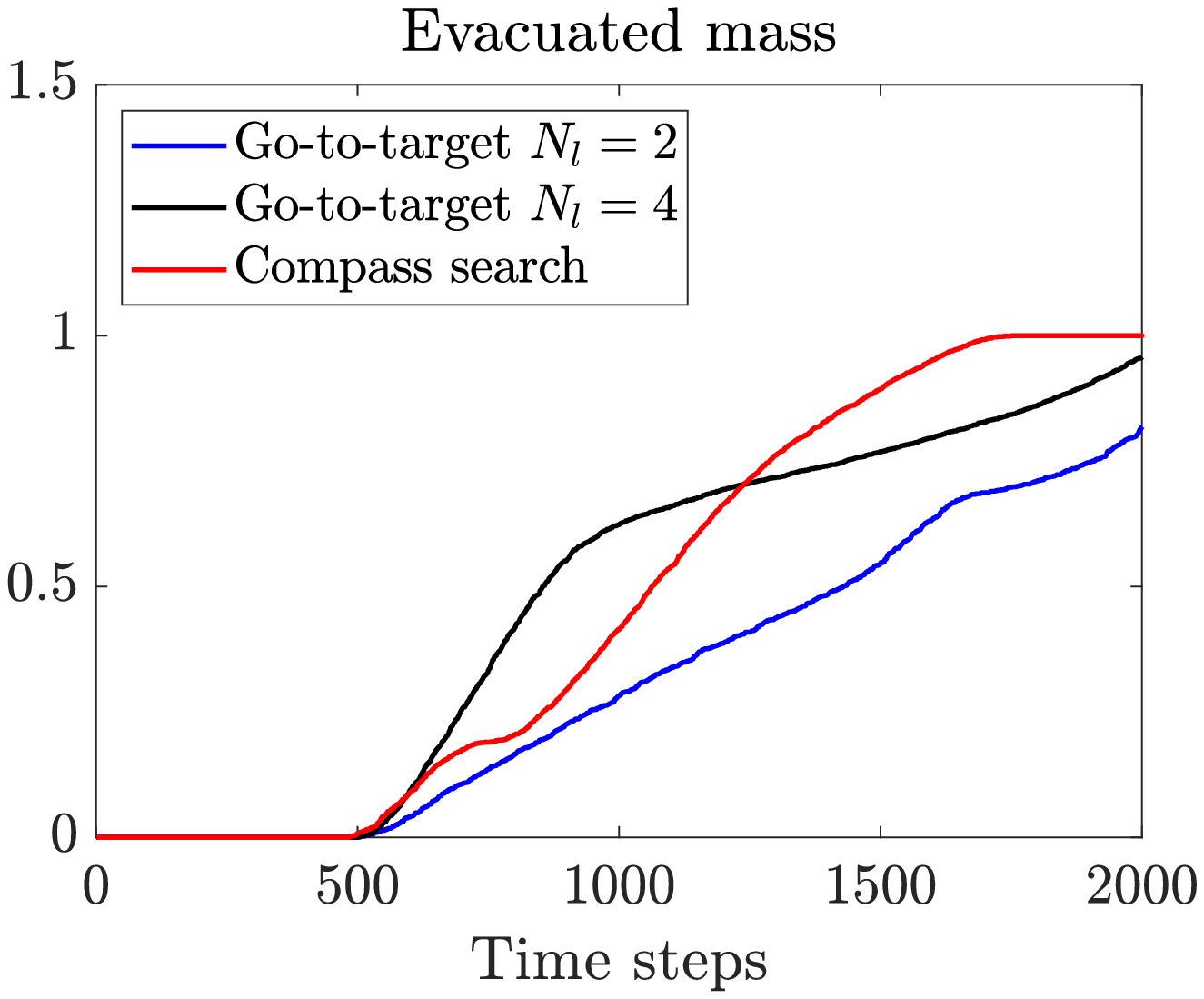}
	\includegraphics[width=0.328\linewidth]{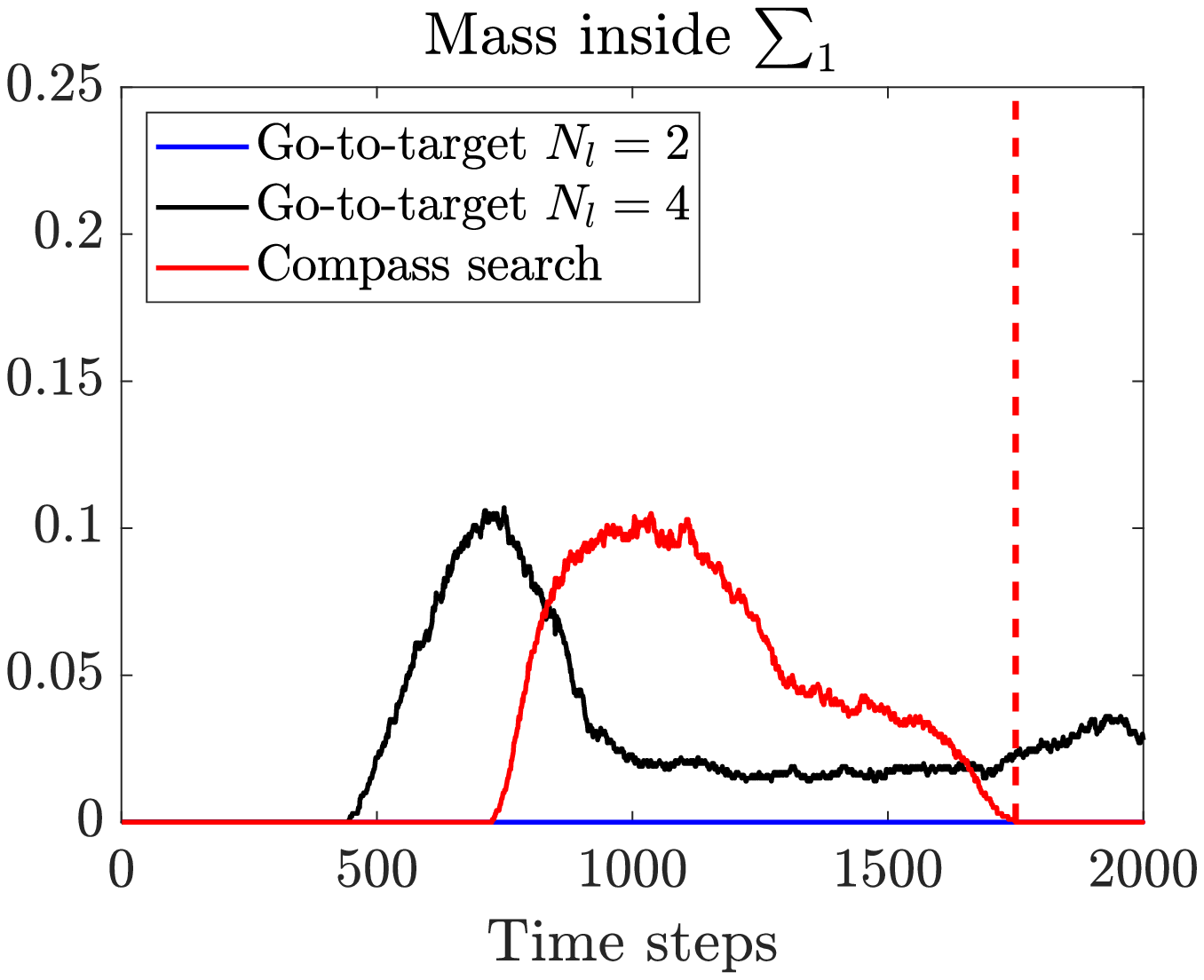}
	\includegraphics[width=0.328\linewidth]{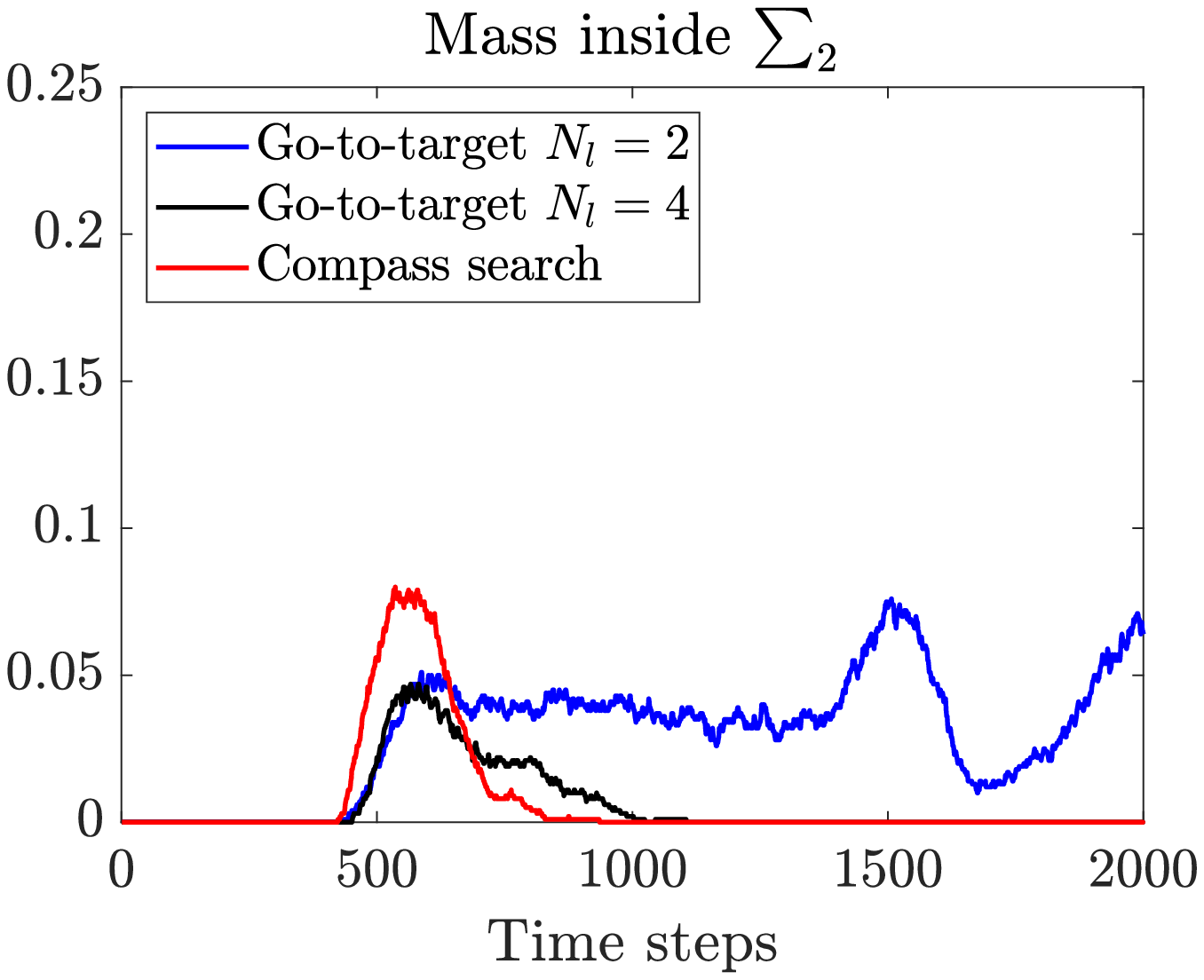}
	\caption{{\em Test 1b.} Mesoscopic case: minimum time evacuation with multiple exits and obstacles. Evacuated mass (left), occupancy of the visibility area $\Sigma_1$ (centre) and $\Sigma_2$ (right) as a function of time for go-to-target and optimal compass search strategies. The red dot line denotes the time step in which the whole mass is evacuated with the optimal compass search strategy.}
	\label{fig:test12exits_time_meso_3}
\end{figure}

{\em Mesoscopic case.} We consider now the mean-field approximation of the microscopic setting. We report in Figure \ref{fig:test12exits_time_meso_1} three snapshots of followers density and trajectories of  leaders, for each scenario. 
In this case, unaware leaders moving selfishly towards exit $x_1^\tau$ are able to influence followers and evacuate $81\%$ at final time, whereas the rest of the mass is congested around the exit.
Introducing two aware leaders with a fixed strategy toward $x_2^\tau$ is not sufficient to reach total evacuation at final time which is and at final time $95\%$ of the mass is evacuated. The bottom row depicts the case with optimized leaders strategies, in this case, the total mass is evacuated at time step $1750$. We summarize the performances of the results in Table \ref{tab_2exits_meso}, and in Figure \ref{fig:test12exits_time_meso_3} we report the occupancy of the visibility areas and the cumulative distribution of mass evacuated as a function of time. 
\begin{table}
	\caption{{\em Test 1b.} Performance of leader strategies over mesoscopic dynamics.}\label{tab_2exits_meso}
\begin{tabular}{ccccc}
& go-to-target $N^L=2$  & go-to-target $N^L=4$ & CS(50 it)\\
 		\cmidrule(r){2-2}\cmidrule(r){3-3}\cmidrule(r){4-4}
Evacuation time (time steps)  & >2000 & >2000& 1750\\
Evacuated mass  & 81\% & 95\% & 100\%\\
\hline
\end{tabular}
\end{table}

\subsubsection{Test 2 : Mass evacuation in presence of obstacles}\label{sec:test2}
\begin{figure}[h!]
	\centering
	\includegraphics[width=0.495\linewidth]{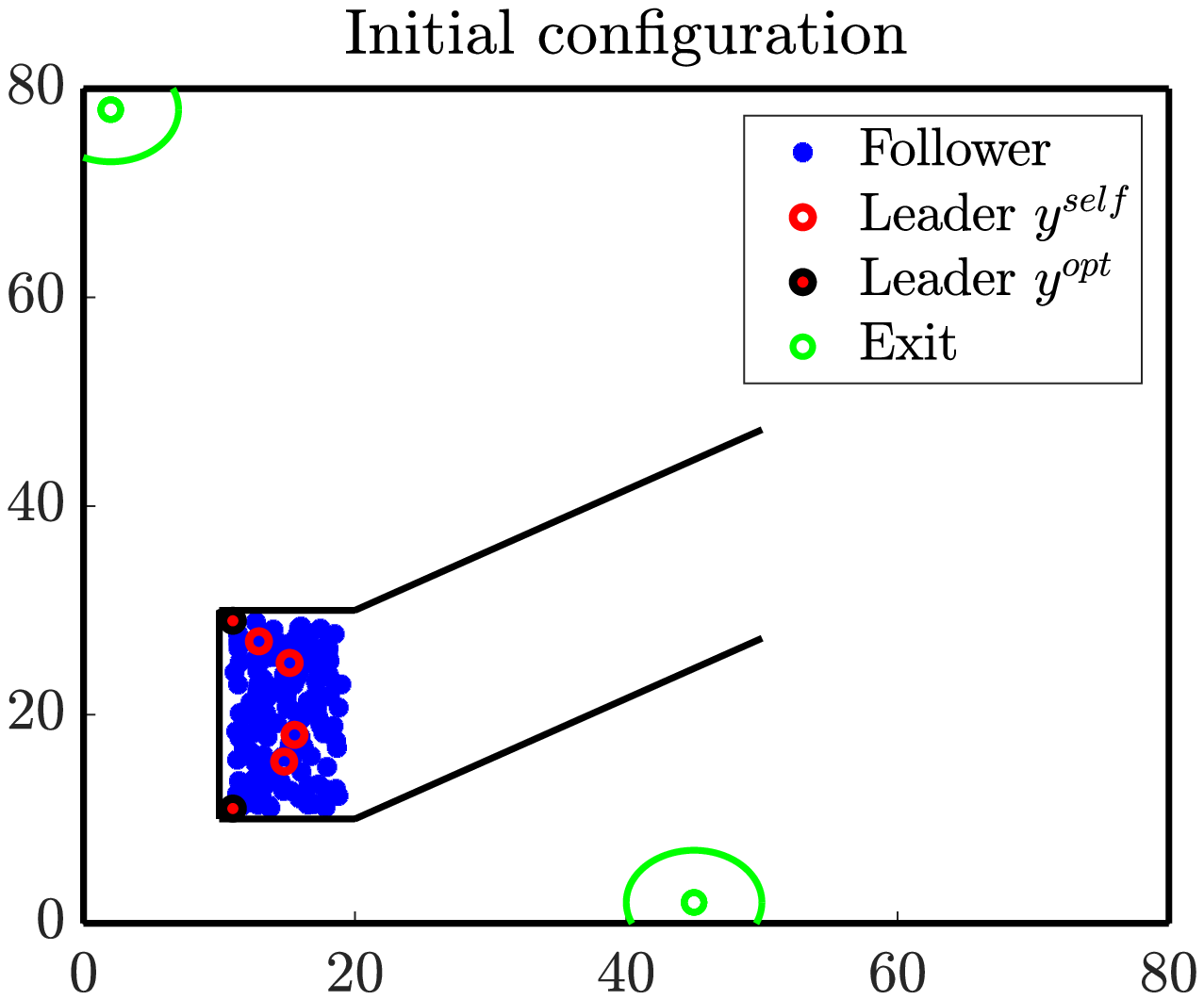}		\includegraphics[width=0.495\linewidth]{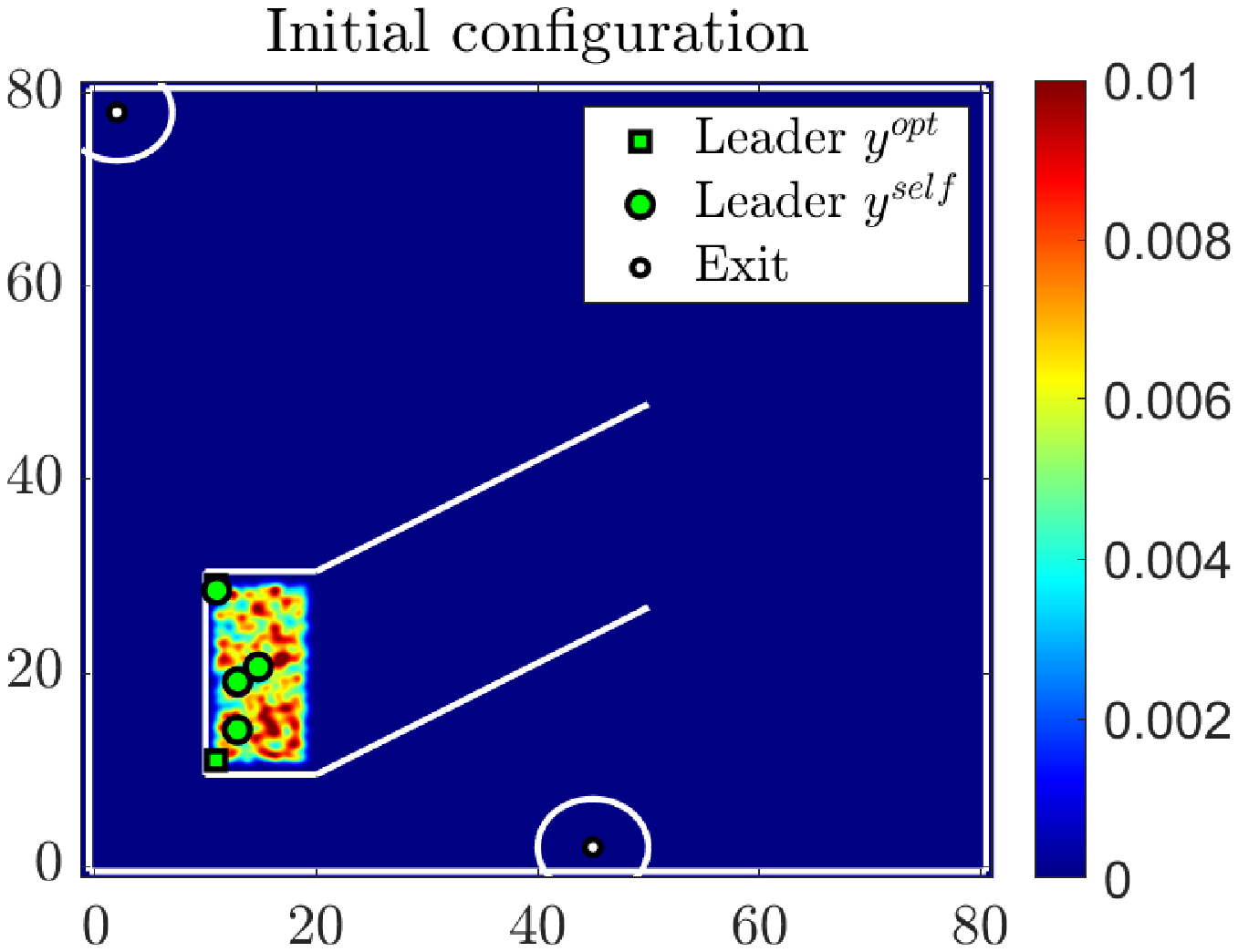}
	\caption{{\em Test 2.} Maximization of mass evacuated in presence of obstacles, initial configuration.}
	\label{fig:test2D_mass_micro_0}
\end{figure}

We consider two rooms, one inside the other, where the internal room is limited by three walls while the external one is bounded by four walls. We assume that walls are nonvisible obstacles, i.e. people can perceive them only by physical contact. This corresponds to an evacuation in case of null visibility (but for the exit points which are still visible from within $\Sigma_1$ and $\Sigma_2$). Consider the case of two exits, $x_1^\tau=(2,78)$ and $x_2^\tau=(45,2)$ positioned in the external room. Figure \ref{fig:test2D_mass_micro_0} provides a description of the initial configuration.  Note that in order to evacuate, people must first leave the inner room, in which they are initially confined, and then search for exits.
Evacuation in presence of obstacles is not always feasible. Instead of minimizing the total evacuation time as in section \ref{sec:test1}, we aim to minimize the total mass inside the domain as reported in \eqref{eq:test2intro} and hence to maximize the total evacuated mass. 

Each leader will move toward one of the exits following a go-to-target, similar to \eqref{eq:gototarget_beta}, and such that it is admissible for the configuration of the obstacles.
We choose $\beta = 1$ for every leader.
The target position is $\Xi_k(t)= x_e^\tau$ for $t>t_*$, while for $t<t_*$ we consider one intermediate point in order to let the leaders to evacuate the inner room.

\begin{figure}[h!]
	\centering
	\includegraphics[width=0.495\linewidth]{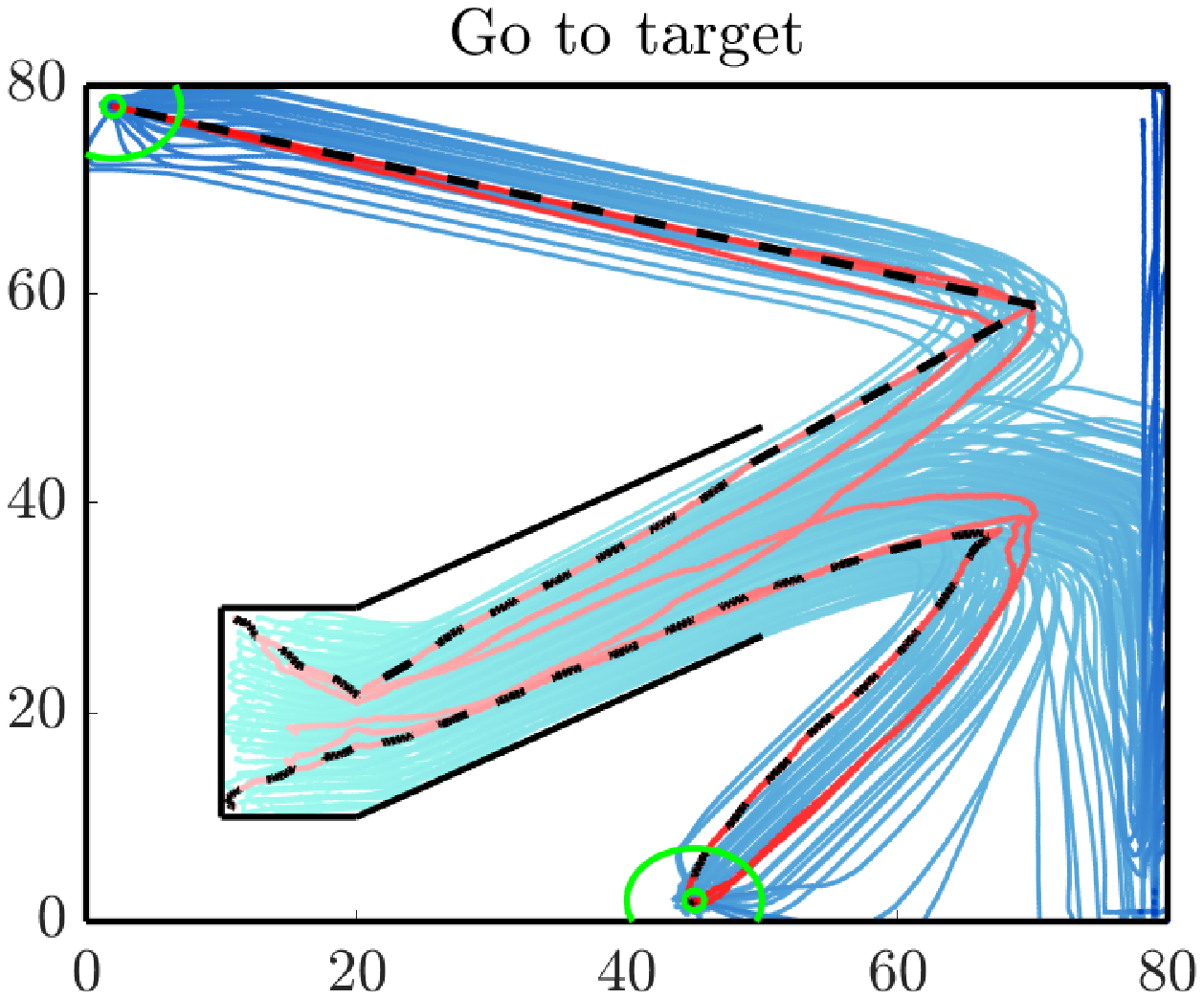}
	\includegraphics[width=0.495\linewidth]{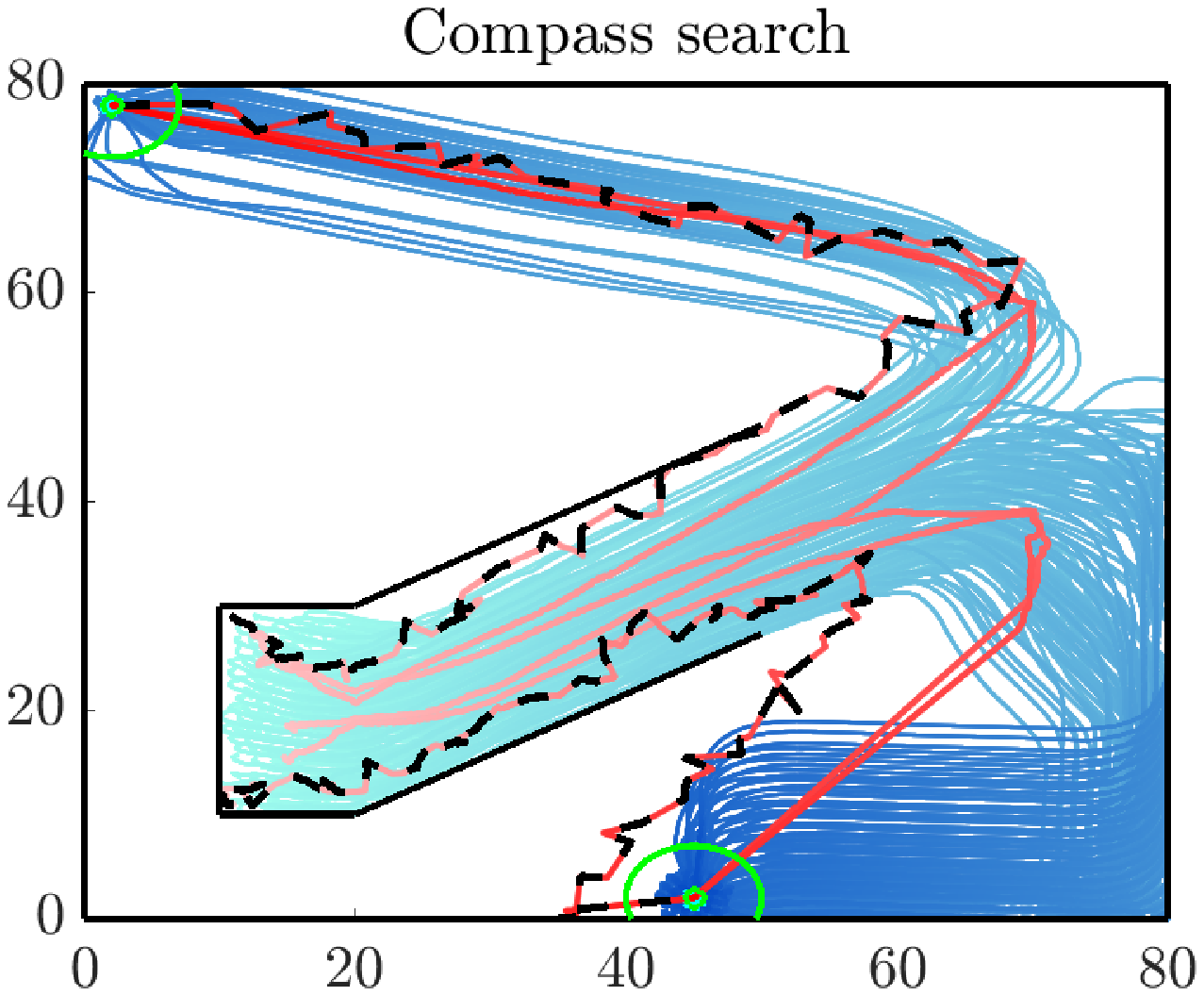}
	\caption{{\em Test 2.} Microscopic case: mass maximization in presence of obstacles. On the left, go-to-target. On the right, optimal compass search.  }
	\label{fig:test2D_mass_micro_1}
\end{figure}

\begin{figure}[h!]
	\centering
	\includegraphics[width=0.328\linewidth]{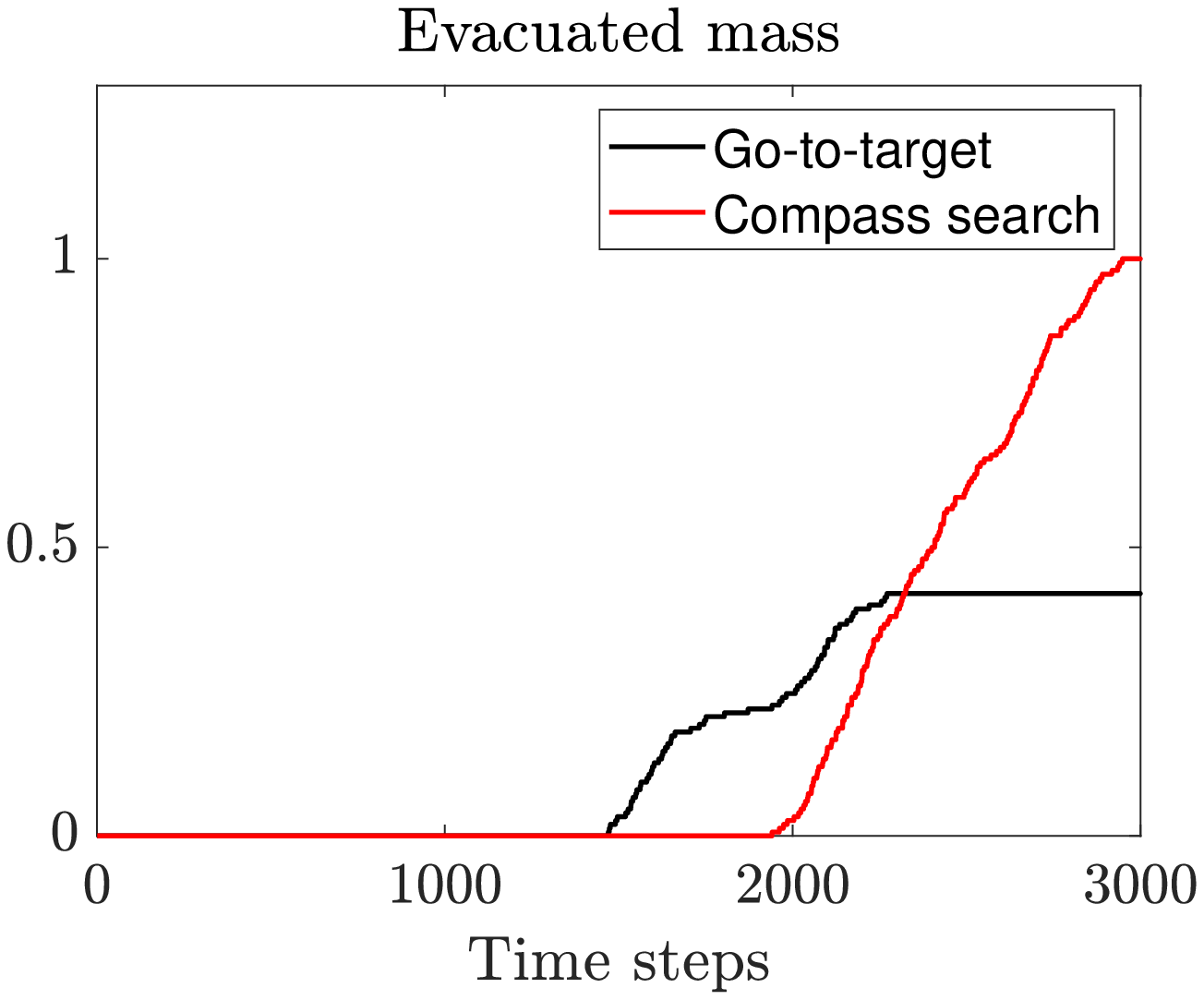}
	\includegraphics[width=0.328\linewidth]{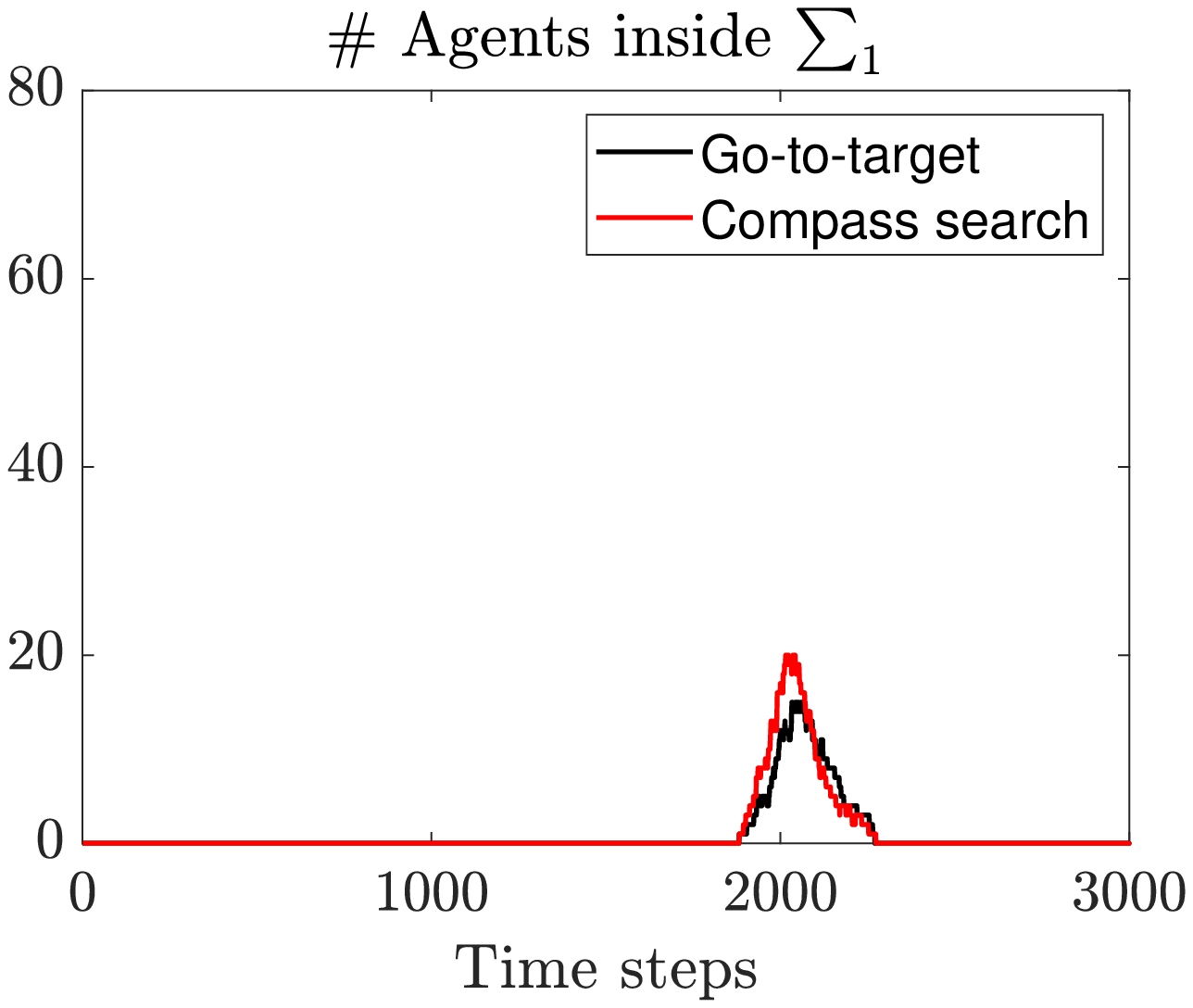}
	\includegraphics[width=0.328\linewidth]{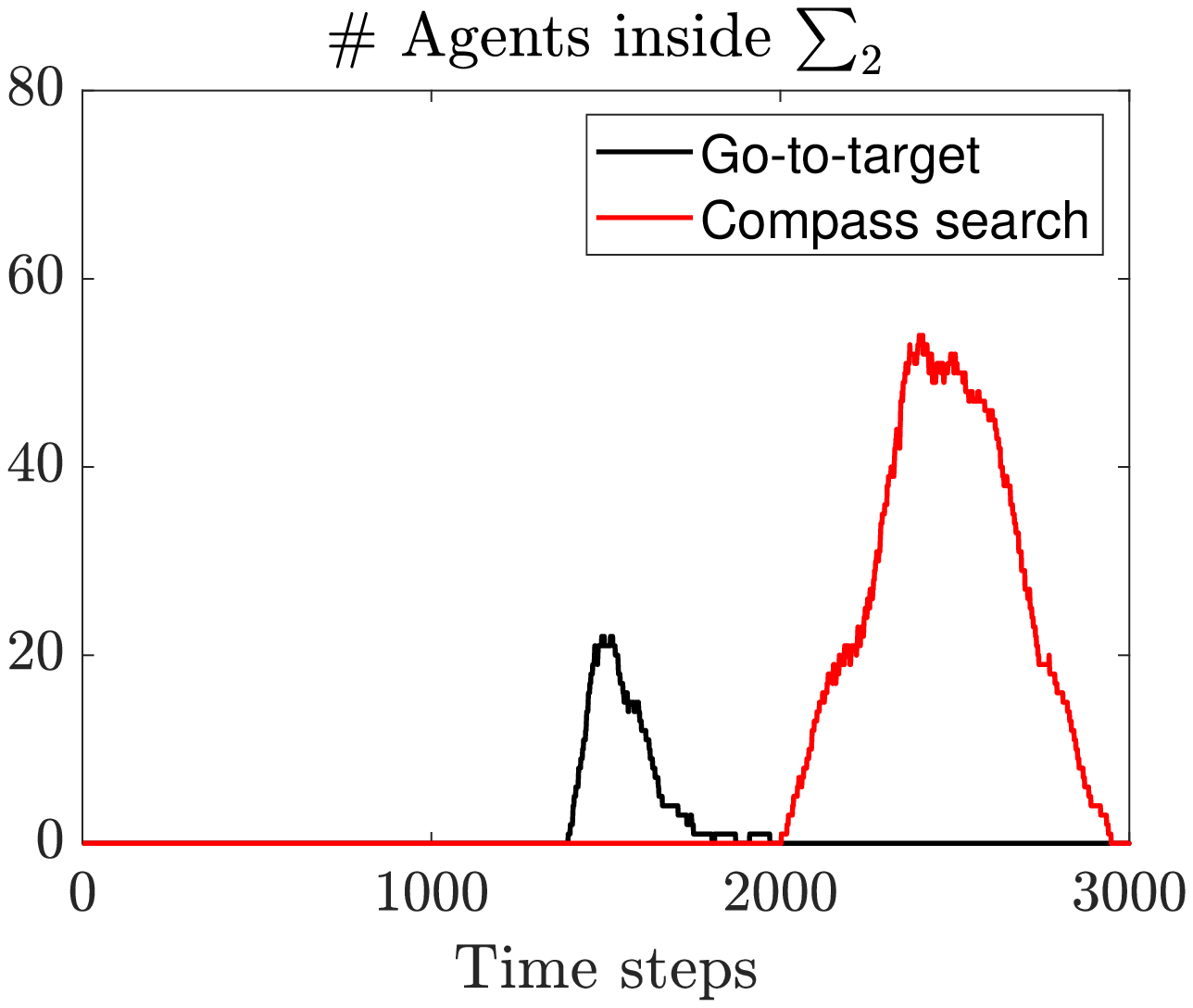}
	\caption{{\em Test 2.} Microscopic case: mass maximization in presence of obstacles. Evacuated mass (left), occupancy of the visibility area $\Sigma_1$ (centre) and $\Sigma_2$ (right) as a function of time for go-to-target and optimal compass search strategies.}
	\label{fig:test2D_mass_micro_2}
\end{figure}

\begin{figure}[h!]
	\centering
	\includegraphics[width=0.328\linewidth]{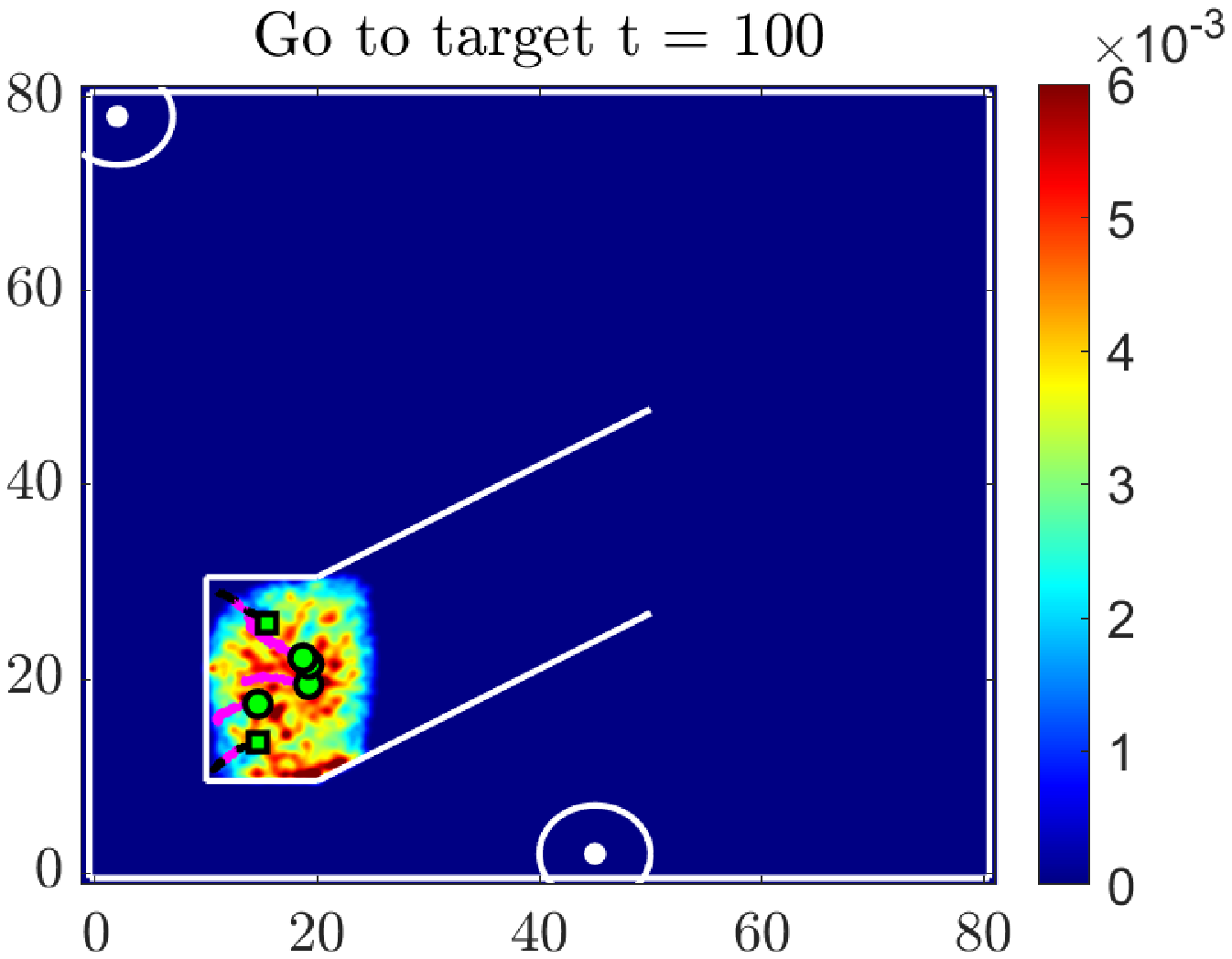}
	\includegraphics[width=0.328\linewidth]{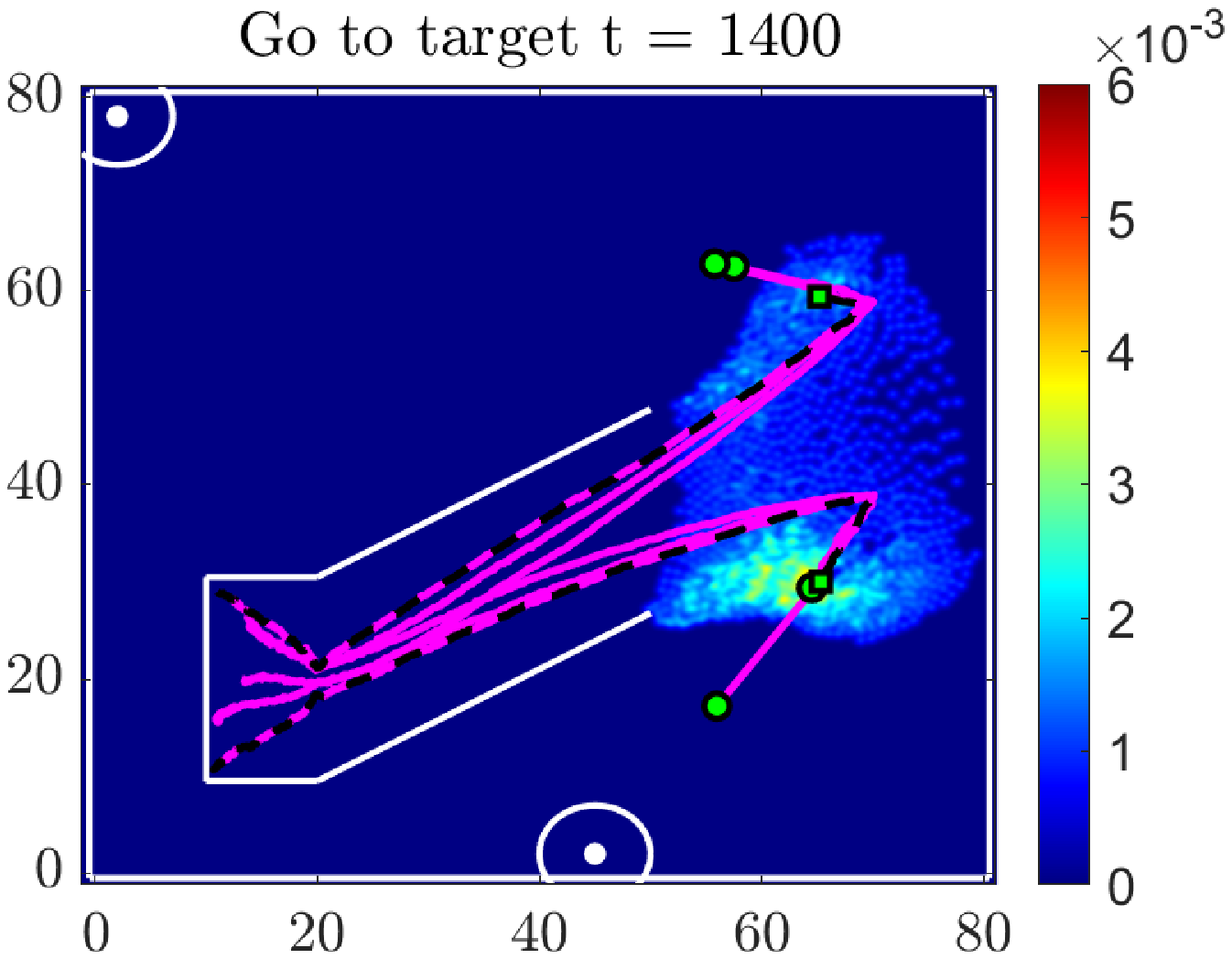}
	\includegraphics[width=0.328\linewidth]{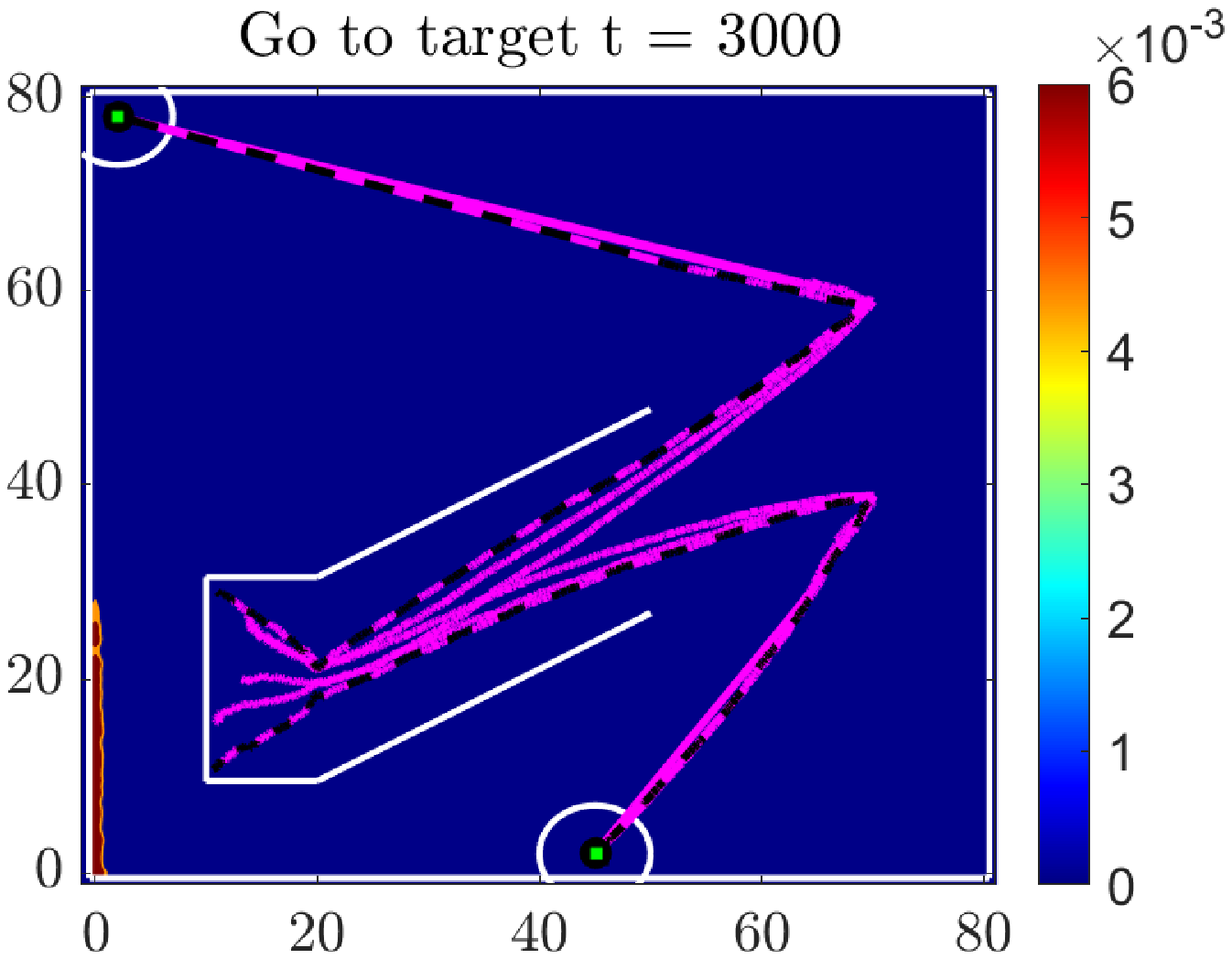}
	\\
	\includegraphics[width=0.328\linewidth]{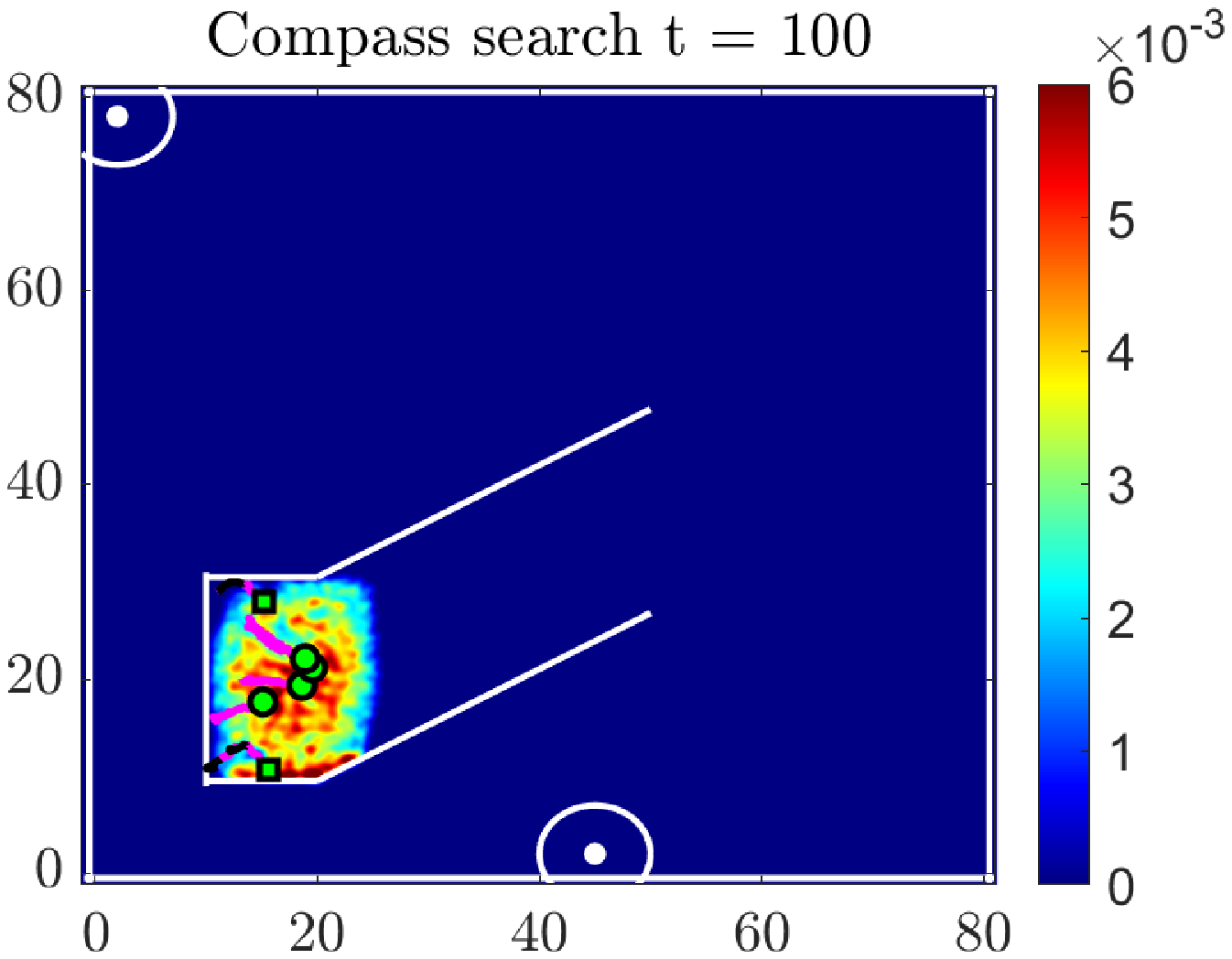}
	\includegraphics[width=0.328\linewidth]{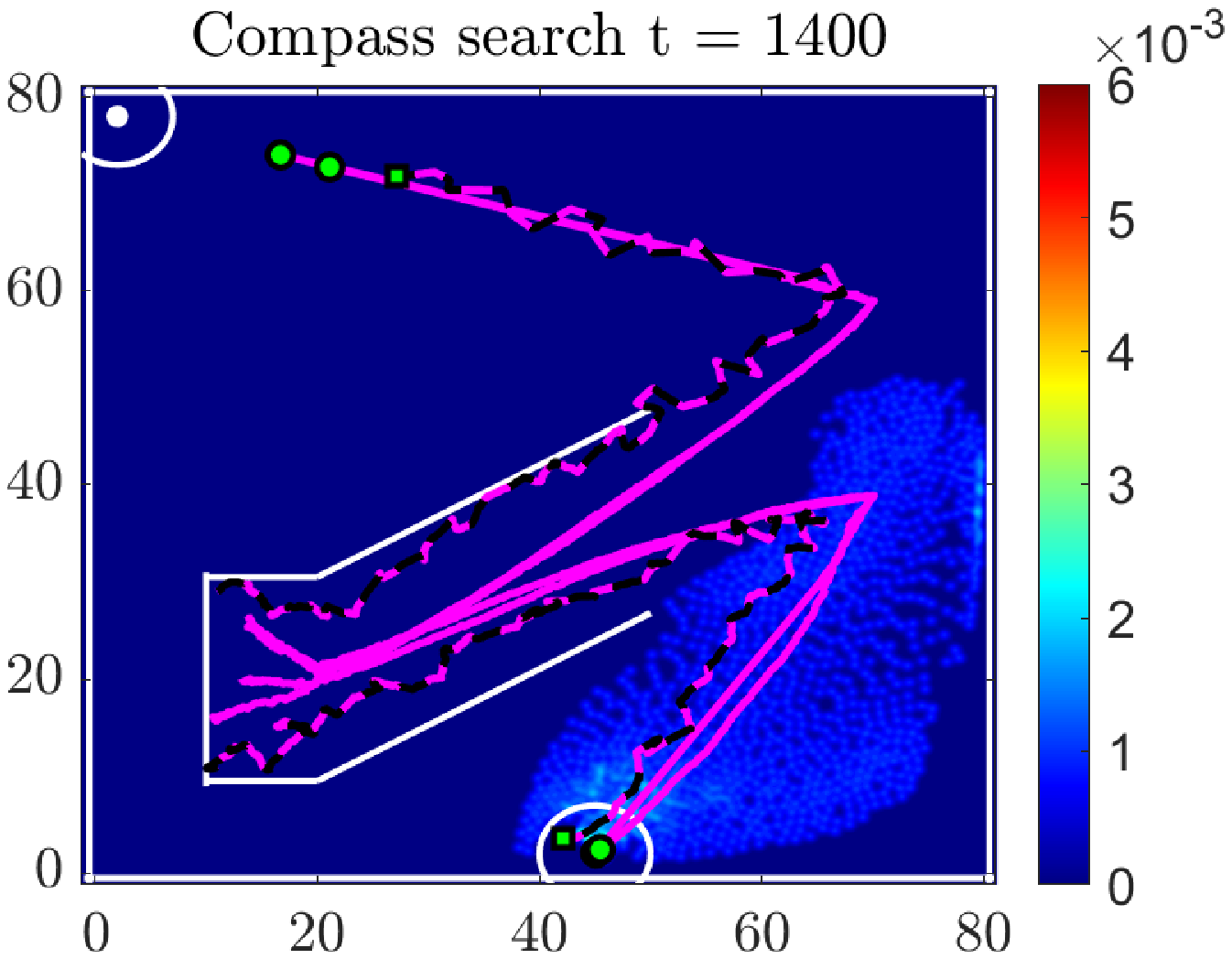}
	\includegraphics[width=0.328\linewidth]{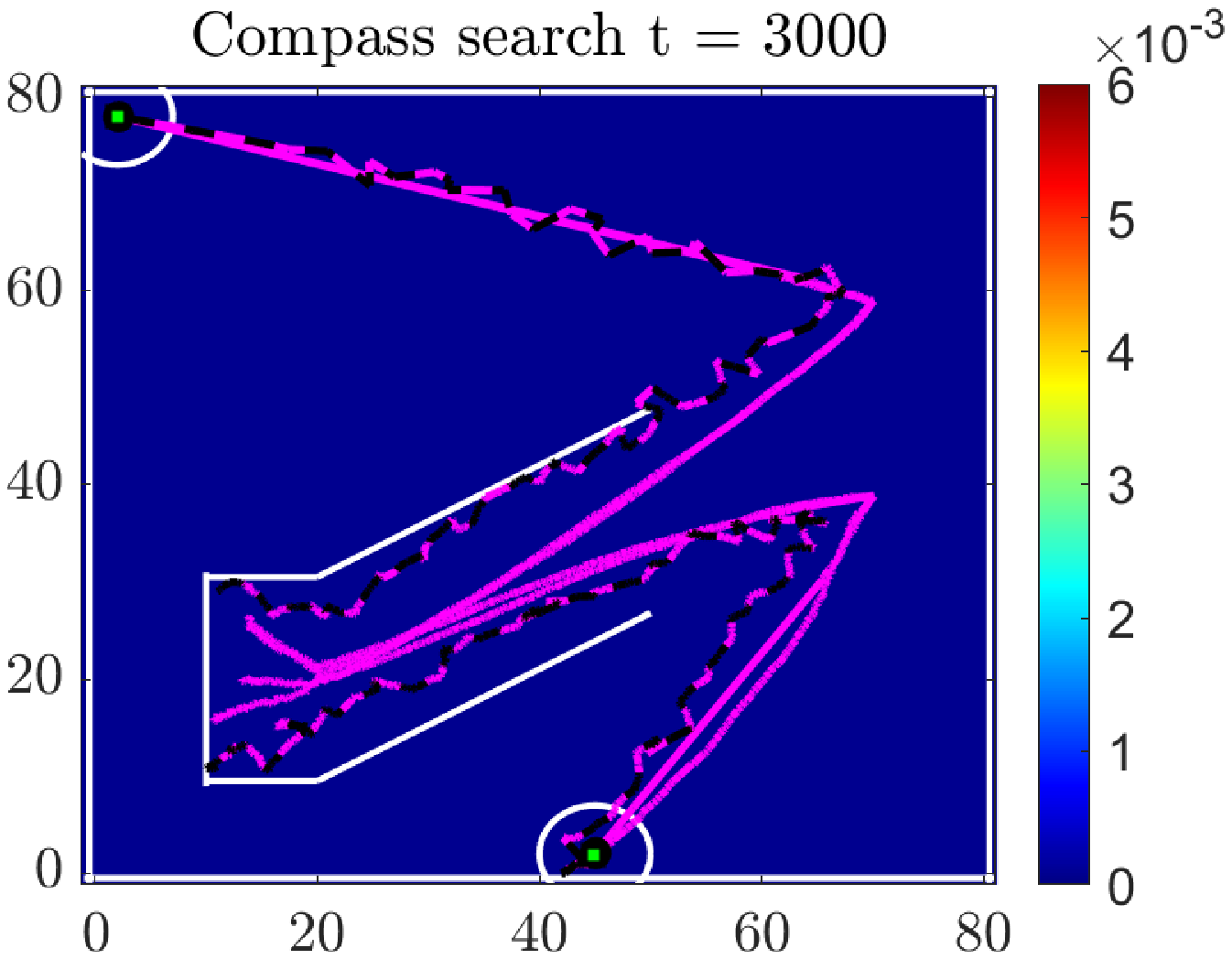}
	\caption{{\em Test 2.} Mesoscopic case: mass maximization in presence of obstacles. Upper row: three snapshots taken at time $t=100$, $t=1400$, $t=3000$ with the go-to-target strategy. Lower row: three snapshots taken at time $t=100$, $t=1400$, $t=3000$ with the optimized compass search strategy.}
	\label{fig:test2D_mass_meso_1}
\end{figure}

\begin{figure}[h!]
	\centering
	\includegraphics[width=0.328\linewidth]{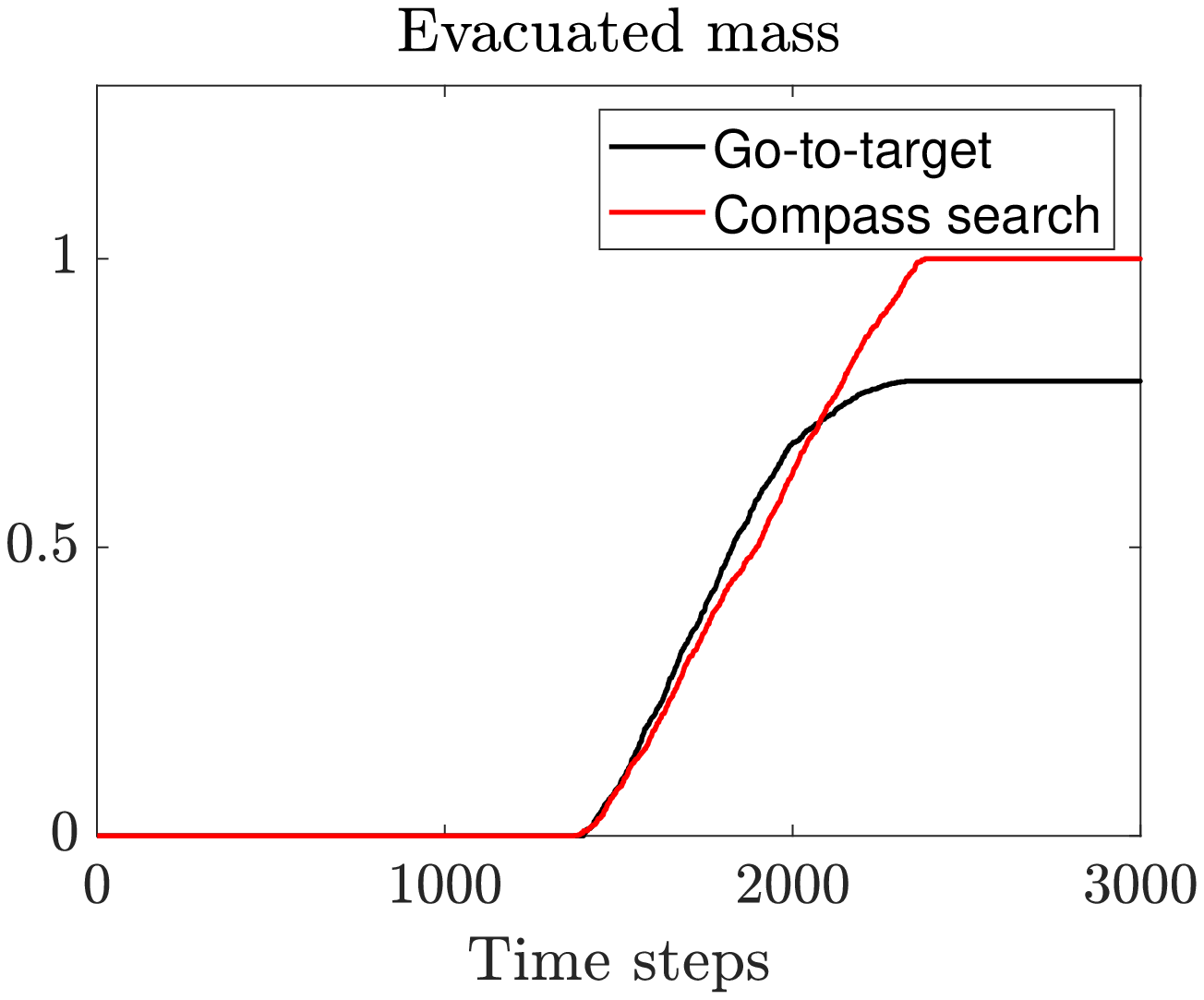}
	\includegraphics[width=0.328\linewidth]{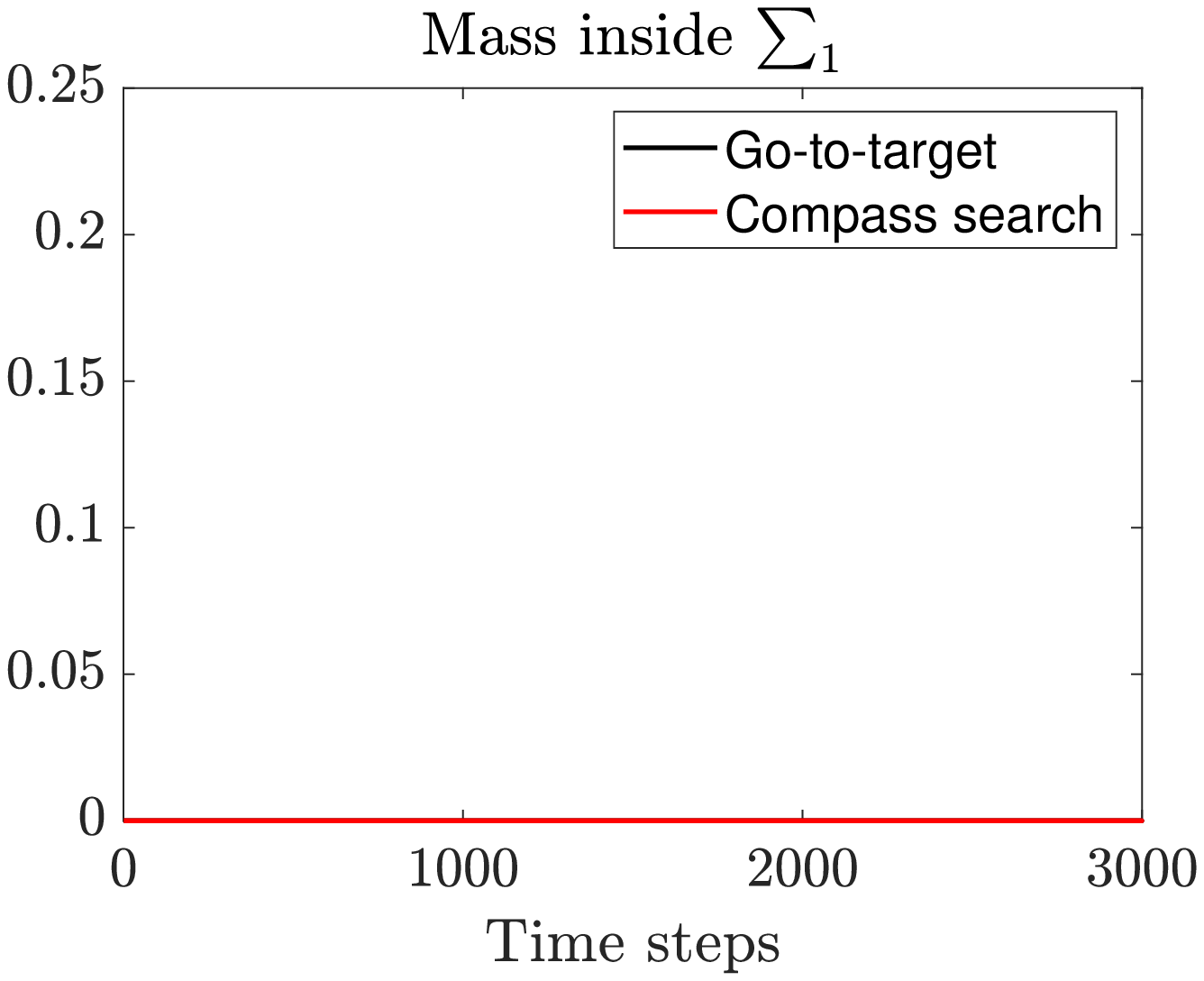}
	\includegraphics[width=0.328\linewidth]{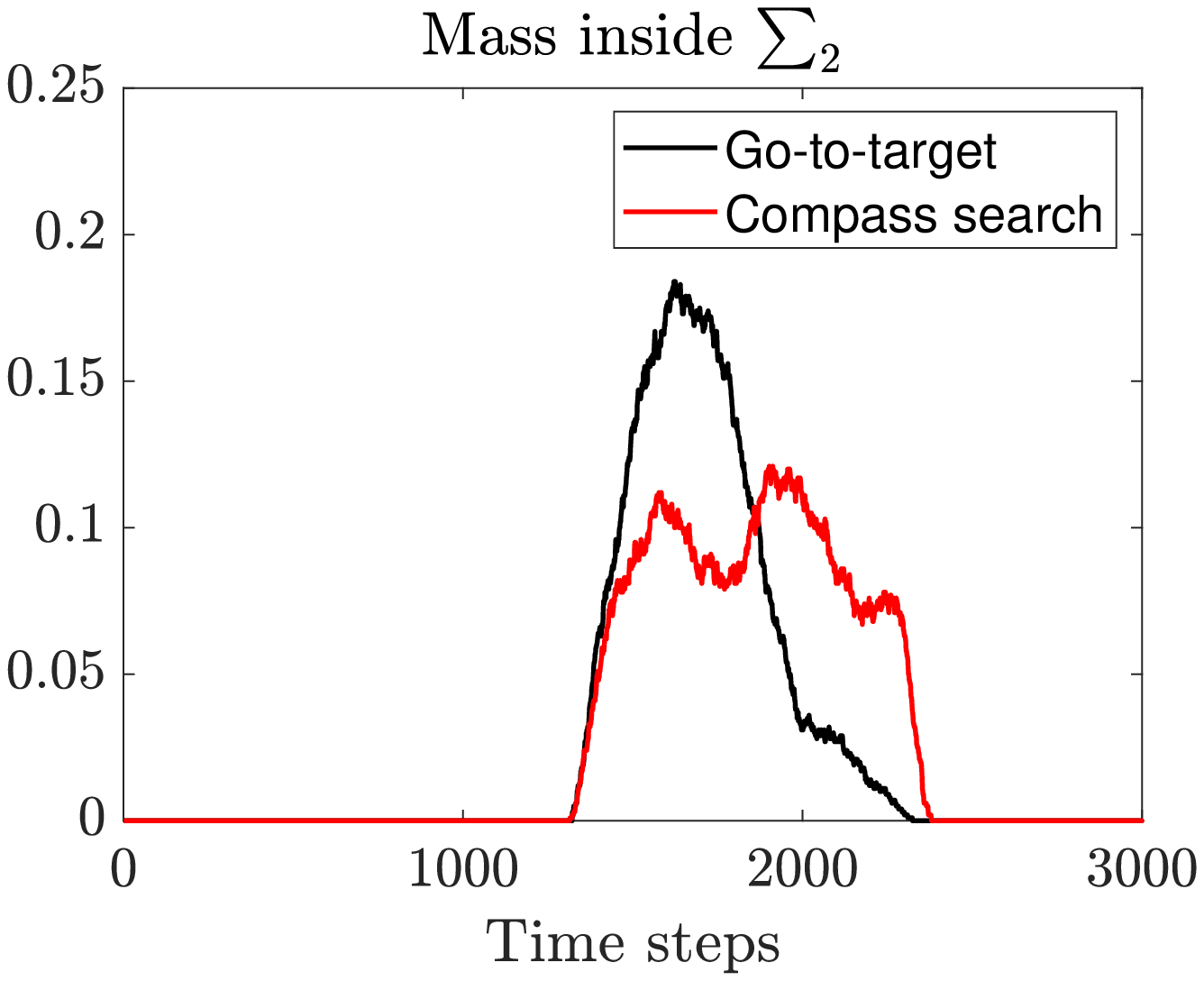}
	\caption{{\em Test 2.} Mesoscopic case: mass maximization in presence of obstacles. Evacuated mass (left), occupancy of the visibility area $\Sigma_1$ (centre) and $\Sigma_2$ (right) as a function of time for go-to-target and optimal compass search strategies.}
	\label{fig:test2D_mass_meso_3}
\end{figure}

{\em Microscopic case.} We consider $N^L=6$ leaders, with two aware leaders. Initially, followers have zero velocity. Three leaders, only one aware, will move towards exit $x_1^\tau$, and the remaining towards exit $x_2^\tau$. We report in
Figure \ref{fig:test2D_mass_micro_1} the evolution with the go-to-target strategy on the left, and with optimized strategies for the two aware leaders on the right. With go-to-target strategy leaders first leave the room and then move towards the exits. Since leaders move rapidly towards the exits, their influence over followers vanishes after a certain time. Indeed, part of the followers hits the right boundary wall and does not reach the exits. Instead, with optimized strategies, leaders are slowed down, as consequence followers are influenced by leaders for a larger amount of time. Table \ref{eq:test2_table1} reports the comparison between two strategies in terms of evacuated mass, where with only three iterations of the optimization method total evacuation is accomplished.\begin{table}
	\centering
	\caption{{\em Test 2.} Performance of various strategies for obstacle case with two exits in the microscopic case.}\label{eq:test2_table1} 
	\begin{tabular}{cccc}
	& go-to-target  & CS (3 it) \\
  \cmidrule(r){2-2}\cmidrule(r){3-3}
  	Evacuation time (time steps) &>3000 & 2948 \\
	Evacuated mass (percentage)  & 42\%  & 100\% \\
	\hline
\end{tabular}
\end{table}
In Figure \ref{fig:test2D_mass_micro_2} we compare the cumulative distribution of evacuated mass and the occupancy of the exits visibility areas as a function of time for go-to-target strategy and optimized strategy. We remark that with minimal change of the fixed strategy we reach evacuation of the total mass.
{\em Mesoscopic case.} Consider now the case of continuous mass of followers, and the equivalent setting as in the microscopic case. Initial configuration is reported in Figure \ref{fig:test2D_mass_micro_0}. We report the evolution of the two scenarios in Figure \ref{fig:test2D_mass_meso_1}, where in the upper row we depict three different time frames of the dynamics obtained with go-to-target strategy. Once leaders have moved outside the inner room, at time $t=1400$, followers mass splits into two parts. However, only leaders moving towards the lower exit $x_2^\tau$ are able of steering the followers towards the target, the rest of the followers moving upwards get lost and at final time $t=3000$ is located close to the left wall. Hence, partial evacuation of followers is achieved, as shown in Table \ref{eq:test2_table3} we retrieve $78.8\%$ of total mass evacuated. Only one exit is used, this may cause problems of heavy congestion around the exits.
Bottom row of Figure \ref{fig:test2D_mass_meso_1} shows the situation with optimized leaders strategy. Differently from the previous case at time $t=2380$ the whole mass has been evacuated, part of the followers mass reaches the lower exits and the remaining mass reaches $x_1^\tau$ after a while. In Table \ref{eq:test2_table3} we reported the performances of the two approaches.
\begin{table}
	\centering
	\caption{{\em Test 2.} Performances of total mass evacuation problems in the mesoscopic case.}\label{eq:test2_table3} 
	\begin{tabular}{cccc}
		& go-to-target  & CS  (5 it)\\
		 		\cmidrule(r){2-2}\cmidrule(r){3-3}
		 		Evacuation time (time steps) &>3000 & 2380\\
		Evacuated mass (percentage)  & 78,8 \%  & 100 \% \\
		\hline
	\end{tabular}
\end{table}
In Figure \ref{fig:test2D_mass_meso_3} we compare the evacuated mass and the occupancy of the exits visibility zone as a function of time for go-to-target strategy and opitmized strategy after 5 iterations of compass search method.  
\subsubsection{Test 3: Optimal mass splitting over multiple exits}\label{sec:test3}
Problems of heavy congestion and overcrowding around the exits arise naturally in evacuation and, in real-life situations, they can cause injuries due to over-compression and suffocation. Instead of maximizing the total evacuated mass or the minimum time, we ask to distribute the total evacuated mass at final time $T$ between all the exits as reported in \eqref{eq:test3intro}. The choice of mass redistribution among the different exits can be done according to the specific application and environment. In what follows we consider two different examples, both with two exits, and we will require that mass splits uniformly between the two targets.
\paragraph{\bf Setting 1) Two exits in a close environment.} As first example we consider the same setting of Test 2, where complete evacuation was achieved, but all followers were directed toward a single exit. In this case we aim to optimize leaders strategies in order to equidistribute the total mass of follower among the two exits.

\begin{figure}[h!]
	\centering
	\includegraphics[width=0.45\linewidth]{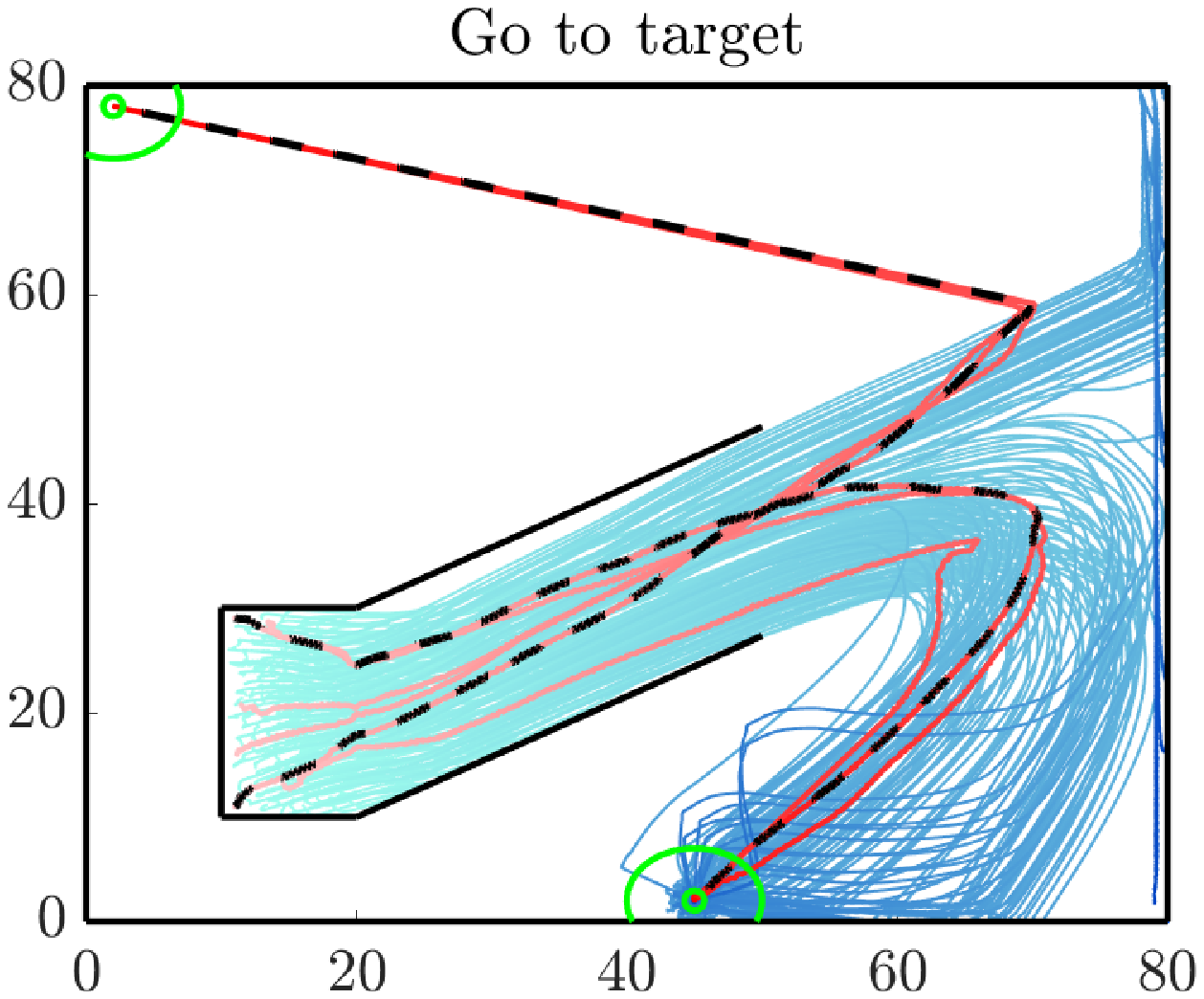}
	\includegraphics[width=0.45\linewidth]{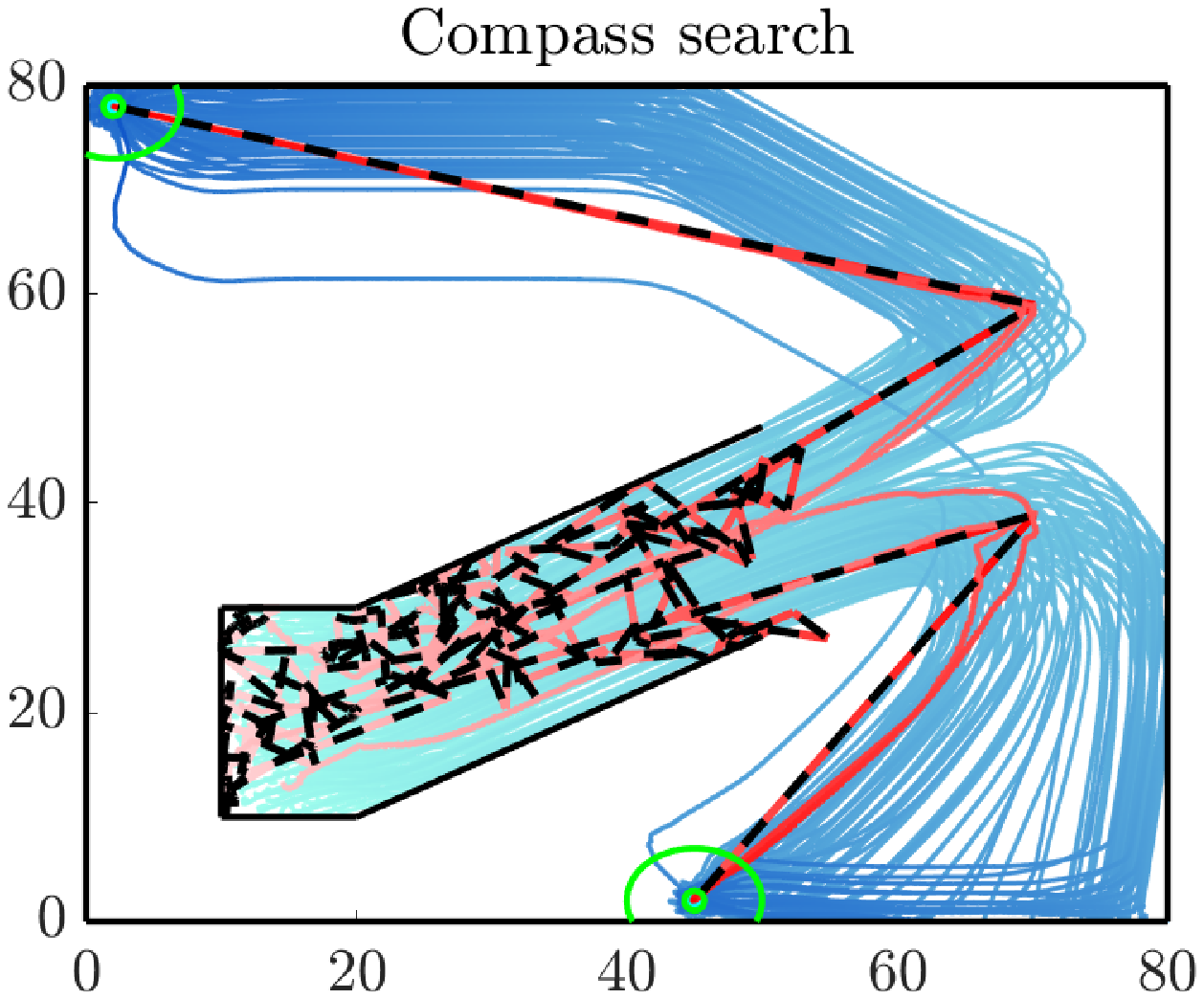}
	\caption{{\em Test 3a.} Microscopic case: mass splitting in presence of obstacles. On the left, go-to-target. On the right, optimal compass search. }
	\label{fig:test3_splitting_micro_1} 
\end{figure}

\begin{figure}[h!]
	\centering
	\includegraphics[width=0.328\linewidth]{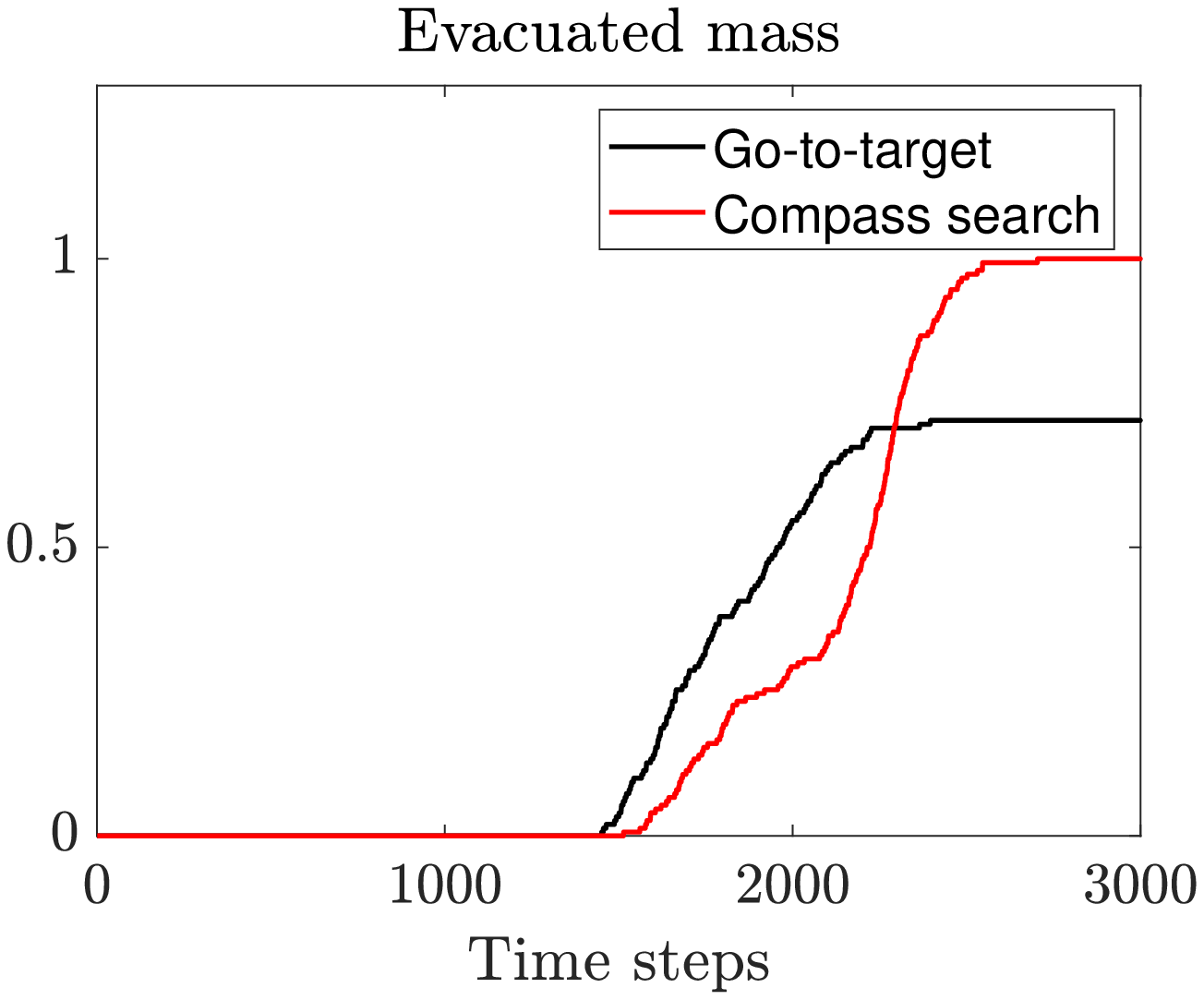}
	\includegraphics[width=0.328\linewidth]{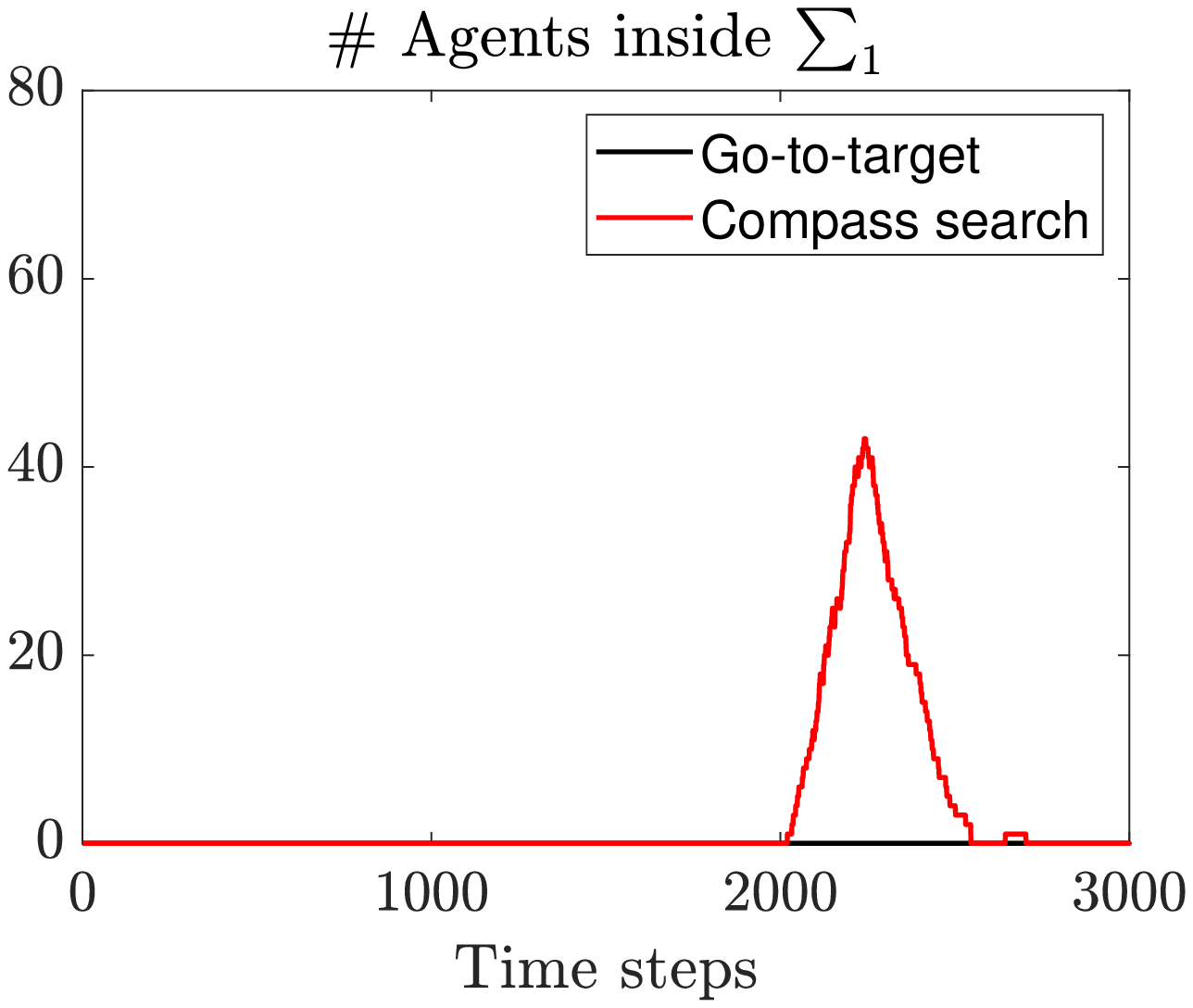}
	\includegraphics[width=0.328\linewidth]{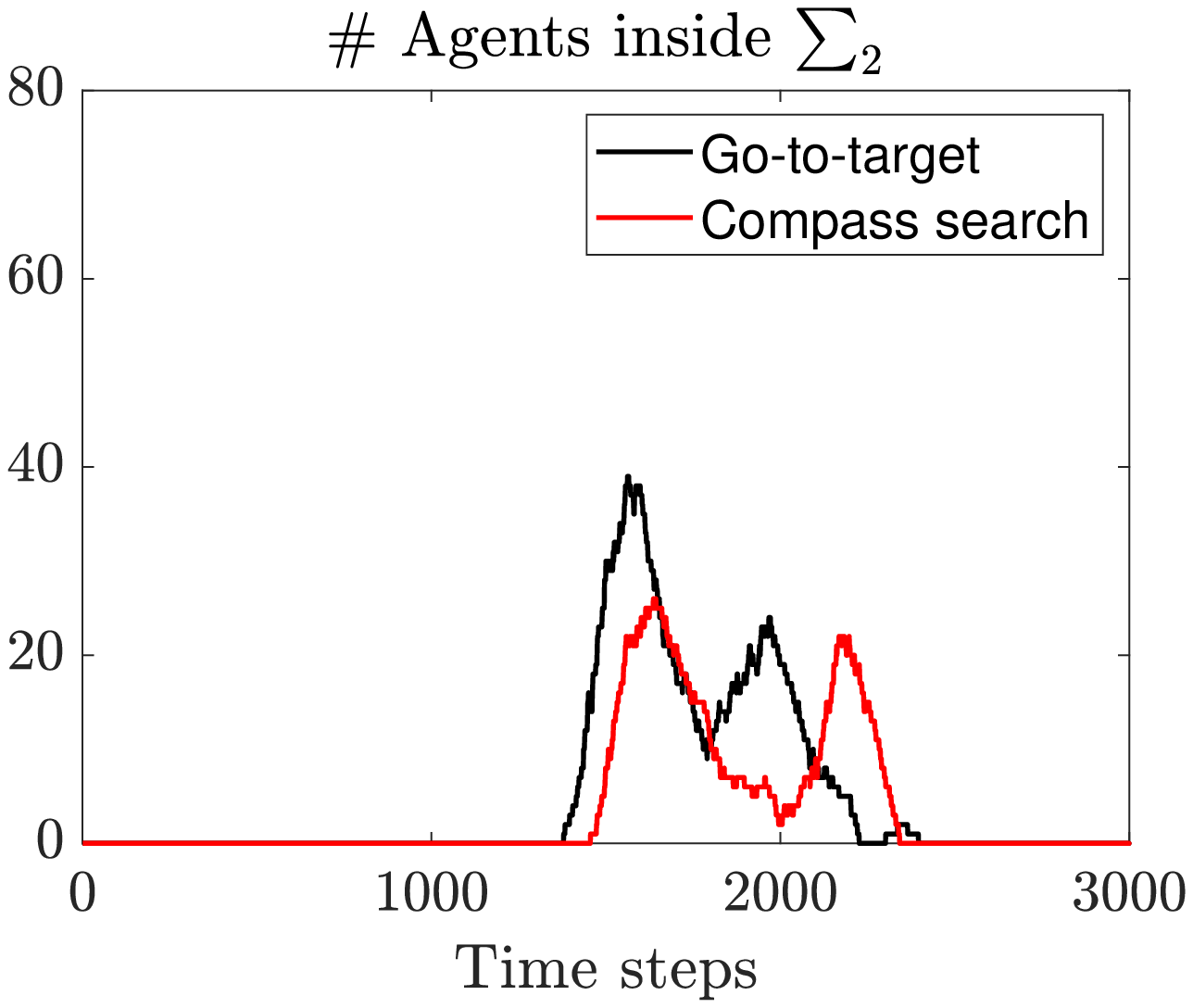}
	\caption{{\em Test 3a.} Microscopic case: mass splitting in presence of obstacles. Evacuated mass (left), occupancy of the visibility area $\Sigma_1$ (centre) and $\Sigma_2$ (right) as a function of time for go-to-target and optimal compass search strategies.}
	\label{fig:test3_splitting_micro_2}
\end{figure}

{\em Microscopic case.} In Figure \ref{fig:test3_splitting_micro_1} we depict the scenario for the fixed strategy  and the optimized one. We observe that again with go-to-target strategy the complete evacuation is not achieved. Moreover, since the vast majority of followers reach the lower exit $x_2^\tau$, heavy congestion is formed in the visibility area $\Sigma_2$. On the other hand with an optimized strategy two aware leaders slow down their motion spending more time inside the inner room.  In this way, followers are split between the two exits, and the entire mass is evacuated at final time. In Table \ref{eq:test3_table1} we report the performances of the two strategies, where for optimized strategy we have $45\%$ of mass in $x_1^\tau$ and $55\%$ in $x_2^\tau$. 
\begin{table}\caption{{\em Test 3a.} Performances of mass splitting  in the microscopic case.}\label{eq:test3_table1} 
	\centering
	\begin{tabular}{ c c c } 
		& go-to-target & CS (50 it)\\
		\cmidrule(r){2-2}\cmidrule(r){3-3}
		Evacuation time (time steps)  & > 3000 & 2704 \\
		Mass evacuated from $E_1$ & $0\%$ & $45\%$\\
		Mass evacuated  from $x_2^\tau$ & $72\%$ & $55\%$\\
	    Total mass evacuated & 72\% & 100\% \\
		\hline
	\end{tabular}
\end{table}
In Figure \ref{fig:test3_splitting_micro_2} we report the evacuated mass and the occupancy of the exits visibility zone as a function
of time for go-to-target strategy and optimal compass search strategy.  Note that, with the compass search strategy, the whole mass is split between the two exits reducing the overcrowding in the visibility region.

\begin{figure}[h]
	\centering
	\includegraphics[width=0.328\linewidth]{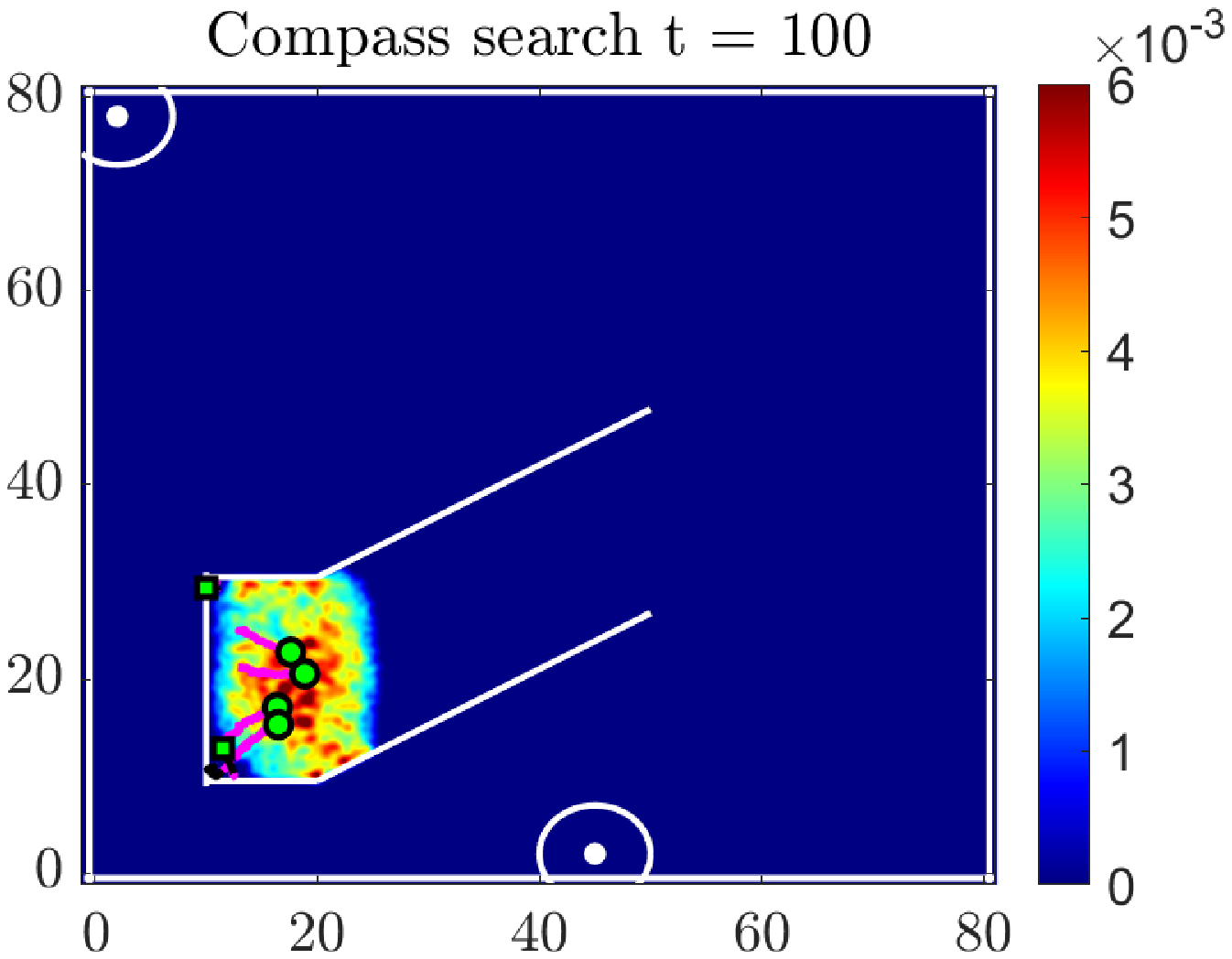}
	\includegraphics[width=0.328\linewidth]{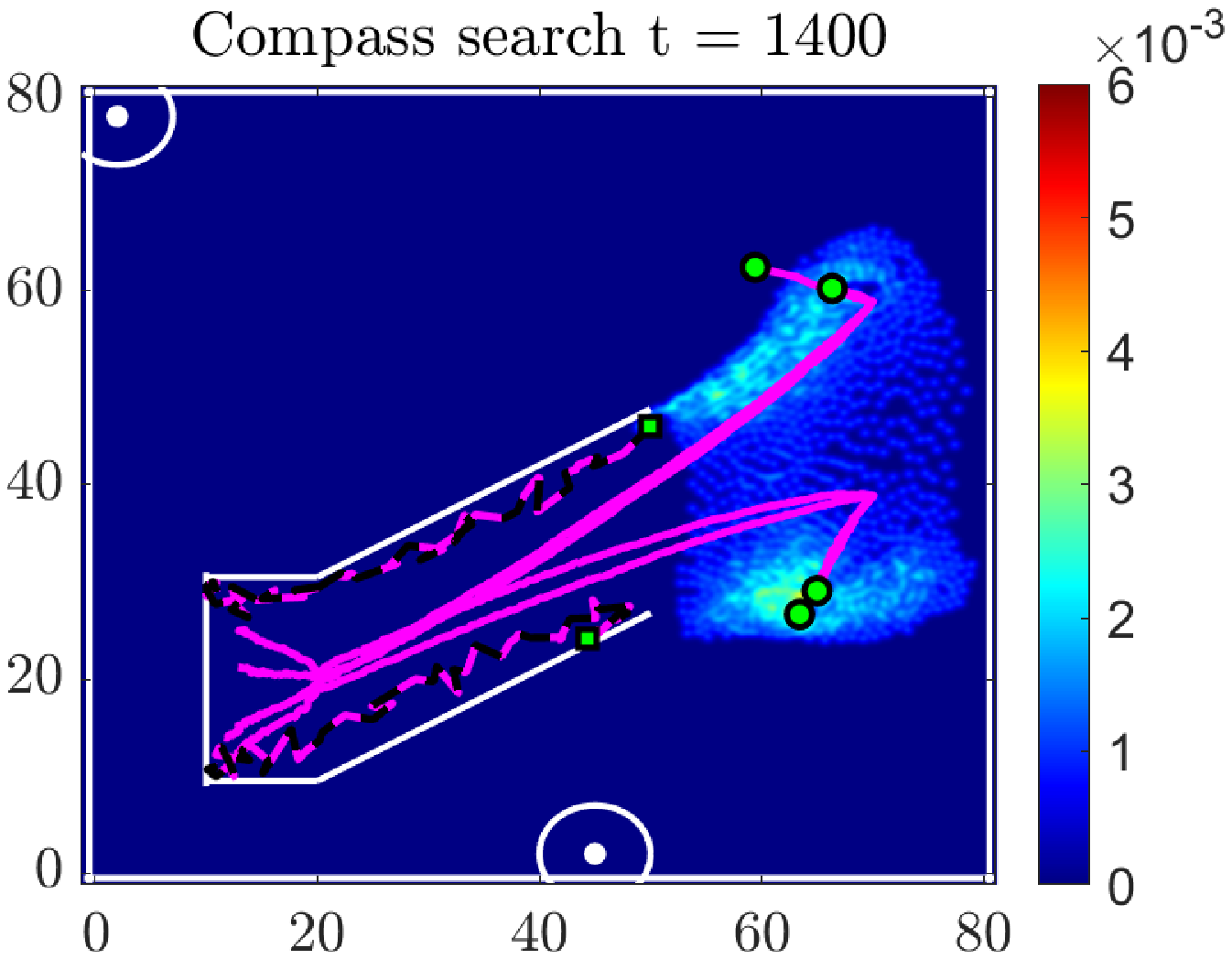}
	\includegraphics[width=0.328\linewidth]{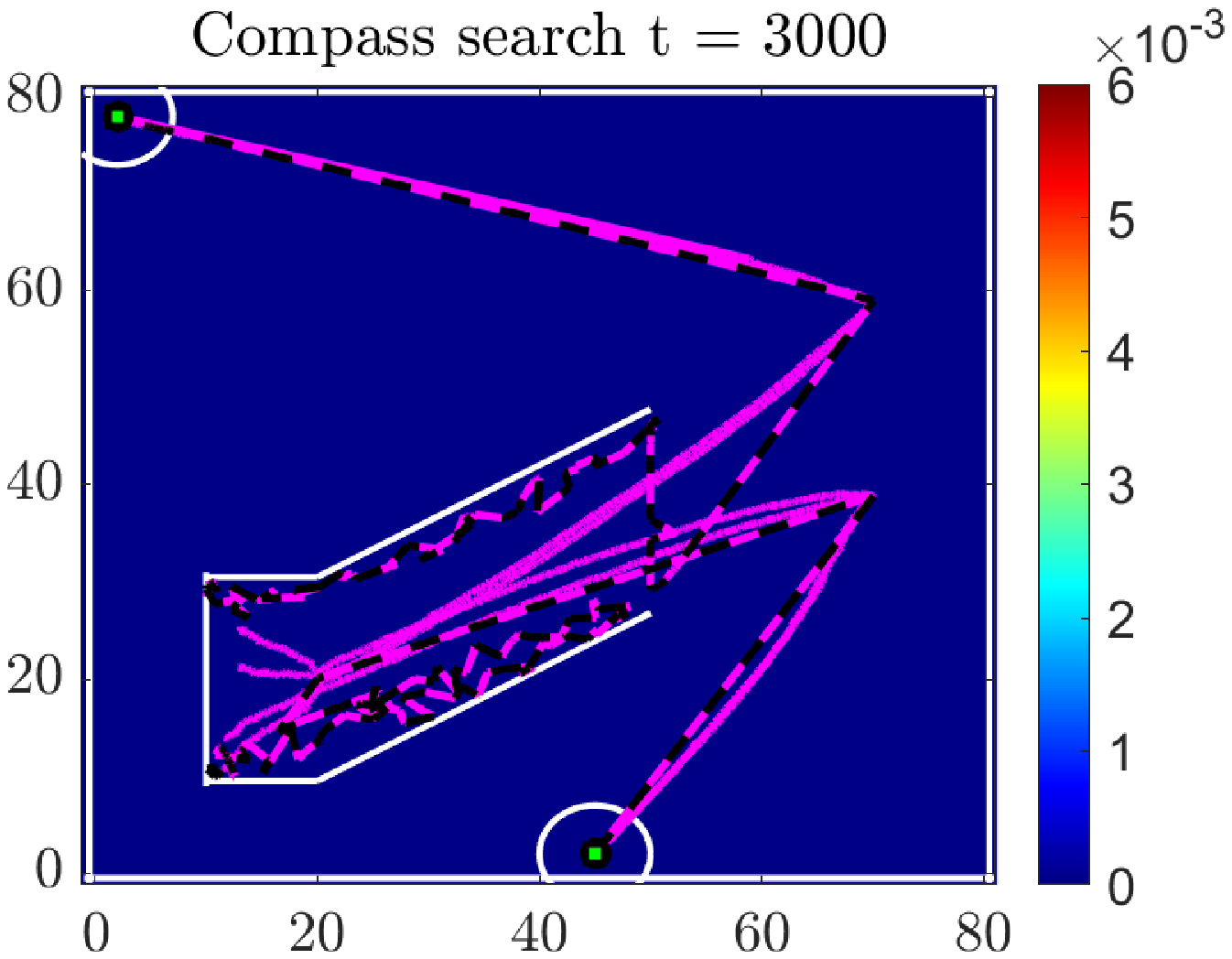}
	\caption{{\em Test 3a.} Mesoscopic case: mass splitting in presence of obstacles. Three snapshots taken at time $t=100$, $t=1400$, $t=3000$ with the optimal compass search strategy. For the go-to-target case we refer to the first row of Figure \ref{fig:test2D_mass_meso_1}.}
	\label{fig:test3_splitting_meso_1}
\end{figure}

{\em Mesoscopic case.}  We report now the case of a continuum density of followers. For the go-to-target strategy, we consider the same dynamics of the previous test, in this case the mass of followers does not split between the two exits, as shown in Figure \ref{fig:test2D_mass_meso_1}, and the $78,8\%$ reaches exit $x_2^\tau$.
In Figure \ref{fig:test3_splitting_meso_1}, three snapshots were taken at three different times with the compass search strategy.
At time $t=100$, leaders move to evacuate the followers mass out of the inner room. At time $t=1400$, the followers mass splits in two masses, one moving towards the upper and the other towards the lower exits. At time $t=3000$, almost all the followers mass is evacuated. The mass is split between the two exits as shown in Table \ref{eq:test3_table2}.
\begin{table}\caption{{\em Test 3a}. Performances of mass splitting in the mesoscopic case.}\label{eq:test3_table2} 	\centering
	\begin{tabular}{ c c c } 
	& go-to-target & CS (50 it)\\
	\cmidrule(r){2-2}\cmidrule(r){3-3}
	Evacuation time (time steps)  & > 3000 & > 3000 \\
	Mass evacuated from $x_1^\tau$ & $0\%$ & $49\%$\\
	Mass evacuated  from $x_2^\tau$ & $78,8\%$ & $50\%$\\
	Total mass evacuated & 78,8\% & 99\% \\
	\hline
\end{tabular}
\end{table}
In Figure \ref{fig:test2D_mass_meso_3} we compare the evacuated mass and the occupancy of the exits visibility zone as a function
of time for go-to-target strategy and optimal compass search strategy. With the compass search technique the occupation of the visibility areas is reduced since the splitting of the total mass between the two exits is optimized. Hence, the risk of injuries due to overcrowding in real-life situations should be reduced.
\begin{figure}[h!]
	\centering
	\includegraphics[width=0.328\linewidth]{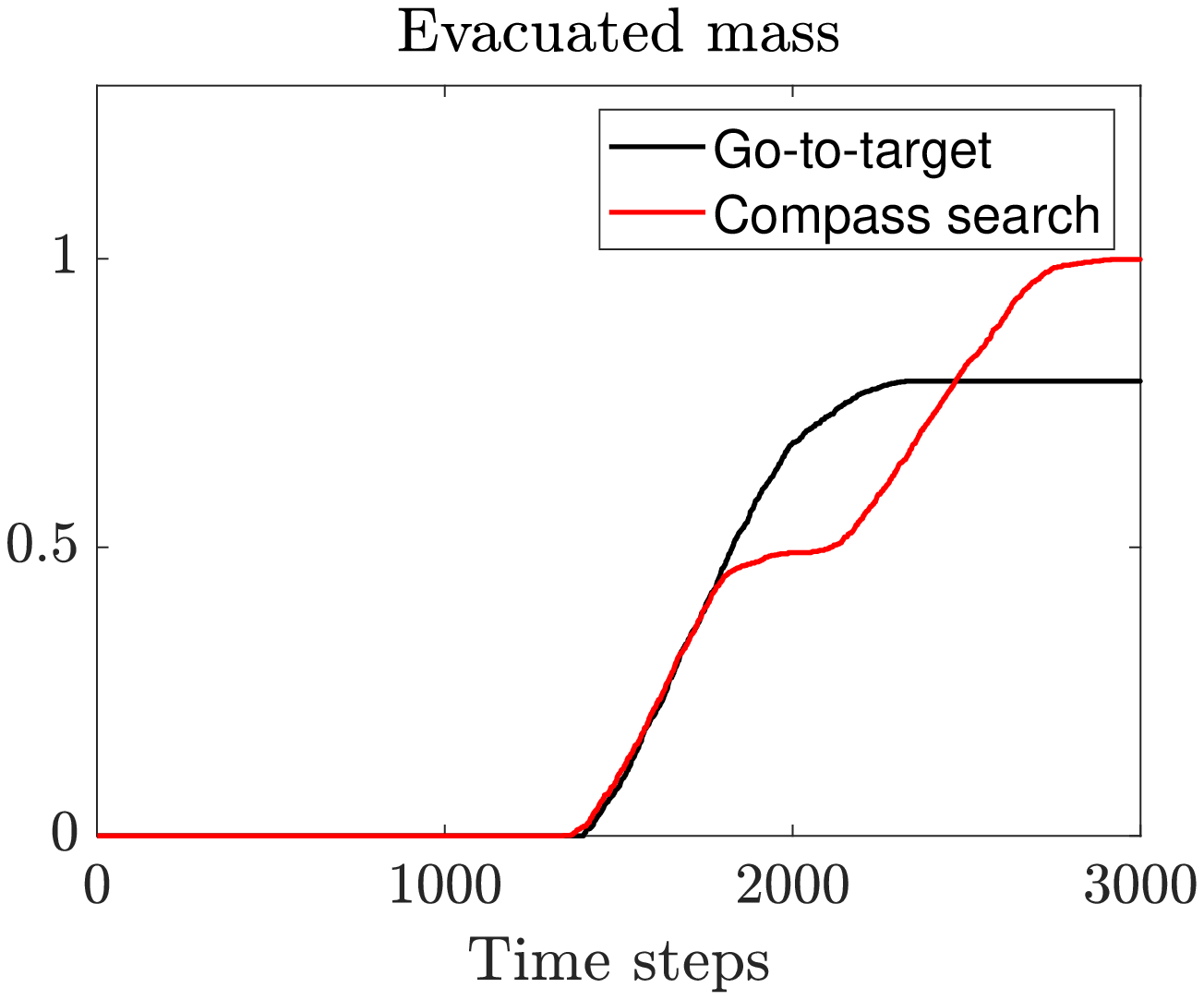}
	\includegraphics[width=0.328\linewidth]{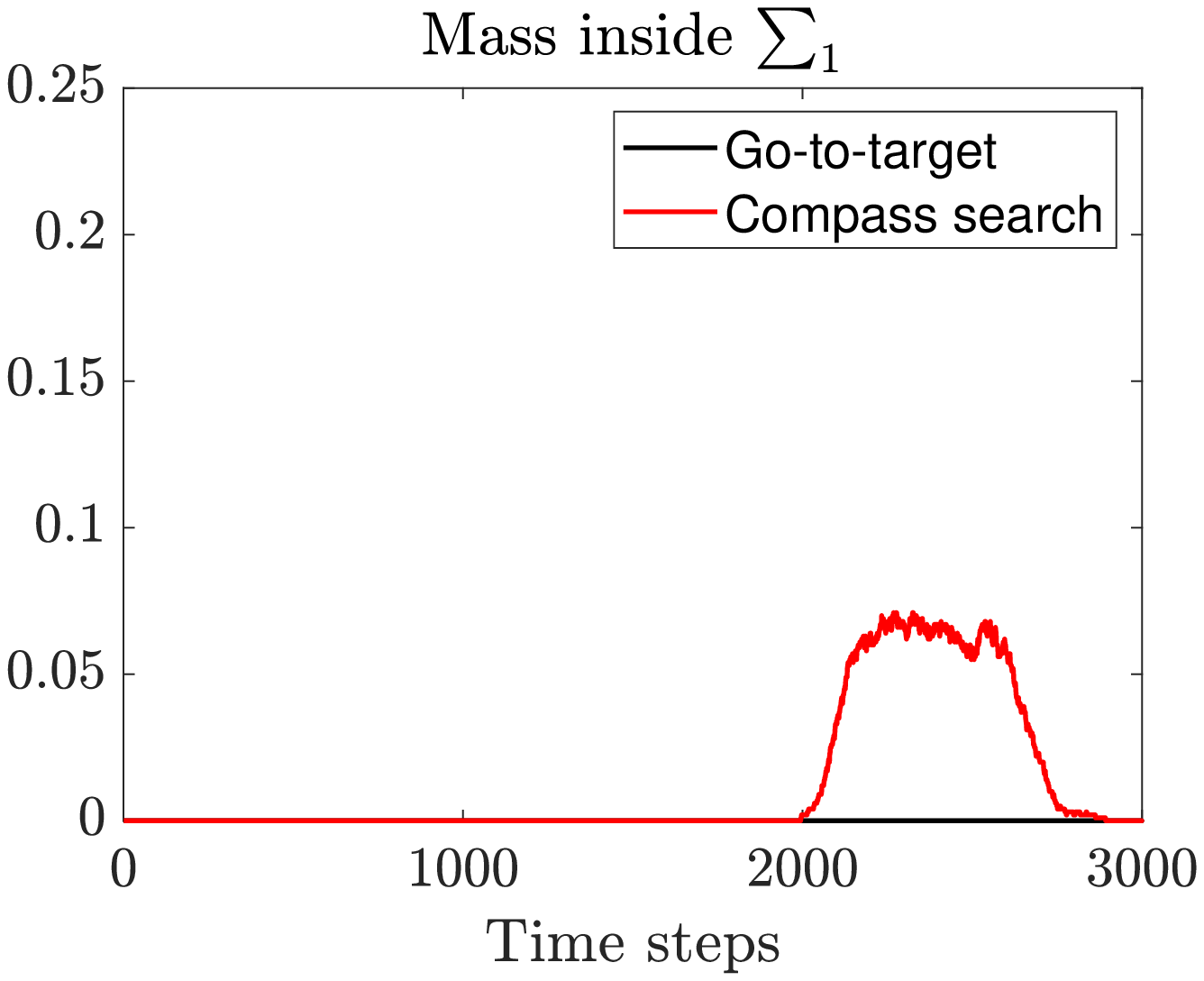}
	\includegraphics[width=0.328\linewidth]{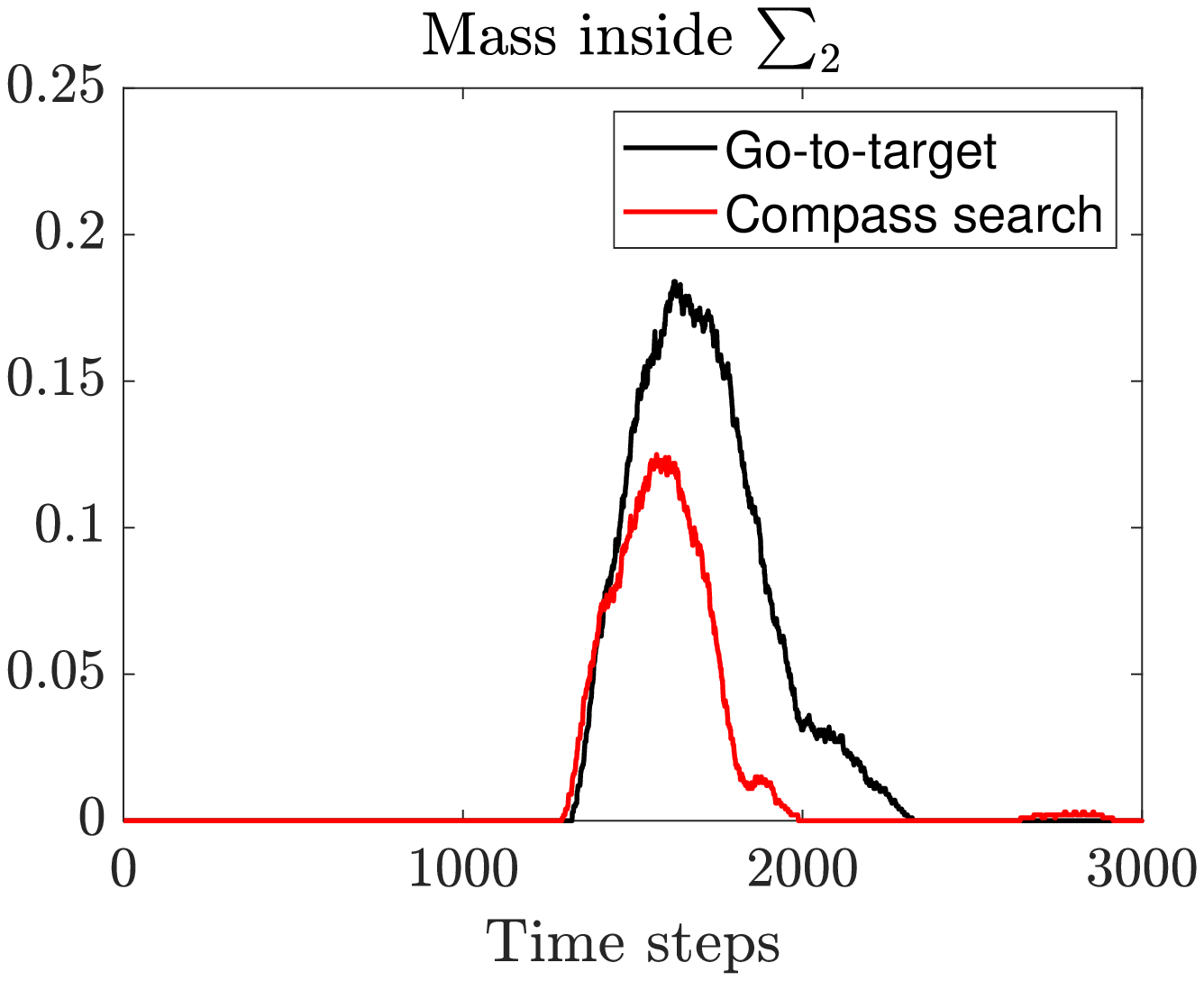}
	\caption{{\em Test 3a.} Mesoscopic case: mass splitting in presence of obstacles. Evacuated mass (left), occupancy of the visibility area $\Sigma_1$ (centre) and $\Sigma_2$ (right) as a function of time for go-to-target and optimal compass search strategies.}
	\label{fig:test3D_splitting_meso_3}
\end{figure}

\begin{figure}[h!]
	\centering
	\includegraphics[width=0.495\linewidth]{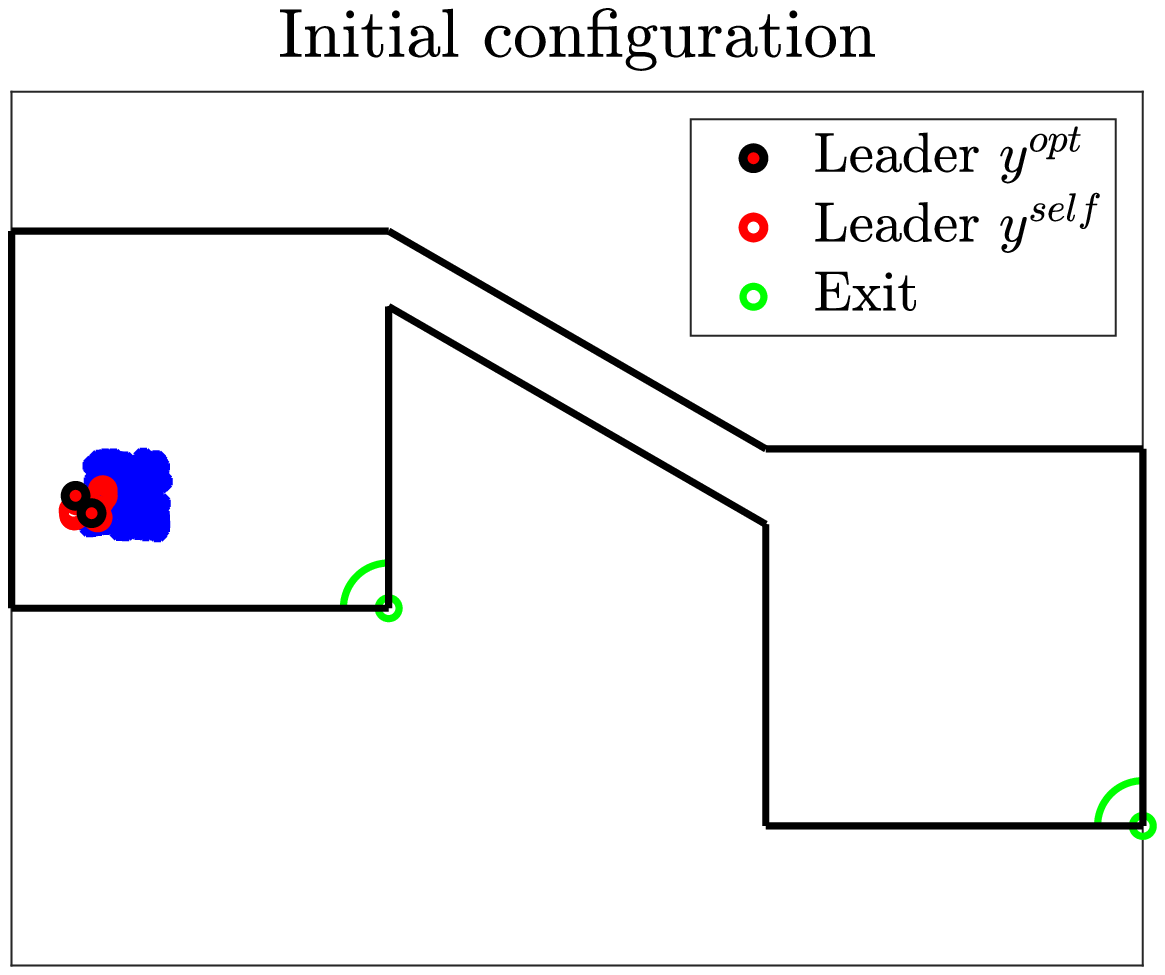}
	\includegraphics[width=0.495\linewidth]{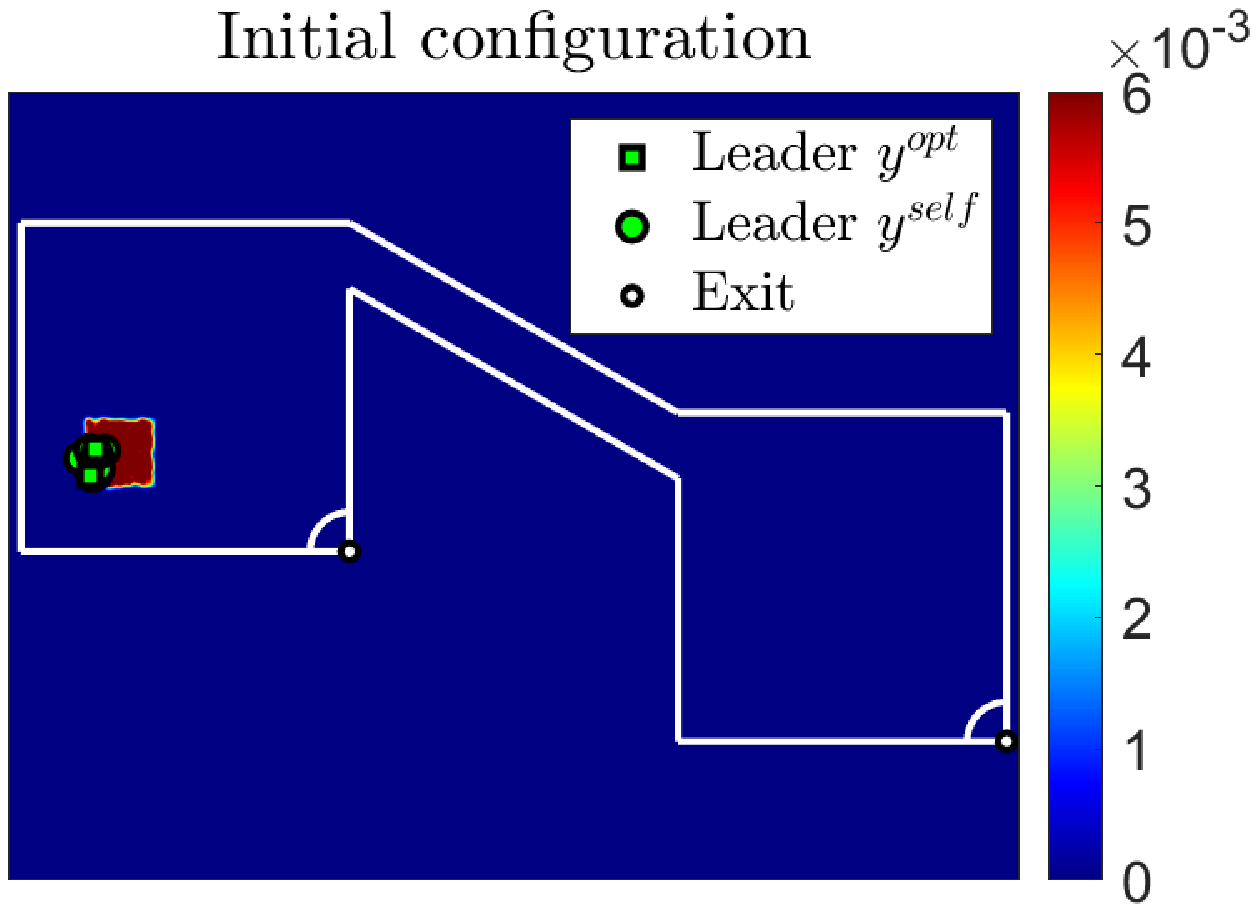}
	\caption{{\em Test 3b.} Mass splitting in presence of staircases, initial configuration.}
	\label{fig:test3_stair_0}
\end{figure}

\paragraph{\bf Setting b) Two exits with staircases.}
Consider two rooms and two exits, limited by walls, positioned at different floors, and connected by a staircase. Each room has an exit located in the bottom right corner. We assume followers and leaders to be uniformly distributed in a square inside the first room. Similar to the previous case, we assume that the model includes eight unaware and two aware leaders in total $N^L=10$.
The admissible leaders trajectories are defined as in Equation \eqref{eq:gototarget_beta}, we choose $\beta = 1$ for every leaders.
The target position is $\Xi_k(t)= x_1^\tau \ \forall t$ for the leaders moving towards the exit in the first room. While for the others is $\Xi_k(t)= x_2^\tau$ for $t>t_*$ and for $t<t_*$ we select two intermediate points in such a way that first leaders reach the staircases and then the second room.
Indeed, to evacuate, agents must either reach the exit in the first room, called exit $x_1^\tau$, or move towards the staircase, reach the second room and then search for the other exit, called exit $x_2^\tau$. The initial configuration is shown in Figure \ref{fig:test3_stair_0}. 

{\em Microscopic case. } Consider $N^L=10$ leaders. Assume that two leaders are aware of their role while the remaining are selfish leaders. Exits for every leader are chosen at time $t=0$ in such a way that five unaware leaders move toward exit $x_1^\tau$ and the remaining toward exit $x_2^\tau$. Among them, one of the two aware leaders moves towards one exit and the other towards the other exit. 

\begin{figure}[h!]
	\centering
	\includegraphics[width=0.45\linewidth]{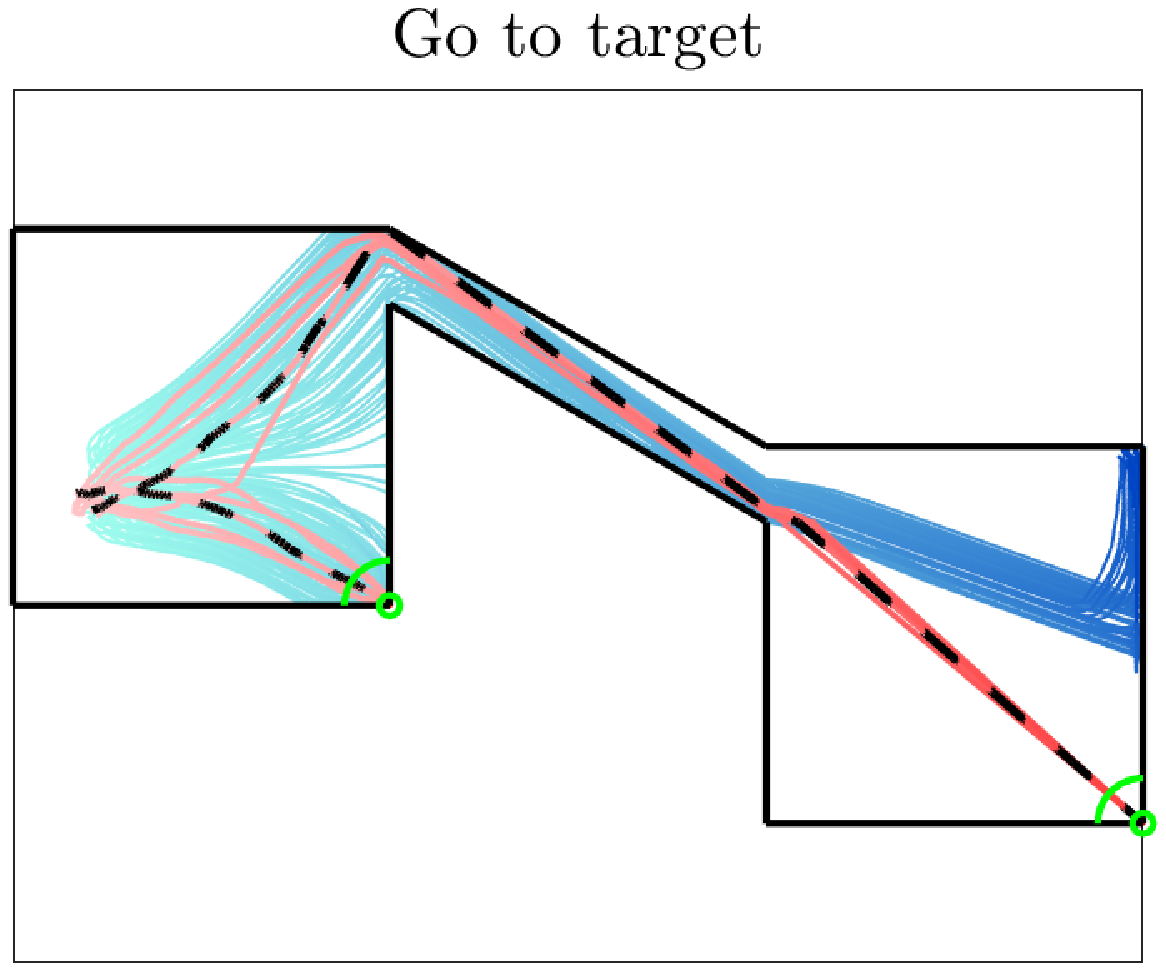}
	\includegraphics[width=0.45\linewidth]{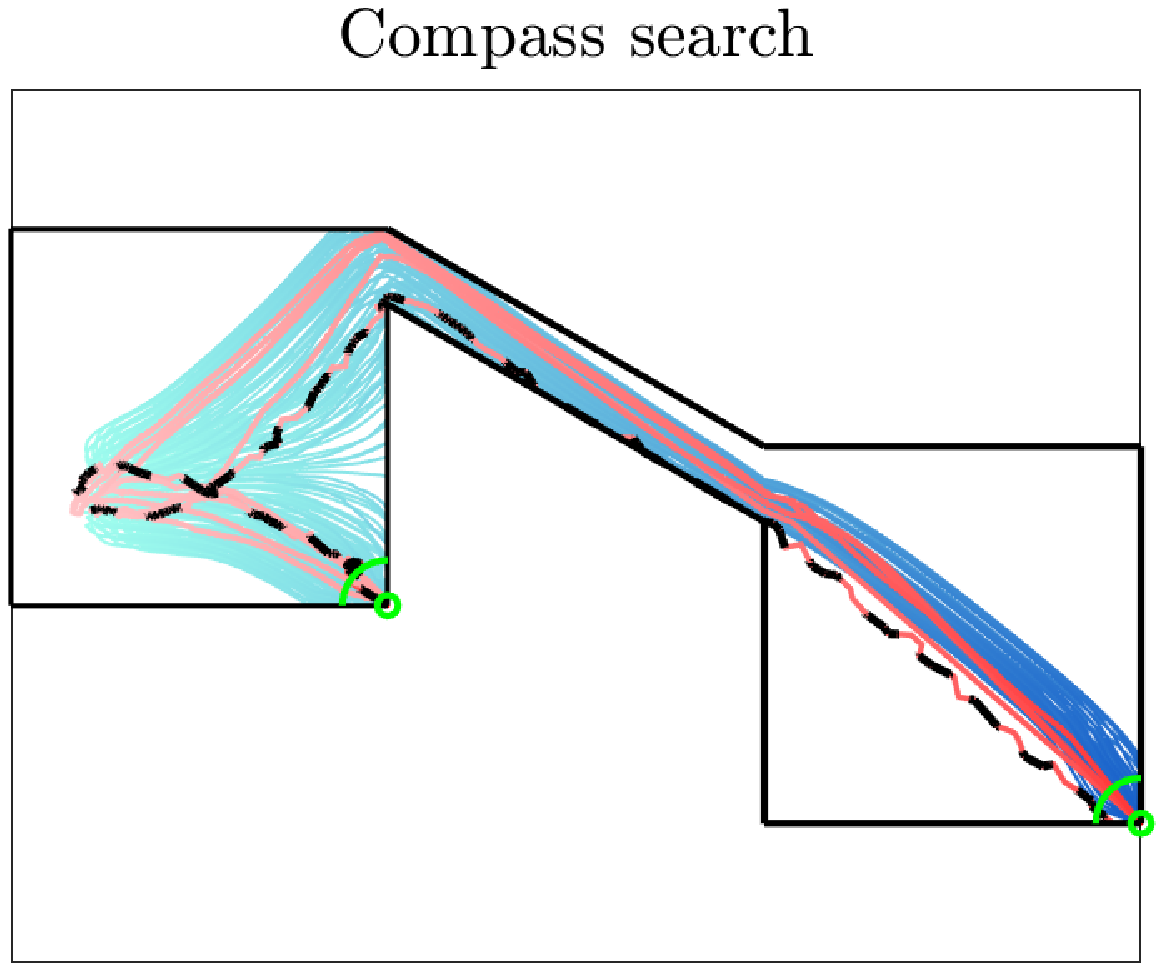}
	\caption{{\em Test 3b.} Microscopic case: mass splitting in presence of staircases. On the left, go-to-target. On the right, compass search. }
	\label{fig:test3_stairs_micro_1} 
\end{figure}

In the case of go-to-target strategy, leaders drive some followers to exit $x_1^\tau$ and some others to the staircase. The ones that reach the staircase move from the upper to the lower room and then are driven by leaders to exit $x_2^\tau$. As shown in Figure \ref{fig:test3_stairs_micro_1}  on the left, some followers are able to reach exit $x_1^\tau$ and some others to reach the second room. However, since the vast majority of leaders are unaware and move selfishly towards the exits, followers do not evacuate completely. Hence the only exit useful for evacuation is the one placed in the first room, exit $x_1^\tau$, whose visibility area is overcrowded.
On the right of Figure \ref{fig:test3_stairs_micro_1} leaders movement follows an optimized strategy allowing followers to split between the two exits. In this case, complete evacuation is achieved. Table \ref{eq:test3_table3} reports the performances of the two strategies. With the go-to-target strategy, all the evacuated followers reach the visibility area $\Sigma_1$ and hence are evacuated from exit $x_1^\tau$. With an optimized strategy instead, a larger amount of followers is evacuated and the overcrowding of the visibility areas is reduced. 

\begin{figure}[h!]
	\centering
	\includegraphics[width=0.328\linewidth]{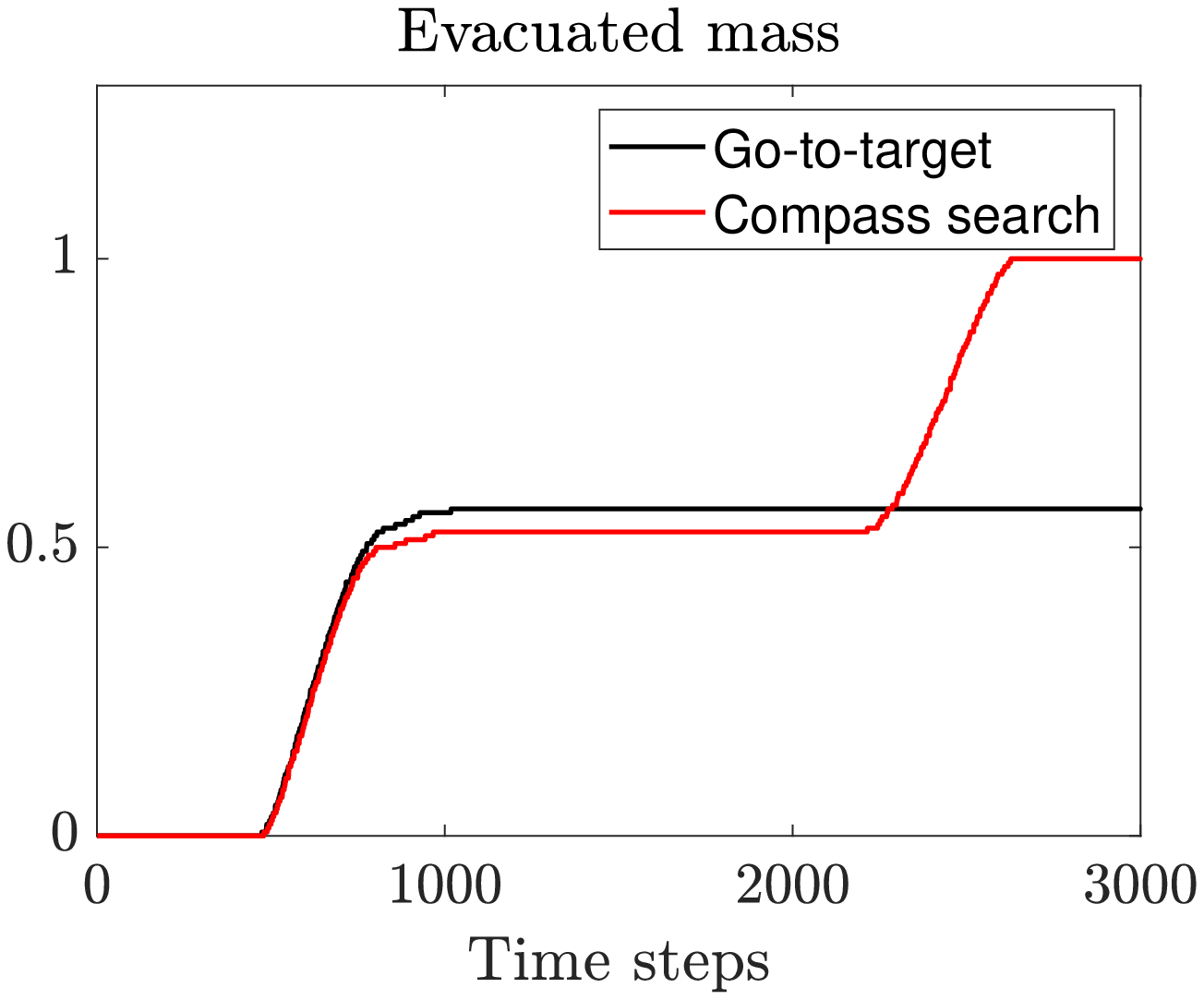}
	\includegraphics[width=0.328\linewidth]{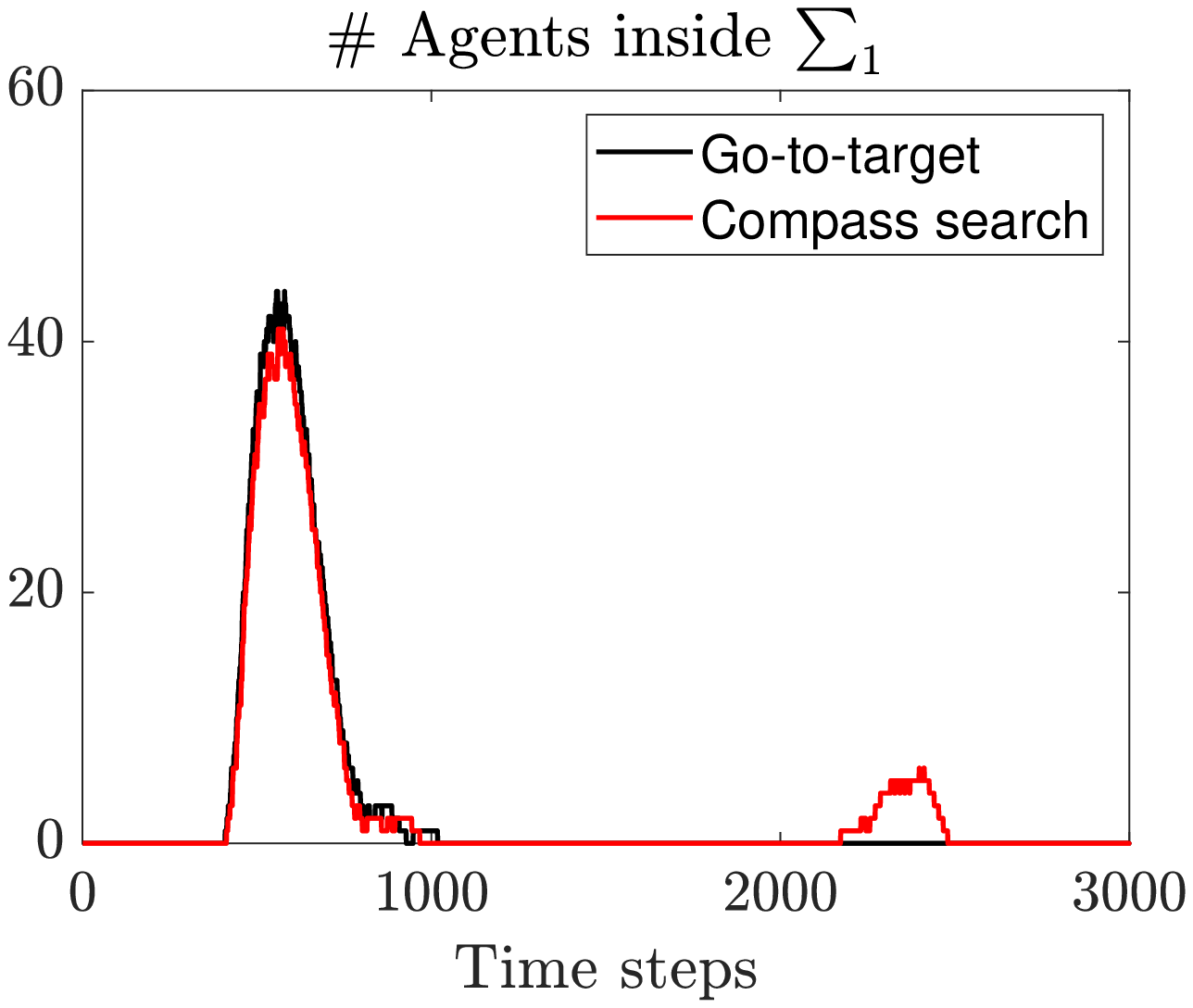}
	\includegraphics[width=0.328\linewidth]{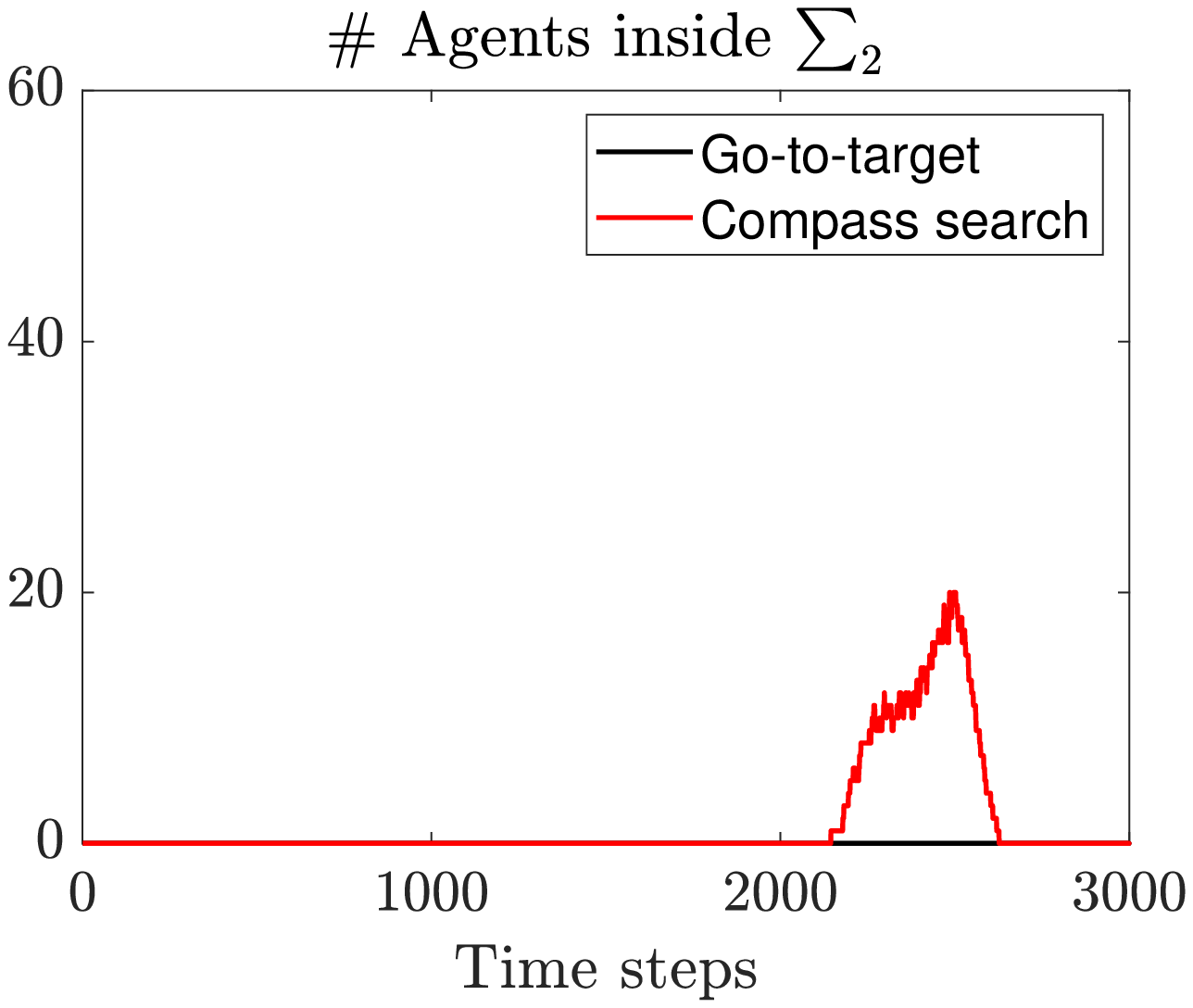}
	\caption{{\em Test 3b.} Microscopic case: mass splitting in presence of obstacles. Evacuated mass (left), occupancy of the visibility area $\Sigma_1$ (centre) and $\Sigma_2$ (right) as a function of time for go-to-target and compass search strategies.}
	\label{fig:test3_stairs_micro_2}
\end{figure}

\begin{table}\caption{{\em Test 3b}. Performances of mass splitting in the microscopic case.}\label{eq:test3_table3} 
	\centering
	\begin{tabular}{ c c c } 
		& go-to-target & CS (50 it)\\
		\cmidrule(r){2-2}\cmidrule(r){3-3}
		Evacuation time (time steps)  & > 3000 & 2627 \\
		Mass evacuated from $x_1^\tau$ & $57\%$ & $62\%$\\
		Mass evacuated  from $x_2^\tau$ & $0\%$ & $38\%$\\
		Total mass evacuated & $57\%$ & $100\%$ \\
		\hline
	\end{tabular}
\end{table}
In Figure \ref{fig:test3_stairs_micro_2}  we compare the evacuated mass and the occupancy of the exits visibility zone as a function of time for go-to-target strategy and optimal compass search strategy. Note that, with the compass search strategy, the whole mass is split between the two exits while with the go-to-target strategy the evacuated mass reaches only exit $x_1^\tau$.

\begin{figure}[h!]
	\centering
	\includegraphics[width=0.328\linewidth]{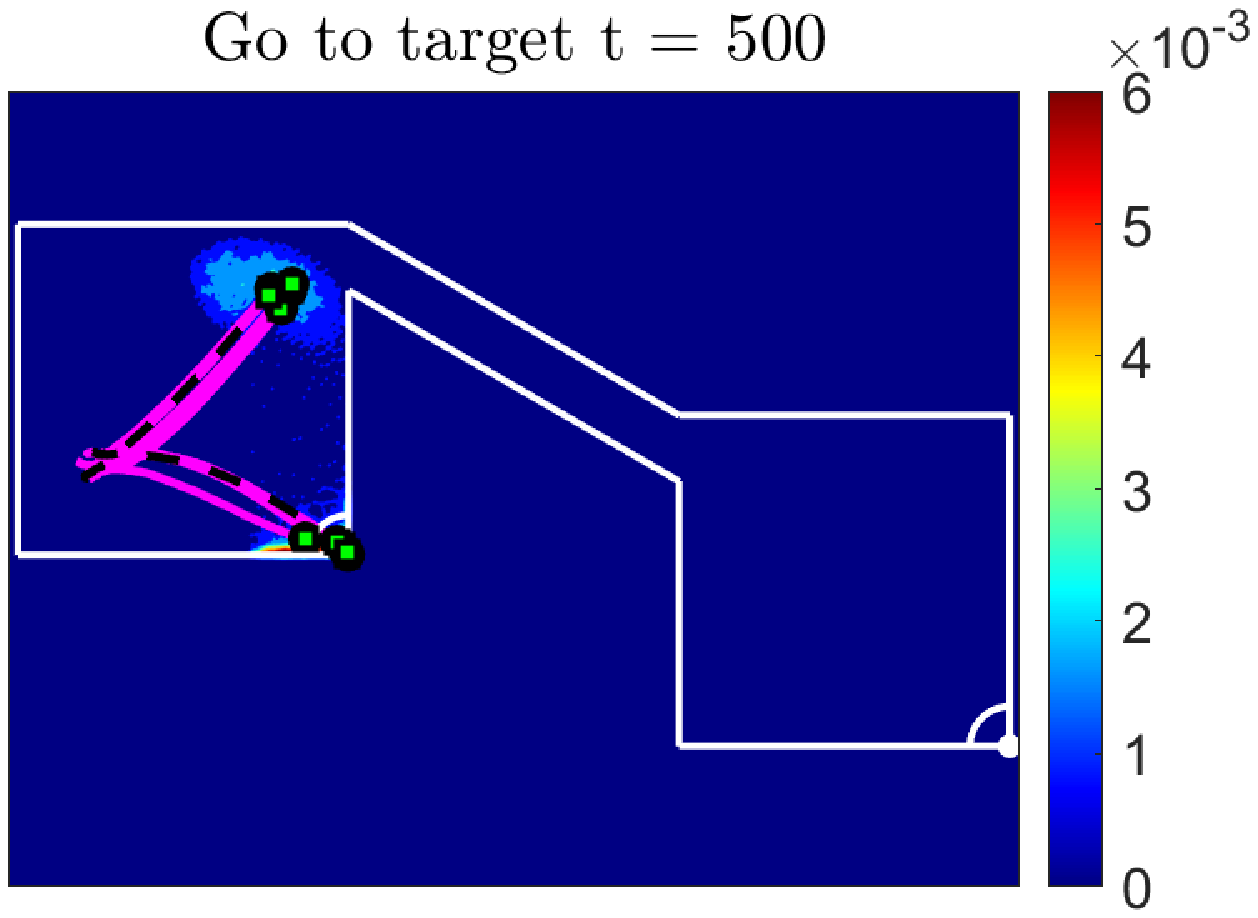}
	\includegraphics[width=0.328\linewidth]{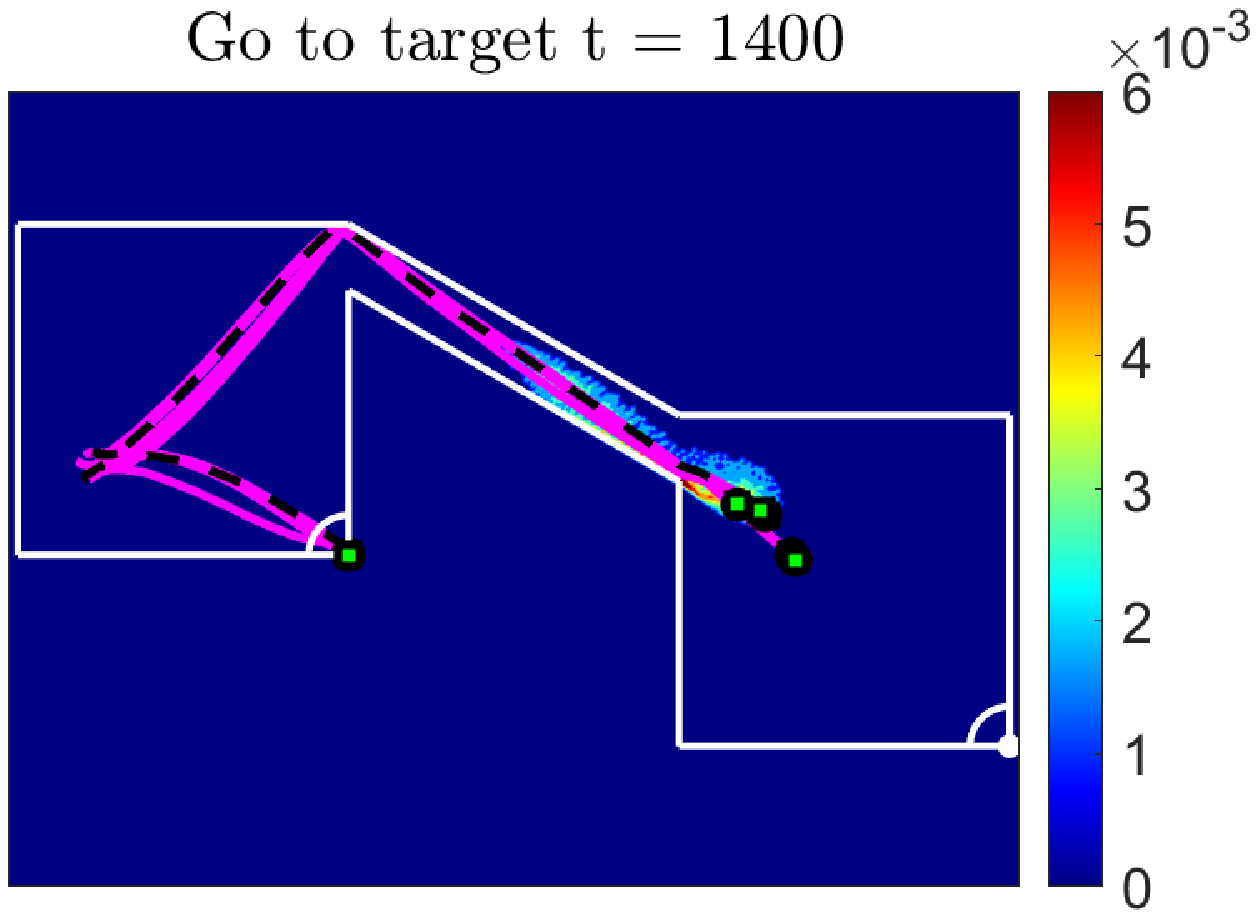}
	\includegraphics[width=0.328\linewidth]{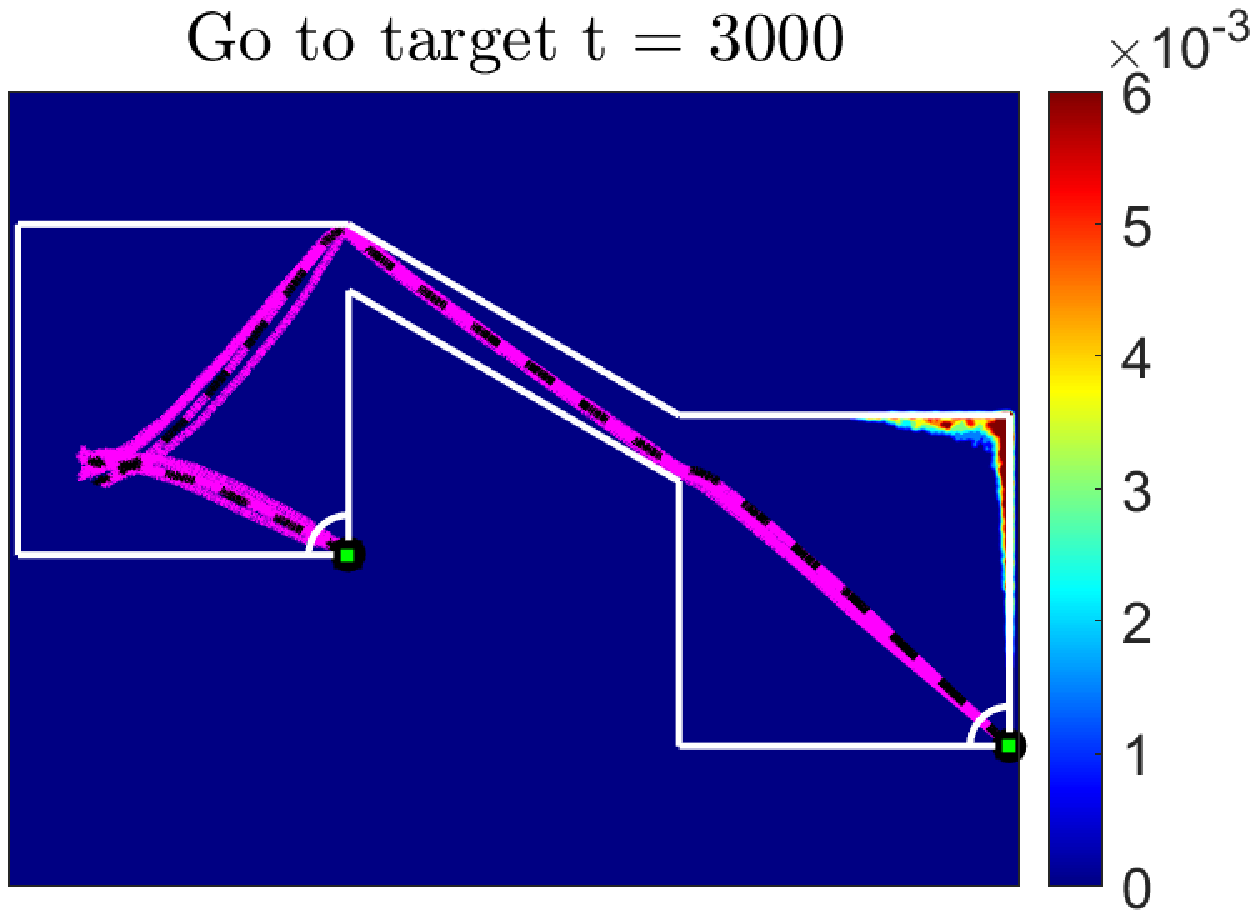}
	\\
	\includegraphics[width=0.328\linewidth]{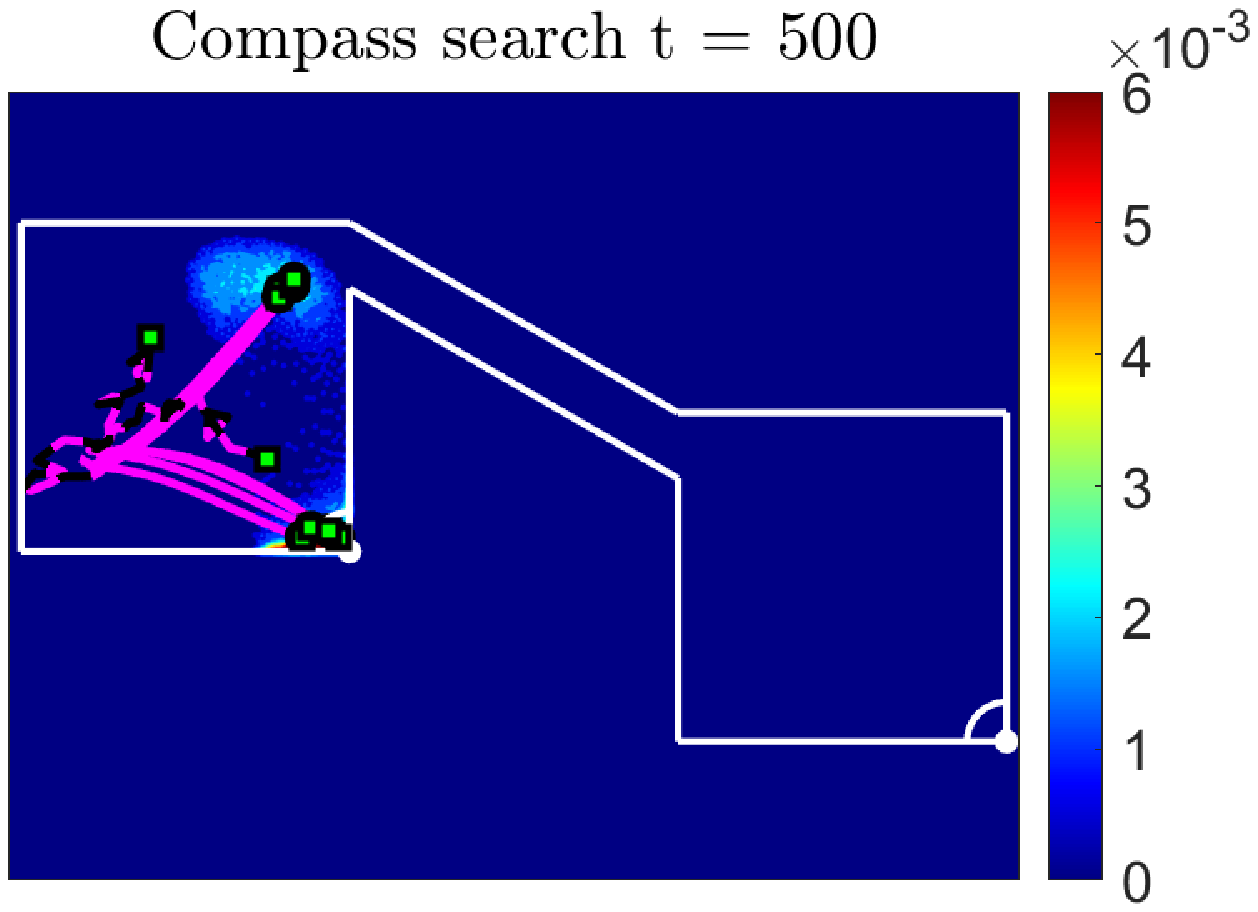}
	\includegraphics[width=0.328\linewidth]{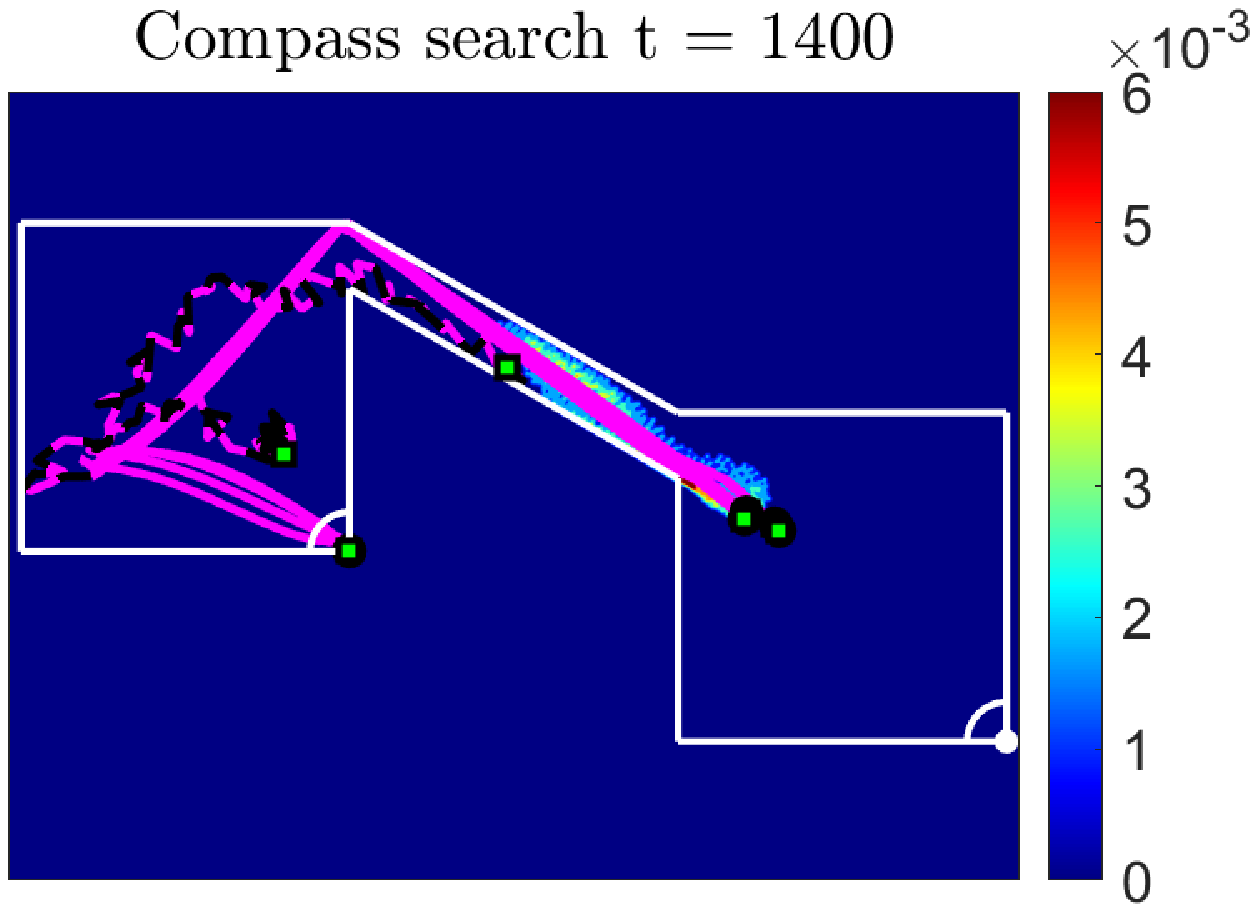}
	\includegraphics[width=0.328\linewidth]{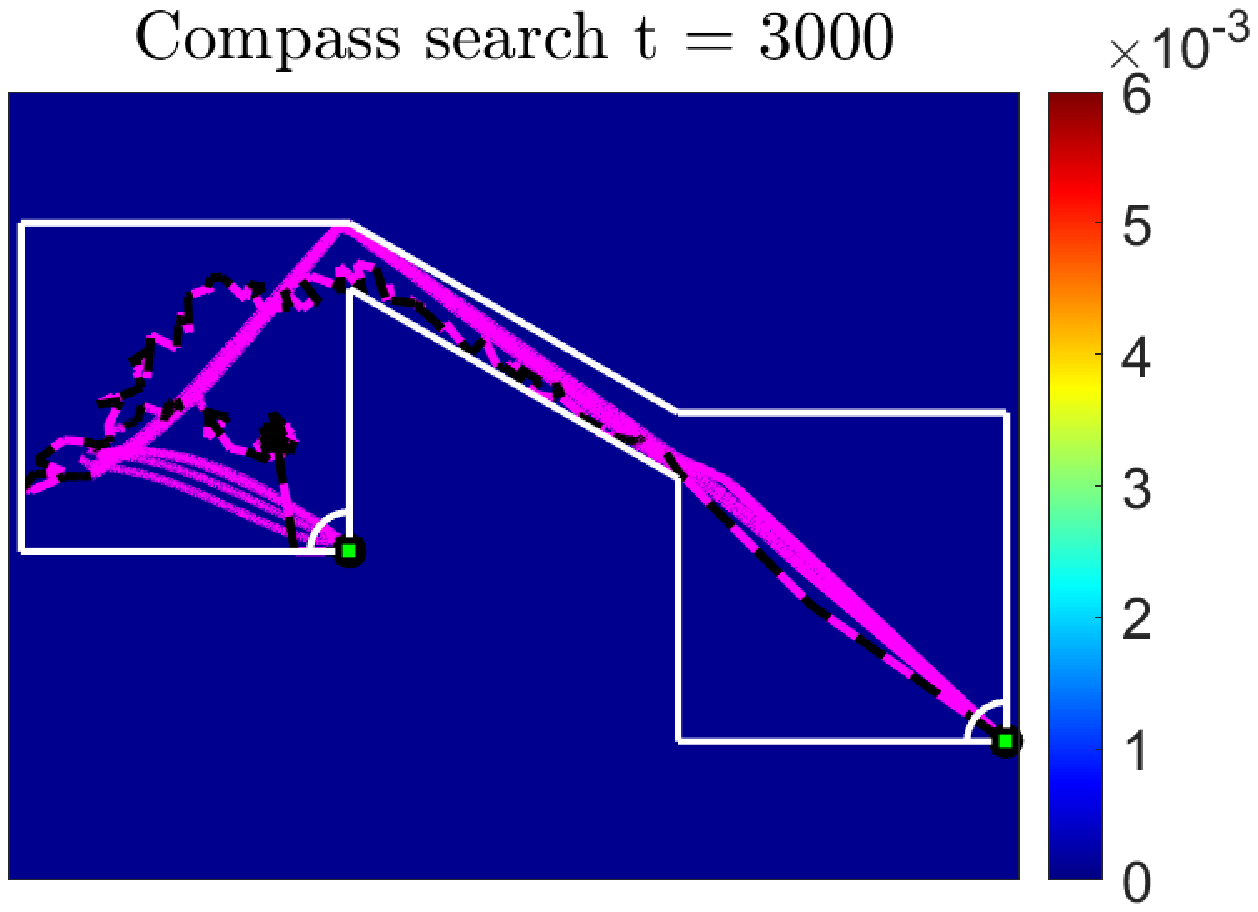}	
	\caption{{\em Test 3b.} Mesoscopic case: mass splitting in presence of staircases. Upper row: three snapshots taken at time $t=500$, $t=1400$, $t=3000$ with the go-to-target strategy. Lower row: three snapshots taken at time $t=500$, $t=1400$, $t=3000$ with the optimized compass search strategy.}
	\label{fig:test3_stairs_meso_1}
\end{figure}

{\em Mesoscopic case. }  Consider the case of a continuous mass of followers. Similar to the microscopic case we observe in Figure \ref{fig:test3_stairs_meso_1} the evolution of the dynamics with fixed strategy and with the optimized one. The upper row shows that with the go-to-target strategy the total evacuated mass is not split between the two exits since just the $1,2\%$ of mass reaches exit $x_2^\tau$. However, as shown in Table \ref{eq:test3_table3} a larger percentage of followers reaches exit $x_1^\tau$ and the remaining part spreads in the second room without evacuate. 

The lower row of Figure  \ref{fig:test3_stairs_meso_1}  shows the dynamics obtained with the optimized compass search strategy. At the time $t=500$ a larger follower mass is moving towards the staircase. At time $t=3000$ almost all the mass is evacuated and split between the two exits. In Table \ref{eq:test3_table3} we compare the two strategies showing that with the compass search technique it is possible to improve the mass splitting.

 \begin{table}\caption{{\em Test 3b}. Performances of mass splitting in the mesoscopic case.}\label{eq:test3_table3} 
	\centering
	\begin{tabular}{ c c c } 
		& go-to-target & CS (50 it)\\
		\cmidrule(r){2-2}\cmidrule(r){3-3}
		Evacuation time (time steps)  & > 3000 & > 3000 \\
		Mass evacuated from $x_1^\tau$ & $46,6\%$ & $51\%$\\
		Mass evacuated  from $x_2^\tau$ & $1,2\%$ & $48\%$\\
		Total mass evacuated & $47,8\%$ & $99\%$ \\
		\hline
	\end{tabular}
\end{table}

\begin{figure}[h!]
	\centering
	\includegraphics[width=0.328\linewidth]{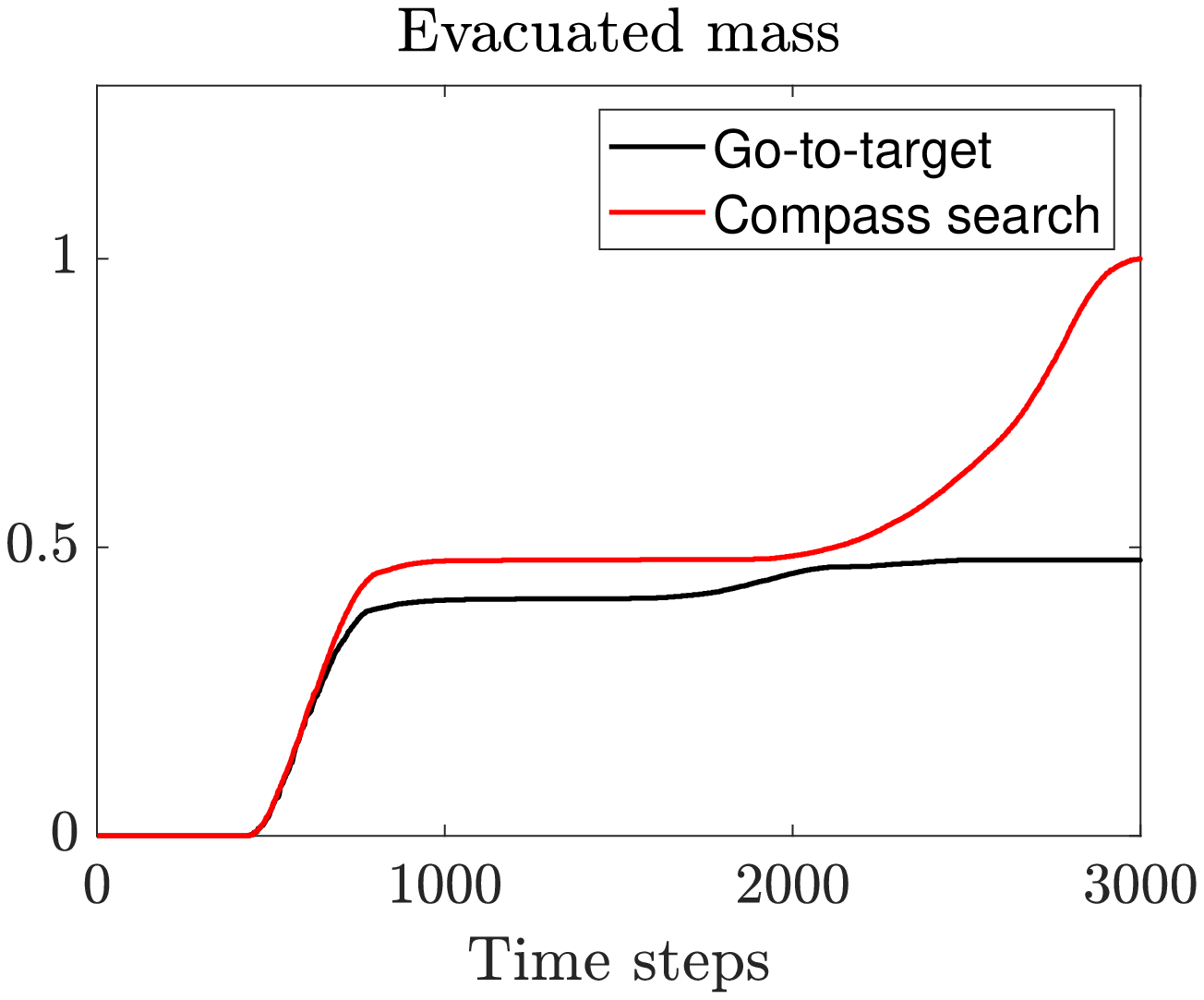}
	\includegraphics[width=0.328\linewidth]{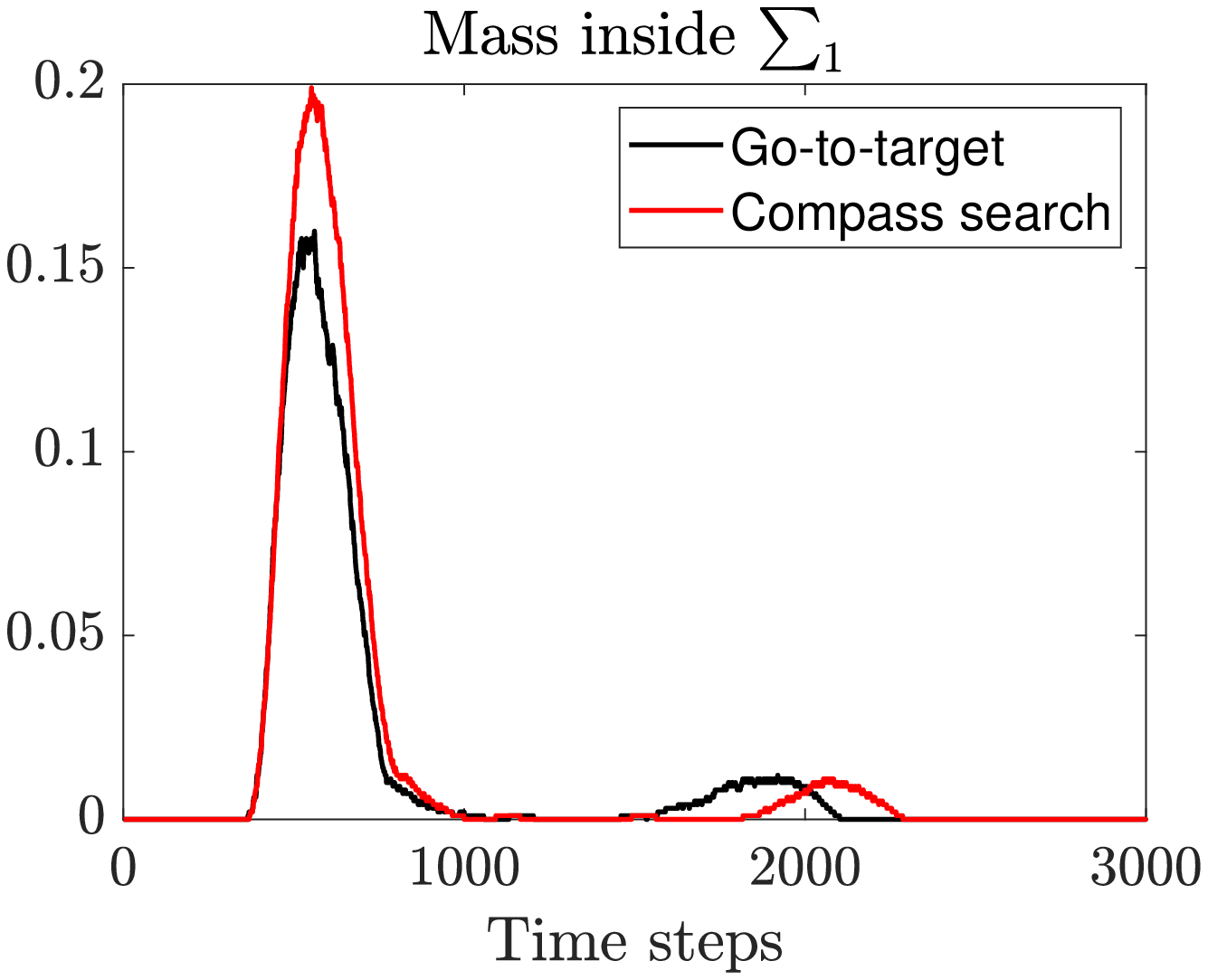}
	\includegraphics[width=0.328\linewidth]{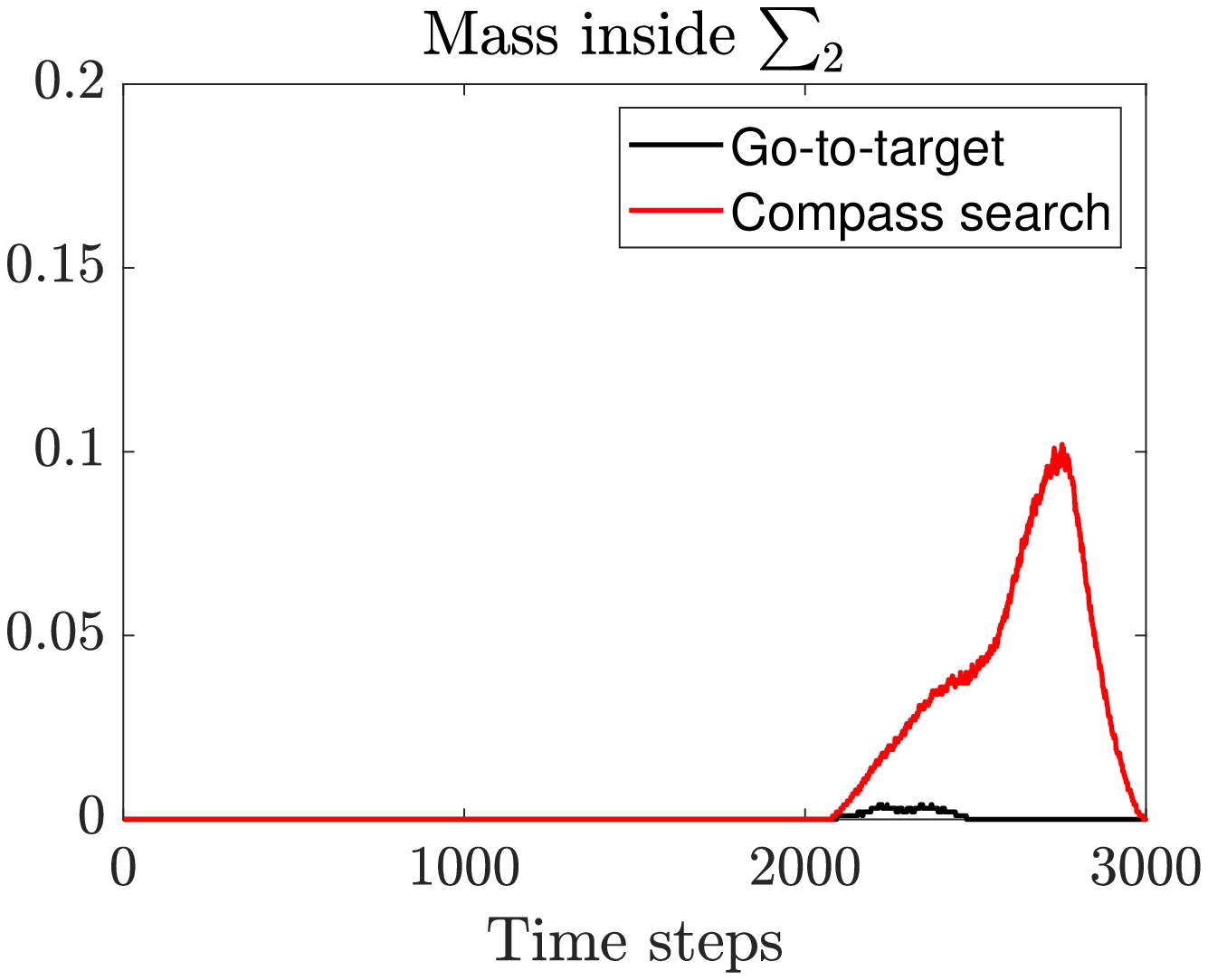}
	\caption{{\em Test 3b.} Mesoscopic case: mass splitting in presence of obstacles. Evacuated mass (left), occupancy of the visibility area $\Sigma_1$ (centre) and $\Sigma_2$ (right) as a function of time for go-to-target and optimal compass search strategies.}
	\label{fig:test3_stairs_meso_3}
\end{figure}

In Figure \ref{fig:test3_stairs_meso_3}  we compare the evacuated mass and the occupancy of the exits visibility zone as a function
of time for go-to-target strategy and optimal compass search strategy.  Note that, with the compass search strategy a larger percentage of mass reaches exit $x_2^\tau$ than with the go-to-target strategy.

\subsection{Discussion and comparison}\label{sec:discussion}
 In the previous tests we have considered different scenarios to create more complex situations in relation to the functionals chosen, \cite{zhou2018optimal}. In general, given a ceratin setting, it is difficult to choose the optimal number of leaders that guarantee evacuation, and a high number of leaders does not necessarily imply better evacuation efficiency, see for example \cite{ma2016effective}. 
 Another challenging aspect is to give an uniform measure of the performance of the different strategies in such different contexts. A viable option is to quantify the congestion around the exits to exclude dangerous situations. 
 Following the idea in \cite{feliciani2018measurement} we consider the congestion value 
 \[
 cong_{\Sigma_i}(t) = \rho_{\Sigma_i}(t) var_{\Sigma_i}(v(t))
 \] 
 where $\rho_{\Sigma_i}(t)$ is the number of agents (mass) in the microscopic (mesoscopic) case inside $\Sigma_i$ at time $t$ and \[
 var_{\Sigma_i}(v(t)) = \frac{1}{\rho_{\Sigma_i}(t)} \sum_{j\in \Sigma_i}\left( {\vert v_j(t) \vert -s}\right) ^2.
 \]   We consider also $m_{\Sigma_i}$ the maximum number of pedestrians over time inside the visibility area $\Sigma_i$ and $l_{\Sigma_i}$ the percentage of time in which the visibility area $\Sigma_i$ is not empty, finally we denote by $M_{\Sigma_i}$ the percentage of mass inside $\Sigma_i$ in the mesoscopic case.
 
 In this way we can compare the congestion  of the various exits for different settings, showing that the more desiderable situations are when $cong_{\Sigma_i}$ and $m_{\Sigma_i}$ ($M_{\Sigma_i}$) are small and $l_{\Sigma_i}$ is high.
 We reported in Table \ref{eq:congestion_micro} and Table \ref{eq:congestion_meso} respectively the values for the microscopic and the mesoscopic setting.

 \begin{table}[H]
 	\centering
 	\caption{Comparison of the congestion in the visibility areas for the microscopic case. In red the maximum value of $cong_{\Sigma_i}$ among the visibility areas $\Sigma_i$.} \label{eq:congestion_micro}
 	\begin{tabular}{c|ccccccccc}
 		&$cong_{\Sigma_1}$& $cong_{\Sigma_2}$&$cong_{\Sigma_3}$&$m_{\Sigma_1}$ &$m_{\Sigma_2}$&$m_{\Sigma_3}$	& $l_{\Sigma_1}$ &$l_{\Sigma_2}$&$l_{\Sigma_3}$  \\
 		\hline
 		{ Test 1a}&\textcolor{red}{0.039}&0.011&0.012 & 40 & 19&17& 0.73 & 0.51 & 0.33 \\
 		{  Test 1b} &\textcolor{red}{0.013}&0.009&-& 27 & 16 &-& 0.36 & 0.22 & -\\
 		{  Test 2} &0.009&\textcolor{red}{0.056}&-& 20 & 54 &-& 0.13 & 0.31 & -\\
 		{  Test 3a}&\textcolor{red}{0.035}&0.027&- & 43 & 26 &- & 0.19 & 0.29 & -\\
 		{  Test 3b}& \textcolor{red}{0.024}&0.006&-& 41 & 20&- & 0.28 & 0.16 & -\\
 		\hline
 	\end{tabular}
 \end{table}

\begin{figure}[h!]
 	\centering
 	\includegraphics[width=0.328\linewidth]{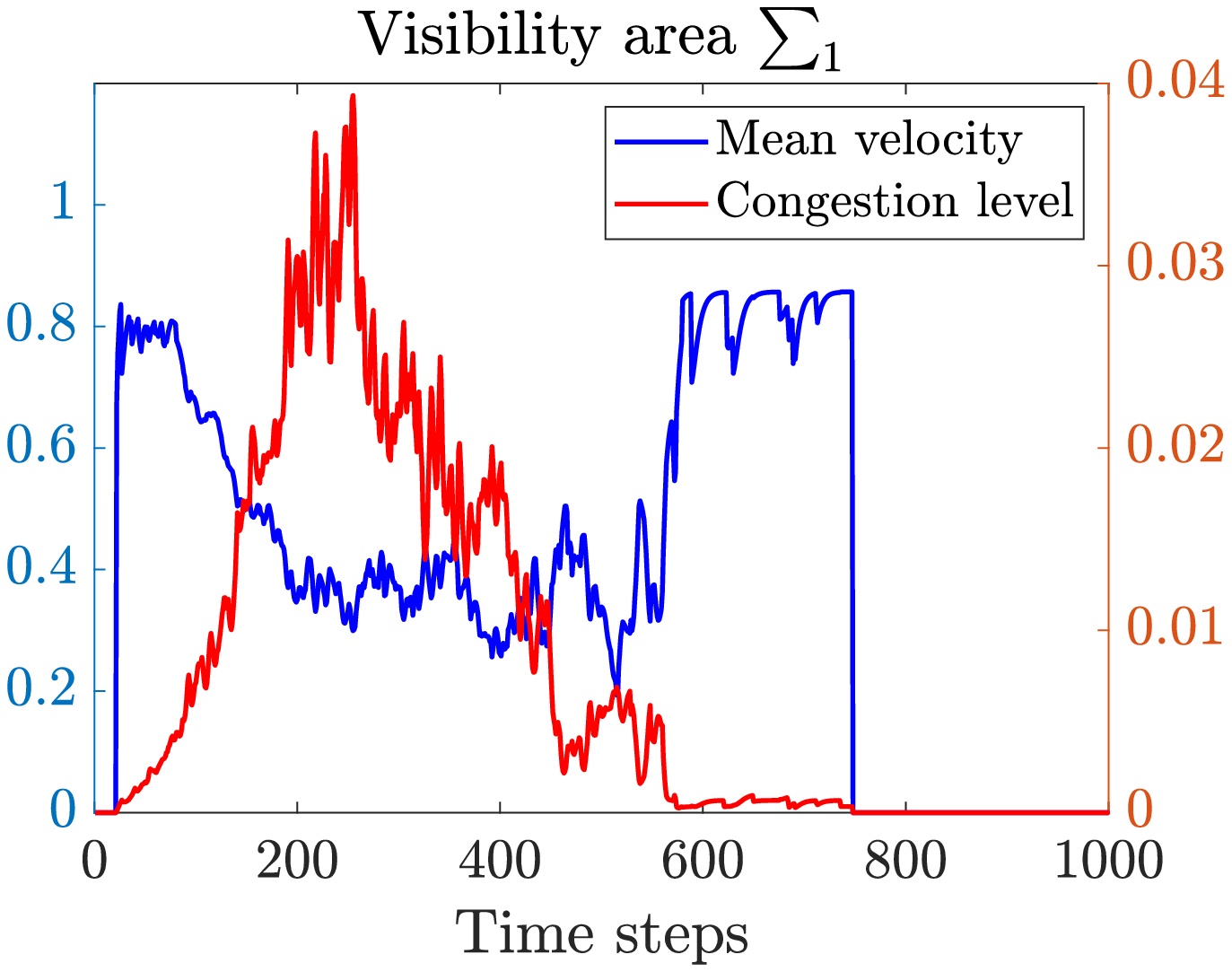}
 	\includegraphics[width=0.328\linewidth]{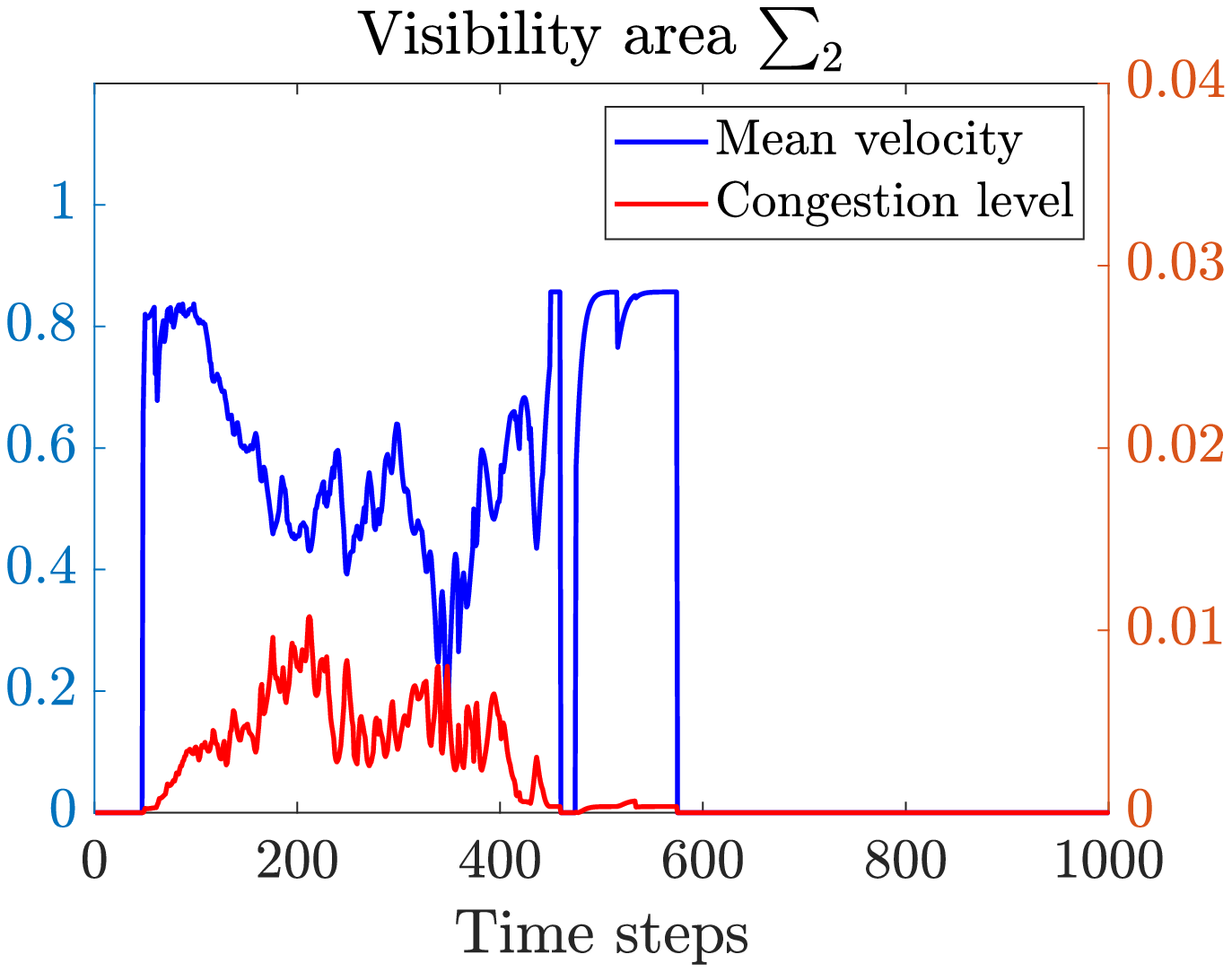}
 	\includegraphics[width=0.328\linewidth]{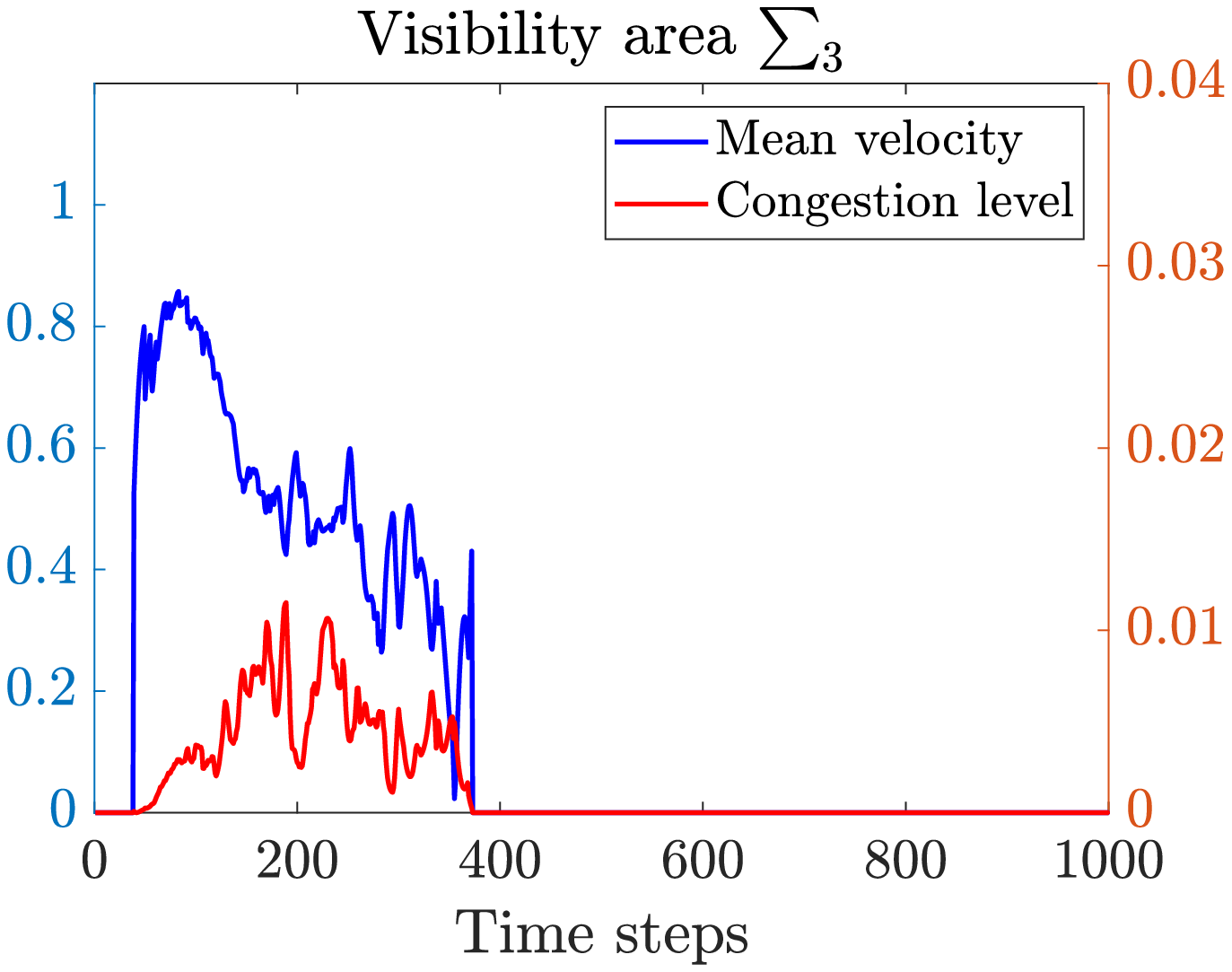}
 	\caption{{  Test 1a}. Microscopic case: number of agents and mean velocity of the visibility areas. }
 	\label{fig:test1_vel_micro}
 \end{figure}

\begin{figure}[h!]
 	\centering
 	\includegraphics[width=0.328\linewidth]{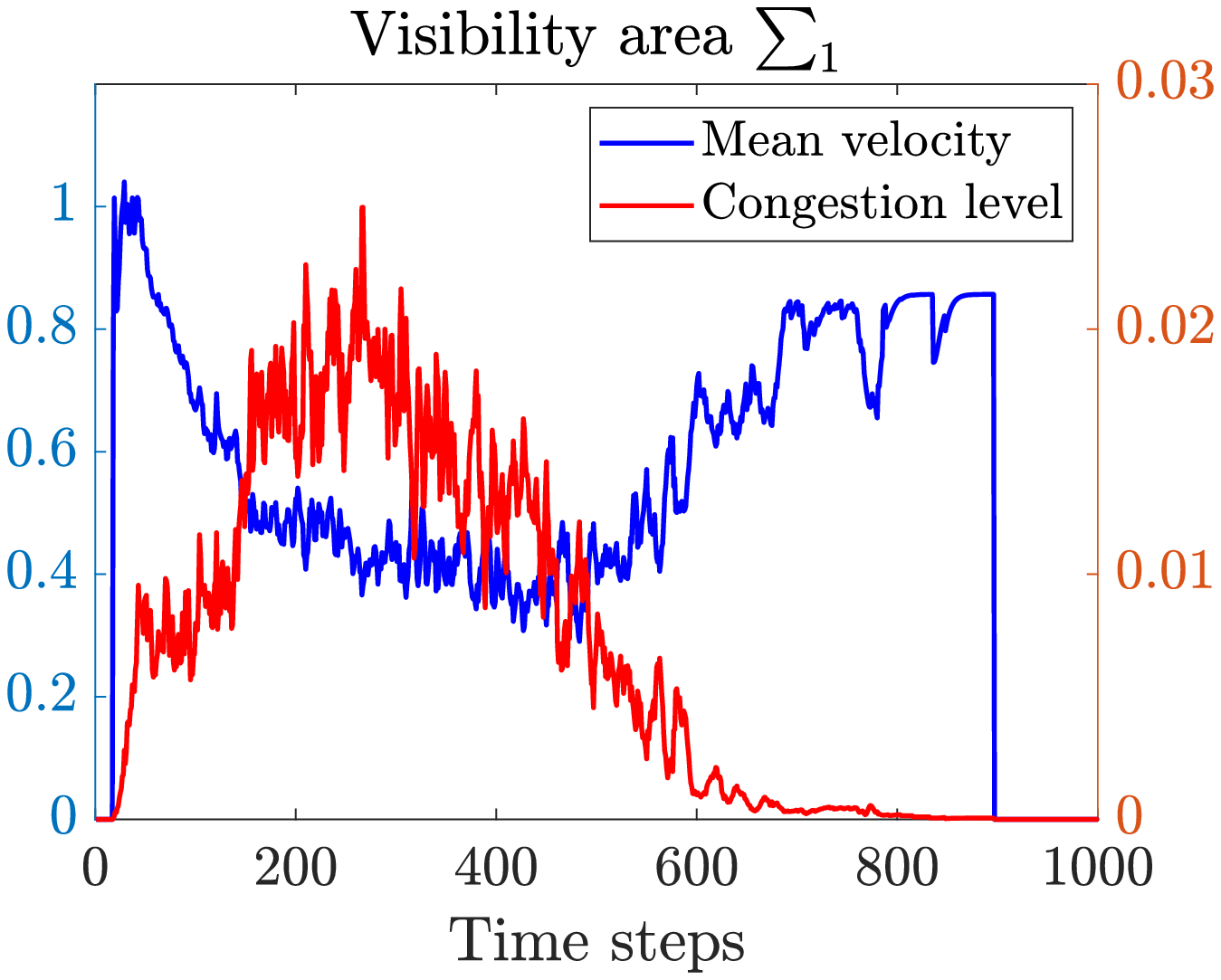}
 	\includegraphics[width=0.328\linewidth]{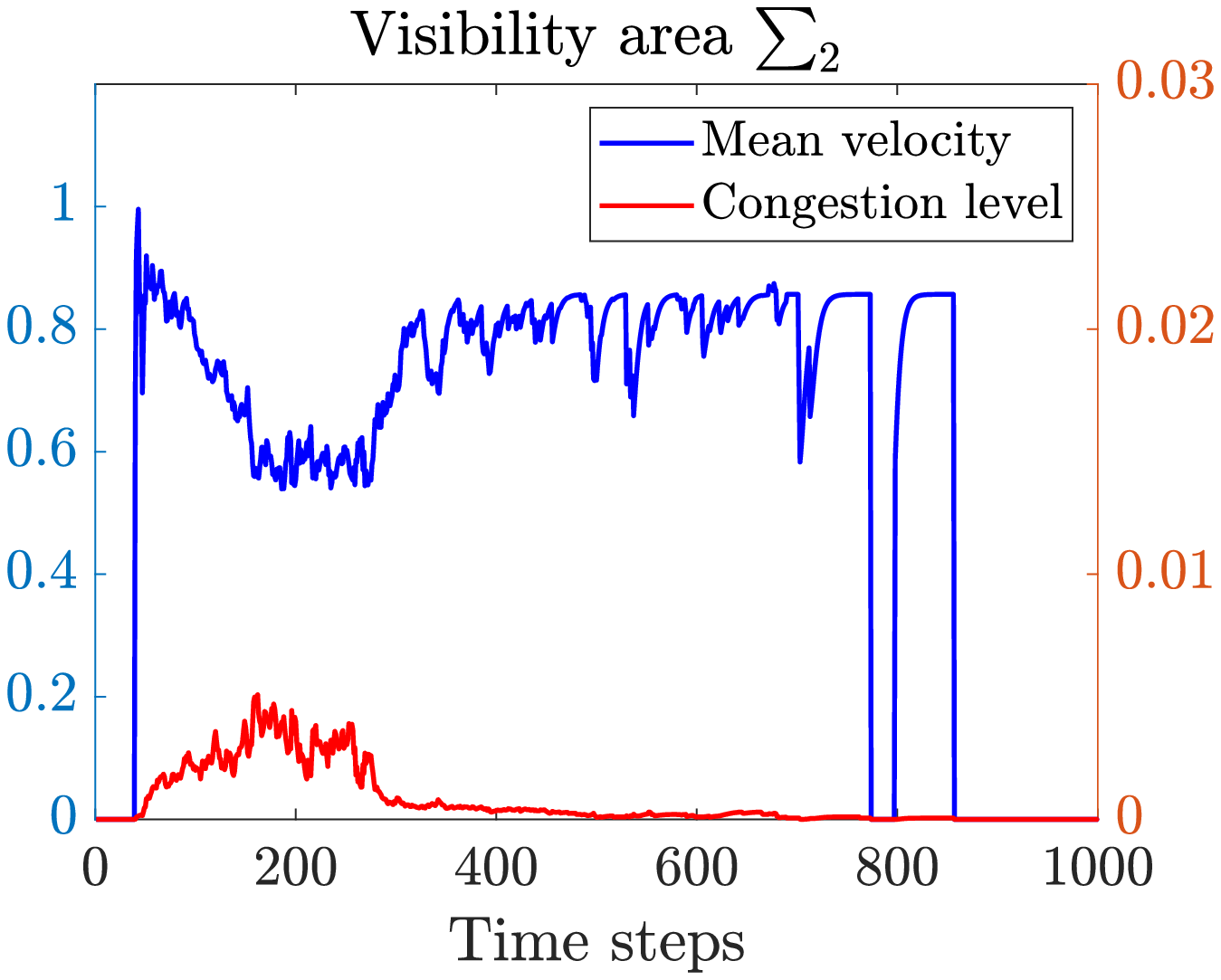}
 	\includegraphics[width=0.328\linewidth]{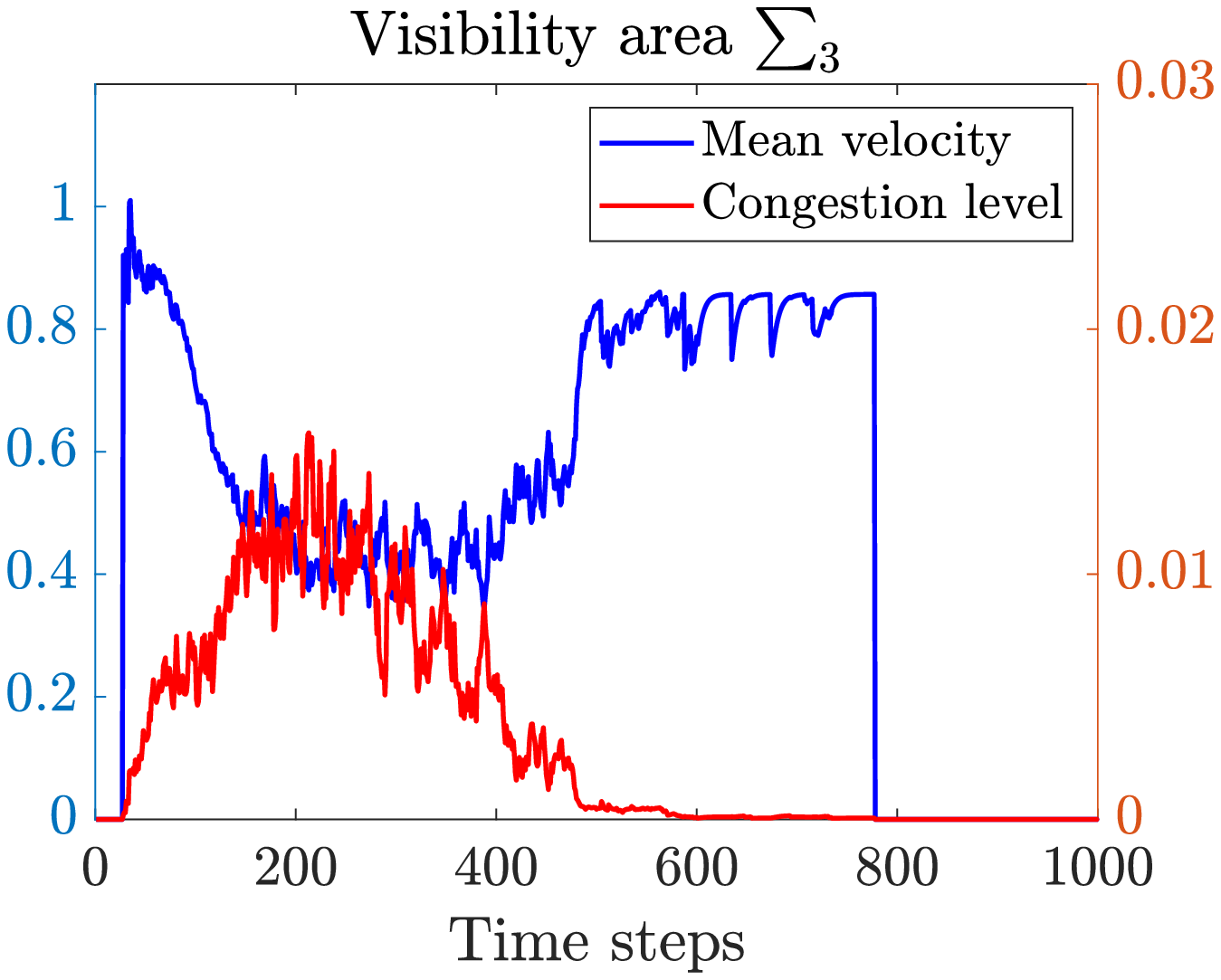}
 	\caption{{  Test 1a}. Mesoscopic case: mass of agents and mean velocity of the visibility areas. }
 	\label{fig:test1_vel_meso}
 \end{figure}

 \begin{table}[H]
 	\centering
 	\caption{Comparison of the congestion in the visibility areas for the mesoscopic case. In red the maximum value of $cong_{\Sigma_i}$ among the visibility areas $\Sigma_i$.} \label{eq:congestion_meso}
  	\begin{tabular}{c|ccccccccc}
  	&$cong_{\Sigma_1}$& $cong_{\Sigma_2}$&$cong_{\Sigma_3}$&$M_{\Sigma_1}$ &$M_{\Sigma_2}$&$M_{\Sigma_3}$	& $l_{\Sigma_1}$ &$l_{\Sigma_2}$&$l_{\Sigma_3}$  \\
  	\hline
  	{  Test 1a}&\textcolor{red}{0.025}&0.005&0.016 & 0.22 & 0.6&0.16& 0.88 & 0.79 & 0.75 \\
  	{  Test 1b} &\textcolor{red}{0.010}&0.005&-& 0.1 & 0.08&-& 0.51 & 0.26 & -\\
  	{  Test 2} &0&\textcolor{red}{0.009}&-& 0 & 0.12 &-& 0 & 0.36 & -\\
  	{  Test 3a}&0.005&\textcolor{red}{0.011}&- & 0.07& 0.12 &- & 0.3 & 0.32 & -\\
  	{  Test 3b}& \textcolor{red}{0.013}&0.004&-& 0.2& 0.1&- & 0.41 & 0.3 & -\\
  	\hline
  \end{tabular}
\end{table}
Finally, Figures \ref{fig:test1_vel_micro}-\ref{fig:test1_vel_meso} show the mean velocity and the congestion level for the case of evacuation with three exits ({  Test 1a}) in the microscopic and mesoscopic case respectively. These plots underline that if the congestion level is higher then the mean velocity is lower.


\section{Conclusions}\label{sec:conclusions}
This work has been devoted to the study of optimized strategies for the control of egressing pedestrians in an unknown environment. In particular, we studied situations with complex environments where multiple exits and obstacles are present. Few informed agents act as controllers over the crowd, without being recognized as such. Indeed it has been shown that minimal intervention can change completely the behavior of a large crowd, and at the same time avoiding adversarial behaviors. 
On the other hand, we observed that if part of the informed agents moves without coordinated action, this may cause critical situations, such as congestion around the exit. Hence it is important to have a clear understanding of different strategies to enhance the safe evacuation of the crowd. To this end, we explored various optimization tasks such as minimum time evacuation, maximization of mass evacuated, and optimal mass distribution among exits.

We investigated these dynamics at the various scales: from the microscopic scale of agent-based systems to the statistical description of the system given by mesoscopic scale. Numerically we proposed an efficient scheme for the simulation of the mean-field dynamics, whereas we use a meta-heuristic approach for the synthesis of optimized leaders strategies. The proposed numerical experiments suggest that the optimization of leaders movements is enough to de-escalate critical situations.

Different questions arise at the level of control through leaders with multiple exits and obstacles. In such a rich environment several research directions can be explored, such as optimal positioning and amount of leaders within the crowd, or different type of cooperative strategies among different groups of leaders to optimally distribute the followers crowd.

\bibliographystyle{spmpsci} 
\bibliography{biblio_crowd}
\end{document}